\documentclass{article}
\usepackage[utf8]{inputenc}
\usepackage[english]{babel}
\usepackage{amsmath}
\usepackage{amsfonts}
\usepackage{amssymb}
\usepackage{amsthm,mathrsfs}
\usepackage{graphicx,xcolor,tikz}
\usepackage{caption,subcaption,multicol,longtable}
\usepackage[left=1.5in, right=1.5in, top=1.75in, bottom=1.5in]{geometry}
\usepackage{enumerate}
\usepackage[hidelinks]{hyperref}
\usepackage[square,numbers]{natbib}
\usepackage[hidelinks]{hyperref}
\usepackage{mathtools} 
\usepackage{cutwin}

\theoremstyle{plain}
\newtheorem{thm}{Theorem}[section]

\newtheorem{question}{Question}[section]
\newtheorem{lem}[thm]{Lemma}
\newtheorem{cor}[thm]{Corollary}
\newtheorem{prop}[thm]{Proposition}
\newtheorem{rem}[thm]{Remark}
\newtheorem{algorithm}[thm]{Algorithm}

\theoremstyle{definition}

\def \N {{\mathbb N}}

\def \sig {{\Sigma}}
\def \s {{\Sigma_{1}}}
\def \t {{\Sigma_{2}}}
\def \mH {{\mathcal{H}}}
\def \qsig {{Q_{\sig}}}

\graphicspath{{figures/}}

\title{Reducing spheres of genus-2 Heegaard splitting of $S^3$}
\author{Sreekrishna Palaparthi$^a$ , Swapnendu Panda$^b$}
\begin{document}
\maketitle

\begin{abstract}
The Goeritz group of the standard genus-g Heegaard splitting of the three sphere, $G_g$, acts on the space of isotopy classes of reducing spheres for this Heegaard splitting. Scharlemann \citep{scharlemann} uses this action to prove that $G_2$ is finitely generated. In this article, we give an algorithm to construct any reducing sphere from a standard reducing sphere for a genus-2 Heegaard splitting of the $S^3$. Using this we give an alternate proof of the finite generation of $G_2$ assuming the finite generation of the stabilizer of the standard reducing sphere.
\end{abstract}

\noindent
\small{$^a${\sc Sreekrishna Palaparthi}, Department of Mathematics, Indian Institute of Technology Guwahati, Assam 781039, India, email : passkrishna@iitg.ac.in\\
\noindent
$^b${\sc Swapnendu Panda}, Department of Mathematics, Indian Institute of Technology Guwahati, Assam 781039, India, email : p.swapnendu@iitg.ac.in

\noindent
\textbf{Acknowledgements:} We are thankful to the Department of Mathematics, Indian Institute of Technology Guwahati for its support while carrying out this work. A part of this work intersects with the thesis work of the second author. \\

\section{Introduction}\label{section:intro}
	The genus-$g$ Goeritz group of $S^3$, denoted $G_g$, is the group of isotopy classes of homeomorphisms of $S^3$ which preserve a genus-$g$ Heegaard surface of $S^3$. By a theorem of Waldhausen \citep{waldhausen}, every Heegaard splitting of $S^3$ of genus greater than $1$ is stabilized. So the group $G_g$ can be taken to be the group of isotopy classes of homeomorphisms of $S^3$ which preserve the standard genus-$g$ Heegaard splitting of $S^3$. Scharlemann \citep{scharlemann} proved that the genus-2 Goeritz group, $G_2$, of $S^3$ is finitely generated by four generators. He did so by considering a certain simplicial complex whose vertices are the isotopy classes of reducing spheres of the standard genus-$2$ Heegaard splitting of $S^3$. Akbas \citep{akbas} used Scharlemann's complex to give a finite presentation for $G_2$. Cho \citep{cho} independently used a certain disk complex to give a finite presentation for $G_2$. For $g \geq 3$, Zupan \citep{zupan} generalized Scharlemann's work and defined a reducing sphere complex for genus-$g$ Heegaard splittings of $S^3$ and showed that the group $G_g$ is finitely generated if and only if this complex is connected. Freedman and Scharlemann \citep{freedman} proved that $G_3$ is also finitely generated. They conjectured that the five generators mentioned in Powell \citep{powell} generate $G_g$ for all $g \geq 3$. Much is unknown about $G_g$ for $g \geq 4$.
	
	The literature on the Goeritz group of $S^3$ shows the importance of describing an arbitrary reducing sphere for reducible Heegaard splittings of $S^3$, and for reducible Heegaard splittings of 3-manifolds in general. In this article, we look at the set of all the reducing spheres for the genus-2 Heegaard splitting of $S^3$. We give structure theorems for such spheres. We then give an algorithm to construct any reducing sphere for the standard genus-2 Heegaard splitting of $S^3$, upto isotopy, from the standard reducing sphere.
	
	The layout of this article is as follows. In section \ref{section:red_curves}, we prove Theorem \ref{arc-lemma} which gives some intersection conditions satisfied by the curve of intersection of a reducing sphere with the standard genus-2 Heegaard surface of $S^3$. In section \ref{section:arcs_pi1}, we use the fundamental group of a genus-2 handlebody to deduce some conditions which apply to the arcs of intersection of a  reducing sphere with a certain component pair of pants of the genus-2 Heegaard surface of $S^3$. We then use these conditions to prove Theorem \ref{thm:aa-atmost-2} for such arcs of a reducing sphere. In section \ref{section:beta_lowering}, we prove Theorem \ref{thm:beta_single_letter} which is a key tool to de-construct a given reducing sphere. In section \ref{section:algorithm}, we give an algorithm to de-construct and transform any reducing sphere into the standard reducing sphere. Using this algorithm and assuming that the stabilizer of the standard reducing sphere of the genus-2 Heegaard splitting of $S^3$ is finitely generated, we show that the Goeritz group $G_2$ is finitely generated.

\section{Reducing Curves - Intersection Numbers}\label{section:red_curves}
	For preliminaries regarding Heegaard splittings of 3-manifolds, one can consult \citep{scharlemann_notes} and for preliminaries on  mapping class groups and simple closed curves on surfaces, one can consult Farb and Margalit \citep{farb_margalit}. We start by considering the standard genus-2 Heegaard splitting of $S^3$, written $S^3=V\cup_{\sig} W$ and denoted by $\mH$, where both $V$ and $W$ are unknotted genus-2 handlebodies embedded in $S^3$ with disjoint interiors glued along their common boundary $\sig$. The Heegaard surface, $\sig$, of the above splitting is a closed orientable surface of genus two embedded in $S^3$. Since $\mH$ is stabilized, one can find a reducing sphere for $\mH$ \textit{i.e.} a sphere which intersects each of $V$ and $W$ in an embedded disk with a common boundary circle on $\sig$. Such a reducing sphere intersects $\sig$ in an essential separating simple closed curve, which we call a reducing curve for $\sig$. A reducing curve on $\sig$ uniquely determines a reducing sphere, upto isotopy in $S^3$. Two reducing spheres of $\mH$ are said to be isotopic with respect to $\mH$, if any of their isotopy in $S^3$ restricts to an isotopy of the corresponding reducing curves on the surface $\sig$.  Likewise, if two reducing curves for $\sig$ are isotopic on $\sig$, this isotopy can be extended to an isotopy of the corresponding reducing spheres with respect to $\mH$. This allows one to study reducing curves on $\sig$ in order to study reducing spheres of $\mH$.

The Heegaard splitting $\mH$ of $S^3$ admits the following four automorphisms, as described in Scharlemann \citep{scharlemann}, which preserve the structure of $\mH$. The automorphisms $\alpha$ and $\gamma$ are the automorphisms of $S^3$ induced by the $\pi$-rotations of $\sig$ about the axes indicated in Figures \ref{fig:alpha} and \ref{fig:gamma}. The automorphism $\beta$ is the half twist of the $\sig$ about the central separating curve as indicated in Figure \ref{fig:beta}. 
		\begin{figure}[h!]\centering
		\begin{subfigure}{.4\linewidth}
			\includegraphics[width=\linewidth]{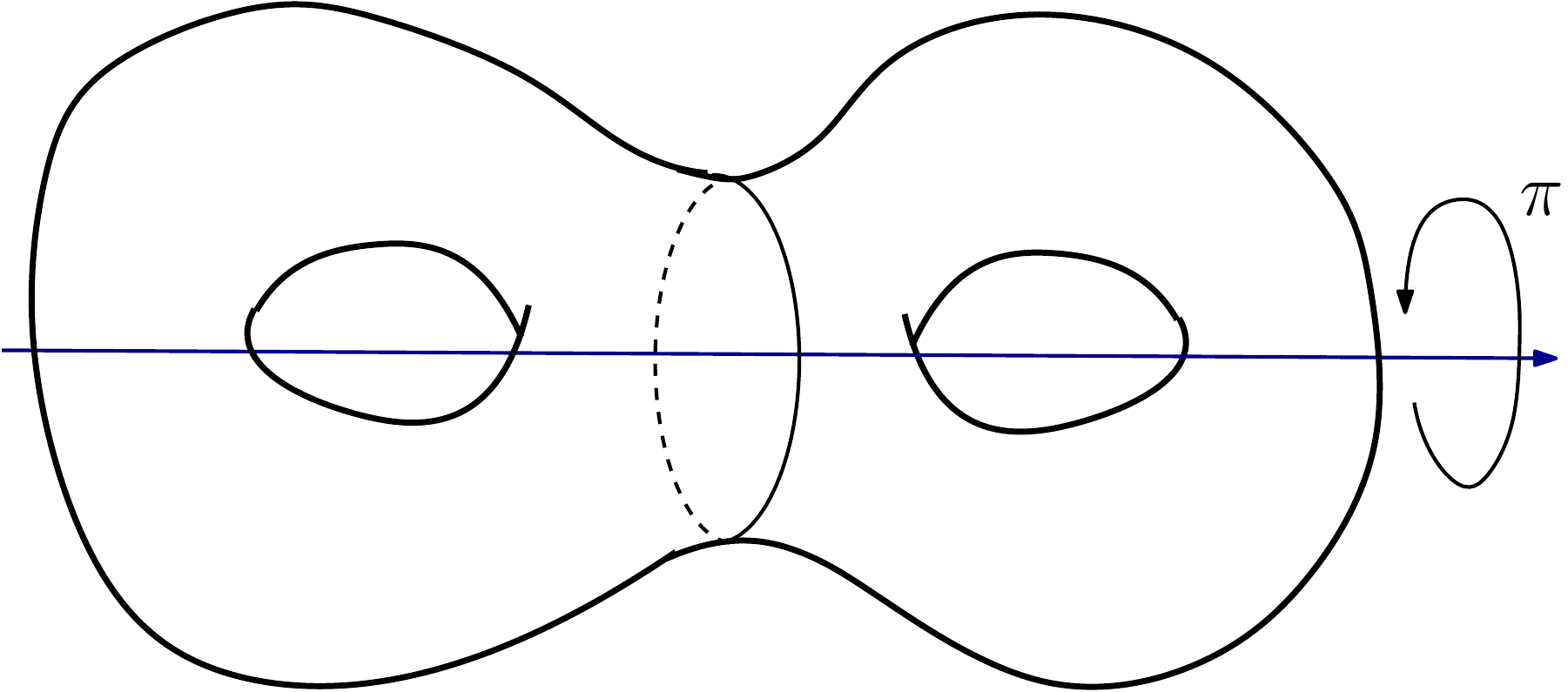}
			\caption{Automorphism $\alpha$}\label{fig:alpha}
		\end{subfigure}
		\begin{subfigure}{.4\linewidth}
			\includegraphics[width=\linewidth]{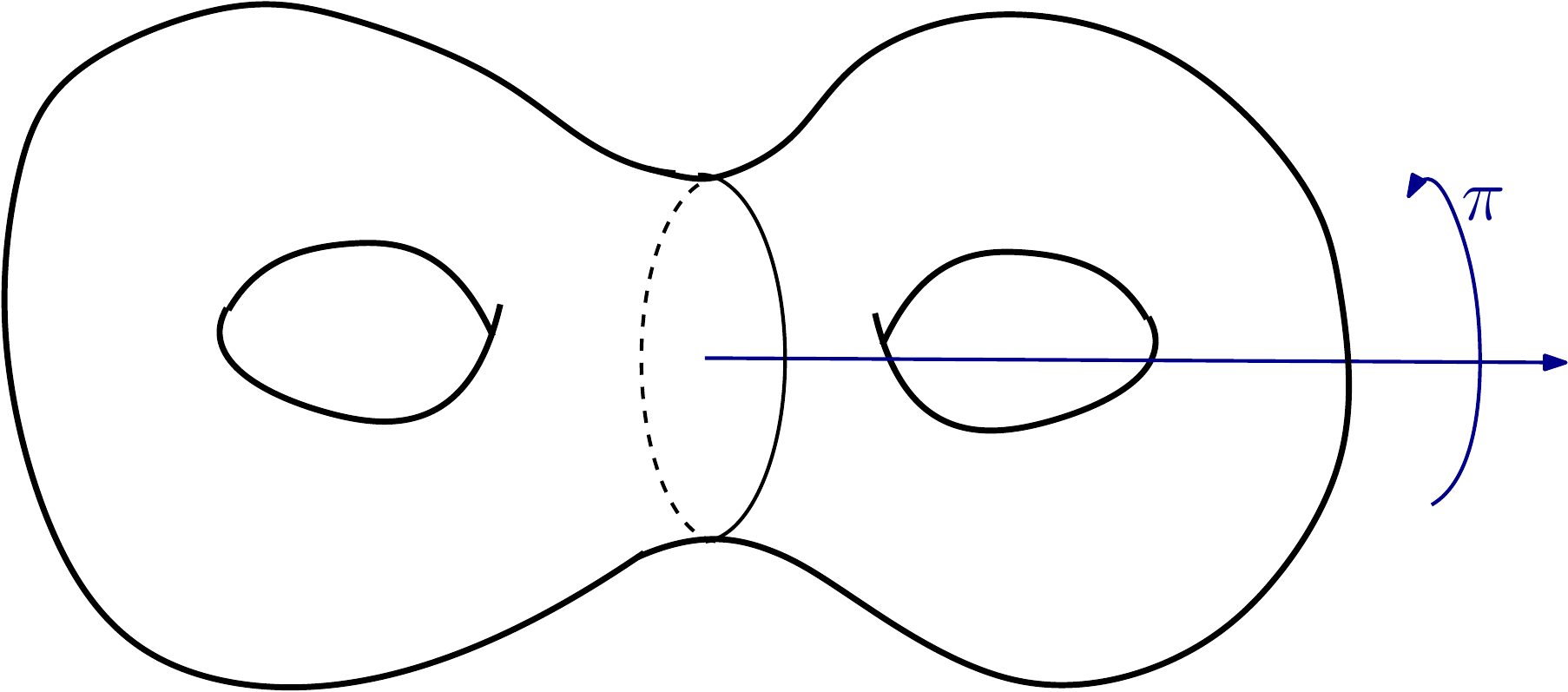}
			\caption{Automorphism $\beta$}\label{fig:beta}
		\end{subfigure}
		\begin{subfigure}{.4\linewidth}
			\includegraphics[width=\linewidth]{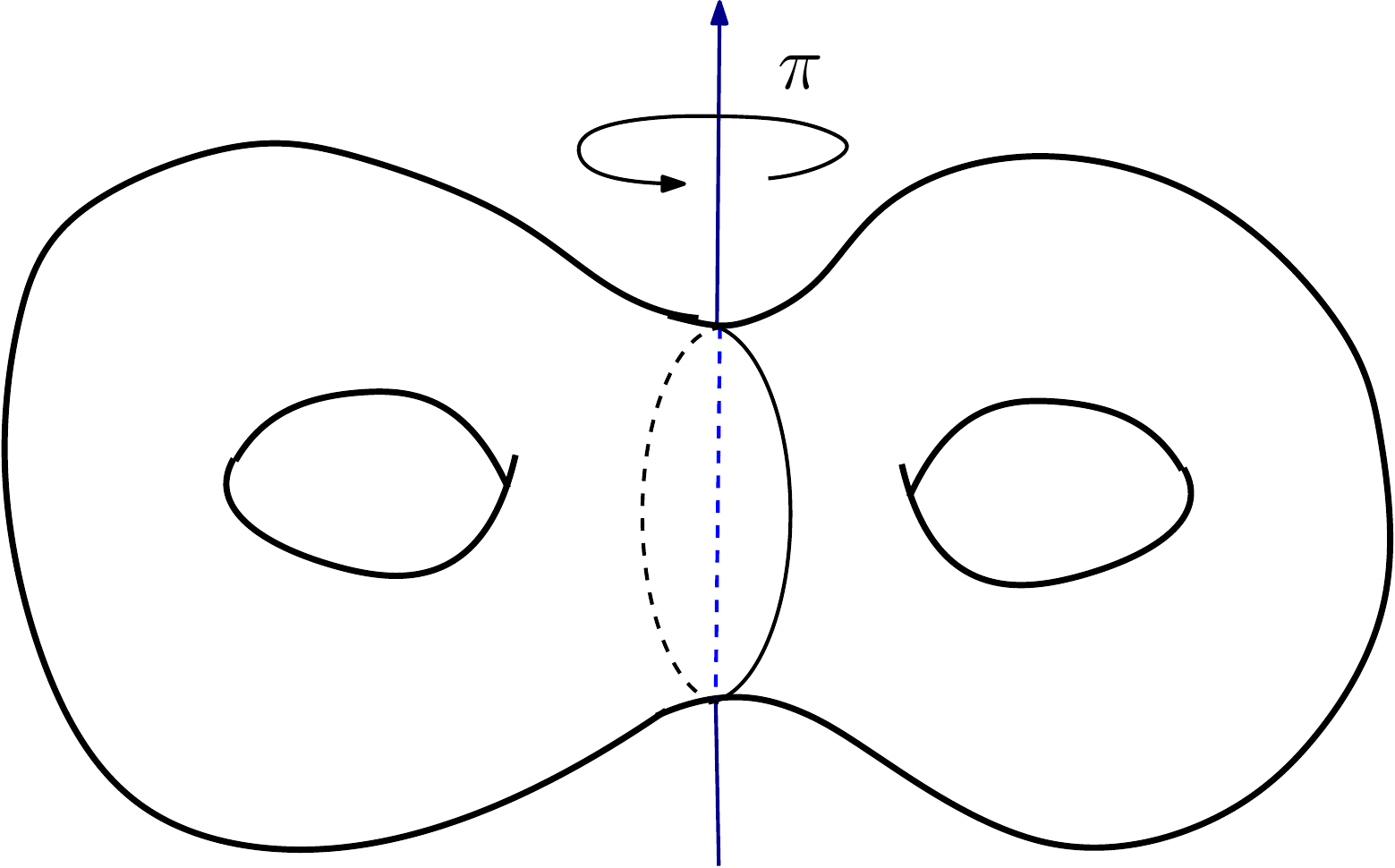}
			\caption{Automorphism $\gamma$}\label{fig:gamma}
		\end{subfigure}
		\begin{subfigure}{.4\linewidth}
			\includegraphics[width=\linewidth]{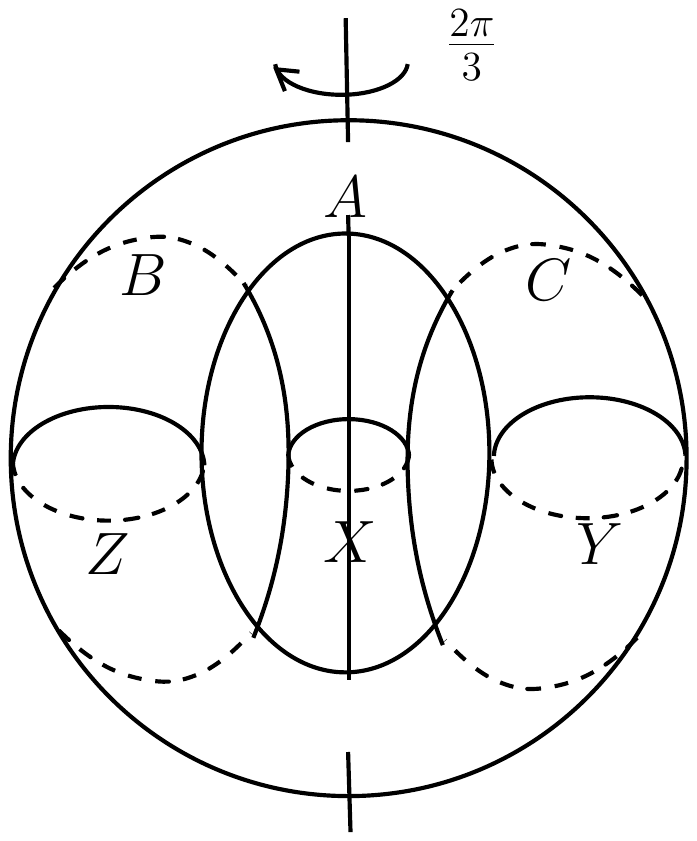}
			\caption{Automorphism $\delta$}\label{fig:delta}
		\end{subfigure}
		\caption{The generators of the genus-2 Goeritz group of $S^3$}\label{fig:automorphisms}
		\end{figure}

The automorphism $\delta$ is of order three which can be described as follows. Consider a graph on the unit sphere in $\mathbb{R}^3$ whose vertices are the north and the south pole and whose edges are three distinct half meridians on the sphere connecting the north pole and the south pole such that the rotation of the sphere by an angle $\frac{2\pi}{3}$ about the axis connecting the north and south pole permutes the three edges cyclically. If $V'$ is a regular neighborhood in $S^3$ of this graph such that $V'$ is homeomorphic to the handlebody $V$, then $\delta$ is the order three automorphism which cyclically permutes the 1-handles of $V'$ by rotation by an angle $\frac{2\pi}{3}$ as shown in Figure \ref{fig:delta}.

	Figure \ref{fig:standard_curves} shows an embedding of the Heegaard surface $\sig$ of $\mH$ in $\mathbb{R}^3$. $P$ is the reducing sphere whose reducing curve, $P_{\Sigma}$, is shown in Figure \ref{fig:standard_curves}. We call $P$ as the standard reducing sphere. Consider the non-separating curves $A,B,C,X,Y,Z$ on $\sig$ as shown in Figure \ref{fig:standard_curves}.
		\begin{figure}[h!]
			\centering\includegraphics[width=.7\linewidth]{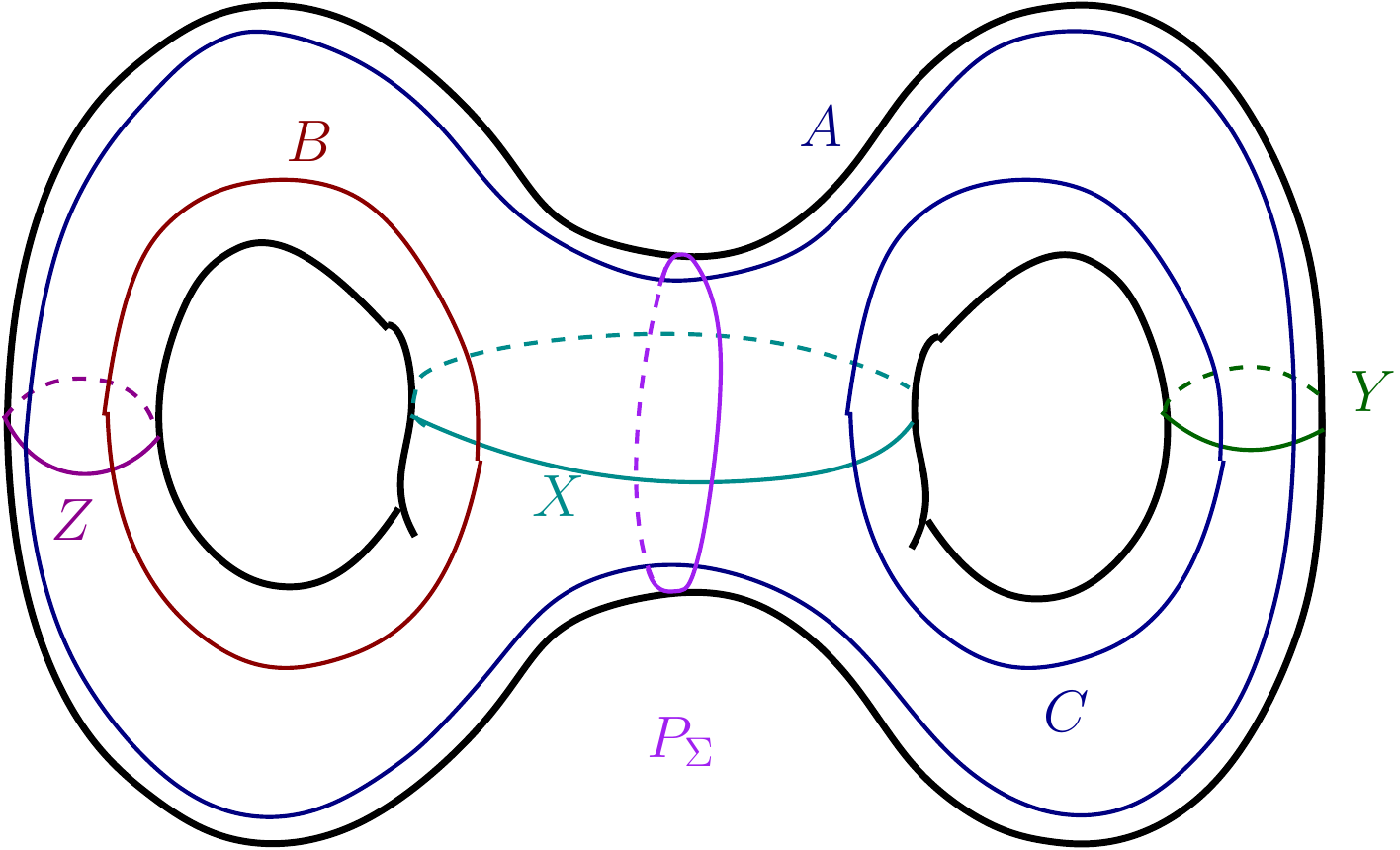}
			\caption{The standard set of curves on $\Sigma$}\label{fig:standard_curves}
		\end{figure}
	 $A\cup B\cup C$ separates $\sig$ into two thrice boundered spheres, call them $\s$ and $\t$. We will refer to the curves in Figure \ref{fig:standard_curves} throughout this article. 
	
	Any reducing sphere $Q$ of $\mH$ intersects $\sig$ in the corresponding reducing curve, which we denote by $Q_{\sig}$. Such a $\qsig$ is an essential separating simple closed curve on $\sig$ which separates $\sig$ into two surfaces, each of which is a surface of genus $1$ with one boundary component. Since $\qsig$ cannot be null-homotopic on $\Sigma$, $\qsig$ will intersect both $\s$ and $\t$ and will also intersect at least one of $A, B$ or $C$.  For any $Q$, we assume that $\qsig$ intersects $A, B, C, X, Y$ and $Z$ minimally and transversally (refer to Figure \ref{fig:standard_curves}). $\qsig \cap \s$  or $\qsig \cap \t$ is a collection of simple essential arcs on the thrice boundered spheres $\s$ or $\t$, respectively. For $D, E \in \{A, B, C\}$, \textit{an arc of type} $DE$, or a $DE$ \textit{arc}, for short, on $\s$ (or on $\t$) is a connected arc of $\qsig$ which is contained in $\s$ (or in $\t$) and which has its end-points on the curves $D$ and $E$. So there are six possible types of arcs of $\qsig$ on $\s$ or $\t$ namely an $AA$ arc, a $BB$ arc, a $CC$ arc, an $AB$ arc, an $AC$ arc and a $BC$ arc. Some arc types intersect the others. So a reducing curve, $\qsig$, cannot have intersecting arc types.
	
	Table \ref{arc-table} shows arc types which necessarily intersect and hence cannot be components of the same reducing curve on $\s$ and likewise on $\t$:
	\begin{table}[h!]
		\centering
		\begin{tabular}{|c|c|c|c|c|c|c|}\hline
			 & $AA$ & $AB$ & $AC$ & $BB$ & $BC$ & $CC$ \\ \hline
			$AA$ &  &  &  & X & X & X\\ \hline
			$AB$ &  &  &  &  &  & X \\ \hline
			$AC$ &  &  &  & X &  &  \\ \hline
			$BB$ & X &  & X &  &  & X \\ \hline
			$BC$ & X &  &  &  &  &  \\ \hline
			$CC$ & X & X &  & X &  &  \\ \hline
		\end{tabular}
		\caption{X indicates arc types necessarily intersect}\label{arc-table}
	\end{table}

    For an arbitrary reducing sphere $Q$, the number of arcs of $\qsig \cap \s$ and of $\qsig \cap \t$ have to be equal, as a mismatch will leave an open end which is impossible for a closed curve. For $D, E \in \{A, B, C\}$, let $|DE|_{Q,1}$ and $|DE|_{Q,2}$ denote the number of arcs of $\qsig$ of type $DE$ on $\s$ and $\t$, respectively. For instance, $|AB|_{Q,1}$ denotes the number of arcs of $\qsig$ of type $AB$ on $\s$ and $|BB|_{Q,2}$ denotes the number of $BB$ arcs of $\qsig$ on $\t$. When the context is clear, we will drop $Q$ in the subscript. For instance, we will just use $|BB|_2$, instead of $|BB|_{Q,2}$. We define $a_Q:= |\qsig \cdot A|, b_Q:= |\qsig \cdot B|$ and $c_Q:= |\qsig \cdot C|$. Since $\qsig$ is a separating curve, $a_Q, b_Q$ and $c_Q$ are non-negative even numbers. 
    
    The four automorphisms of $\mH$ described above, namely, $\alpha, \beta, \gamma$ and $\delta$ affect the numbers $a_Q, b_Q$ and $c_Q$ of a reducing sphere $Q$ of $\mH$ as follows: 
\begin{enumerate}
\item $\delta$ cyclically permutes $a_Q, b_Q$ and $c_Q$,
\item $\gamma$ permutes $b_Q$ and $c_Q$ and keeps $a_Q$ unchanged,
\item $\alpha$ keeps $a_Q, b_Q$ and $c_Q$ unchanged and
\item $\beta$ or its inverse increase or decrease $a_Q$, while keeping $b_Q$ and $c_Q$ unchanged.
\end{enumerate}
    
    Table \ref{tab:arc-contribution} lists the contribution of each arc type to $a_Q$, $b_Q$ and $c_Q$.
		\begin{table}[!h]
			\centering
			\begin{tabular}{|c|c|c|c|}\hline
				Type of arc & \multicolumn{3}{c|}{Each arc's contribution to}\\ \cline{2-4}
				& ~$a_Q$~ & ~$b_Q$~ & ~$c_Q$ ~\\ \hline
				$AA$ & +2 & 0 & 0 \\ \hline
				$AB$ & +1 & +1 & 0 \\ \hline
				$AC$ & +1 & 0 & +1 \\ \hline
				$BB$ & 0 & +2 & 0 \\ \hline
				$BC$ & 0 & +1 & +1 \\ \hline
				$CC$ & 0 & 0 & +2\\
				\hline
			\end{tabular}
			\caption{Contribution of each type of arc to $a_Q,b_Q$ and $c_Q$}
			\label{tab:arc-contribution}
		\end{table}

\begin{lem}\label{aQ-ne-bQpluscQ}
	For any reducing sphere $Q$, let $l, m$ and $n$ denote any distinct non-negative integers in the set $\{a_Q, b_Q, c_Q\}$. Then $ l \ne m + n $.
\end{lem}
	\begin{proof}
		By the order three symmetry of $\mH$, it suffices to prove the statement for $l = a_Q, m = b_Q$ and $n = c_Q$.
		If $Q$ is isotopic to $P$ with respect to $\mH$, i.e. $\qsig$ is isotopic to $P_{\sig}$ on $\sig$ then $a_Q=2$, $b_Q+c_Q=0$ and the lemma holds.\\
		Now suppose $\qsig$ is not isotopic to $P_{\sig}$. Using Table \ref{tab:arc-contribution} for $\qsig$ on $\s$ we have, $a_Q = 2|AA|_1 + |AB|_1 + |AC|_1$, $b_Q = |AB|_1 + 2|BB|_1 + |BC|_1$ and $c_Q = |AC|_1 + |BC|_1 + 2|CC|_1$. If $a_Q=b_Q+c_Q$ then $ |AA|_1 = |BB|_1+|CC|_1+|BC|_1$.
		So if one of $|BB|_1, |CC|_1$ or $|BC|_1$ is non-zero, then $|AA|_1 \neq 0$ and vice versa which means that $\qsig$ has an $AA$ arc along with one of $BB$, $BC$ or $CC$ arc which is impossible by table \ref{arc-table}. We can also arrive at this contradiction by counting arcs on $\t$ instead of $\s$. Therefore, we should have $|AA|_i = |BB|_i = |CC|_i = |BC|_i = 0$ for $i = 1, 2$.

Therefore, for $a_Q$ to be equal to $b_Q+c_Q$, $\qsig$ can only have $AB$ and $AC$ type arcs on both $\s$ and $\t$. But, if these are the only type of arcs of $\qsig$, then following $\qsig$ starting on any point of $\qsig \cap A$, an $AB$ arc should be followed by a $BA$ arc and an $AC$ arc should be followed by a $CA$ arc and so $\qsig$ will intersect $A$ in only one orientation. Hence the absolute value of the algebraic intersection number of $\qsig$ with $A$, $B$ or $C$ will be non-zero, contradicting the fact that $\qsig$ is a separating curve on $\sig$. This shows that $a_Q \neq b_Q+c_Q$.
	\end{proof}
	
\begin{thm}\label{arc-lemma}
            For any reducing sphere $Q$ for $\mH$, we have:
		 \begin{enumerate}[(i)]
		 \item $a_Q > b_Q + c_Q$ if and only if $\qsig$ contains at least one $AA$ arc on each of $\s$ and $\t$.
		 \item $b_Q > c_Q + a_Q$  if and only if $\qsig$ contains at least one $BB$ arc on each of $\s$ and $\t$,
		\item $c_Q > a_Q + b_Q$  if and only if $\qsig$ contains at least one $CC$ arc on each of $\s$ and $\t$.
		 \end{enumerate}
		 Furthermore, exactly one of the above three holds for $Q$.
	\end{thm}
	\begin{proof}	
	By the order three symmetry of $\mH$, it suffices to prove statement (i). The other two statements follow.
	Let $i=1$ or $i=2$ for this proof. As in the proof of Lemma \ref{aQ-ne-bQpluscQ}, $a_Q = 2|AA|_i + |AB|_i + |AC|_i$, $b_Q = |AB|_i + 2|BB|_i + |BC|_i$ and $c_Q = |AC|_i + |BC|_i + 2|CC|_i$. So, $a_Q > b_Q + c_Q$ implies $|AA|_i > |BB|_i + |CC|_i + |BC|_i$. The strict inequality implies that $|AA|_i \ge 1$. Thus $a_Q > b_Q + c_Q$ implies $\qsig$ has an $AA$ arc on both $\s$ and $\t$.
		
		Conversely, suppose $\qsig$ has an $AA$ arc on $\sig_i$. Then, $|AA|_i \ne 0$ and using Table \ref{arc-table}, $\qsig$ cannot simultaneously have any of $BB, CC$ or $BC$ arcs on $\sig_i$. This implies  $|BB|_i= |BC|_i = |CC|_i = 0$. Then, 
		\[
		b_Q + c_Q = |AB|_i + |AC|_i < 2|AA|_i + |AB|_i + |AC|_i = a_Q
		\]
		
		This proves statement (i) of the lemma. 

 Now we show that exactly one of the three inequalities hold for any $Q$. For that, we consider the complement of the case $a_Q > b_Q + c_Q$ and show that $b_Q > c_Q + a_Q$ or $c_Q > a_Q + b_Q$ must hold. By Lemma \ref{aQ-ne-bQpluscQ} $a_Q \neq b_Q + c_Q$. So the complement of $a_Q > b_Q + c_Q$ is when $a_Q < b_Q + c_Q$. If $a_Q < b_Q + c_Q$, by lemma \ref{arc-lemma}, $\qsig$ cannot have any $AA$ arc on $\s$ or on $\t$. By the work of Volodin et al. \citep{volodin}, since $\qsig$ bounds a disk in the handlebody $V$, $\qsig$ must have a `wave' on $\partial V$ with respect to $Y$ or $Z$ (also refer \citep{akbas}). In other words, $\qsig$ must either have a $YBY$ arc (see figure \ref{volodin-wave-v1}) disjoint from $X\cup Z$ or a $ZCZ$ arc disjoint from $X\cup Y$. Likewise, since $\qsig$ must also bound a disk in $W$, $c_Q$ must also have a `wave' on $\partial W$ with respect to $B$ or $C$ \textit{i.e.} a $BYB$ arc disjoint from $X\cup Z$ or a $CZC$ arc disjoint from $X\cup Y$ respectively. The `wave' on $\partial W$ with respect to $B$ when viewed on $\partial V$, is a $BB$ arc on $\partial V$, since there are no $AA$ type arcs on $\s$ or on $\t$. The `wave' on $\partial W$ with respect to $B$ is shown in figure \ref{volodin-wave-w1}, \ref{volodin-wave-w2} and \ref{volodin-wave-w3} respectively.
	\begin{figure}[!h]
	\begin{subfigure}{.49\linewidth}
	\centering
	\includegraphics[width=.65\linewidth]{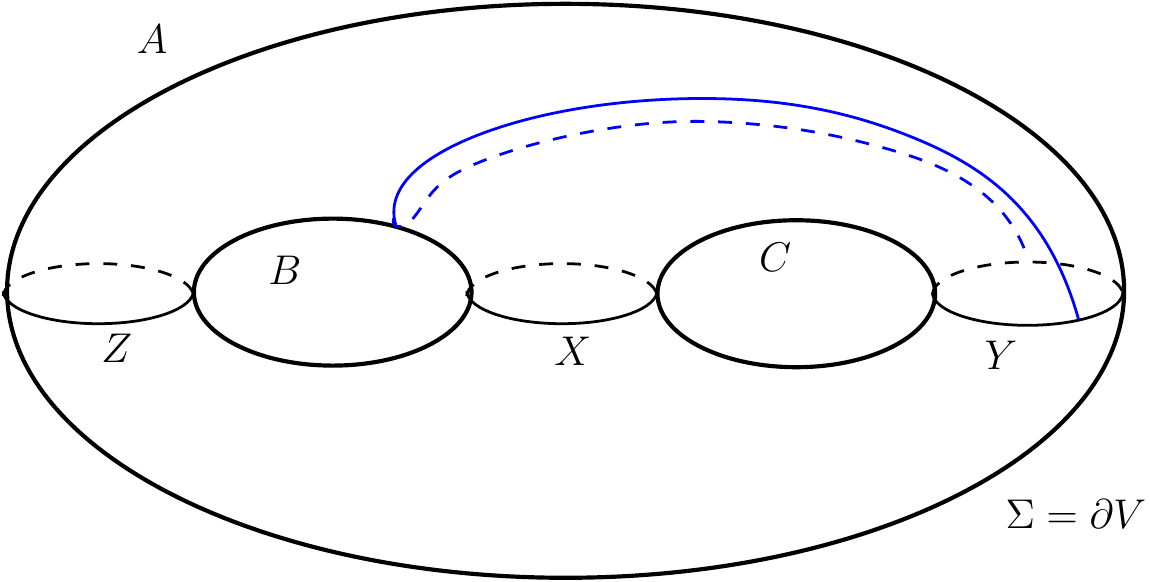}
	\caption{Wave on $\partial V$ with respect to $Y$}
	\label{volodin-wave-v1}
	\end{subfigure}\hfill
	\begin{subfigure}{.49\linewidth}
	\centering
	\includegraphics[width=\linewidth]{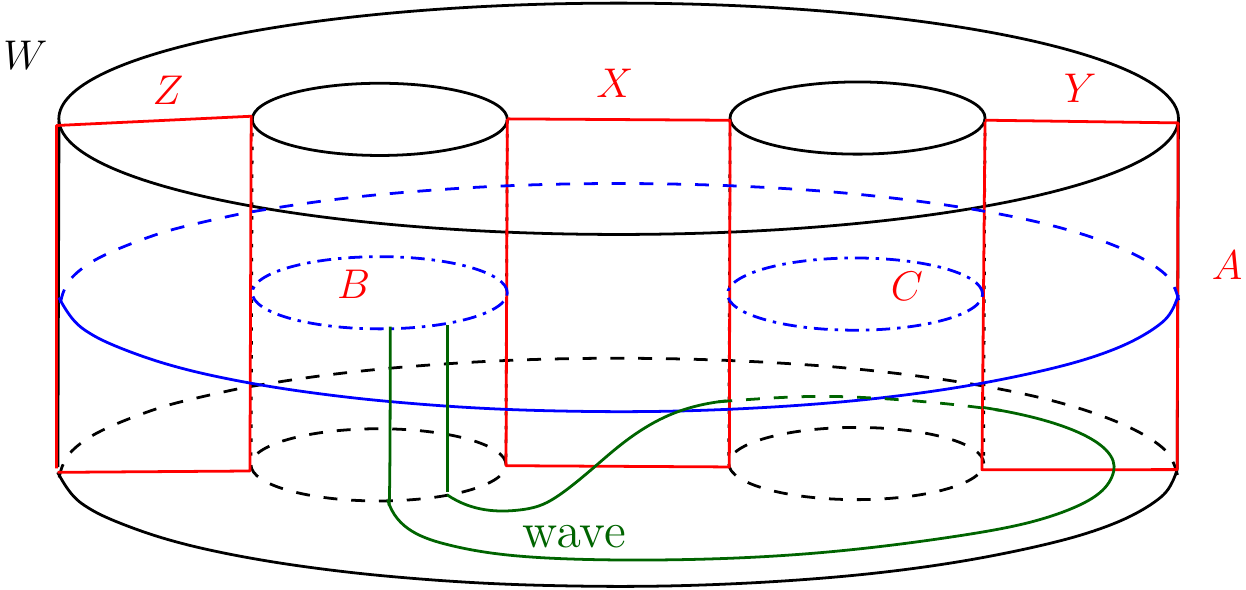}
	\caption{Wave on $\partial W$ with respect to $B$}
	\label{volodin-wave-w1}
	\end{subfigure}
	\begin{subfigure}{.49\linewidth}
	\centering
	\includegraphics[width=.65\linewidth]{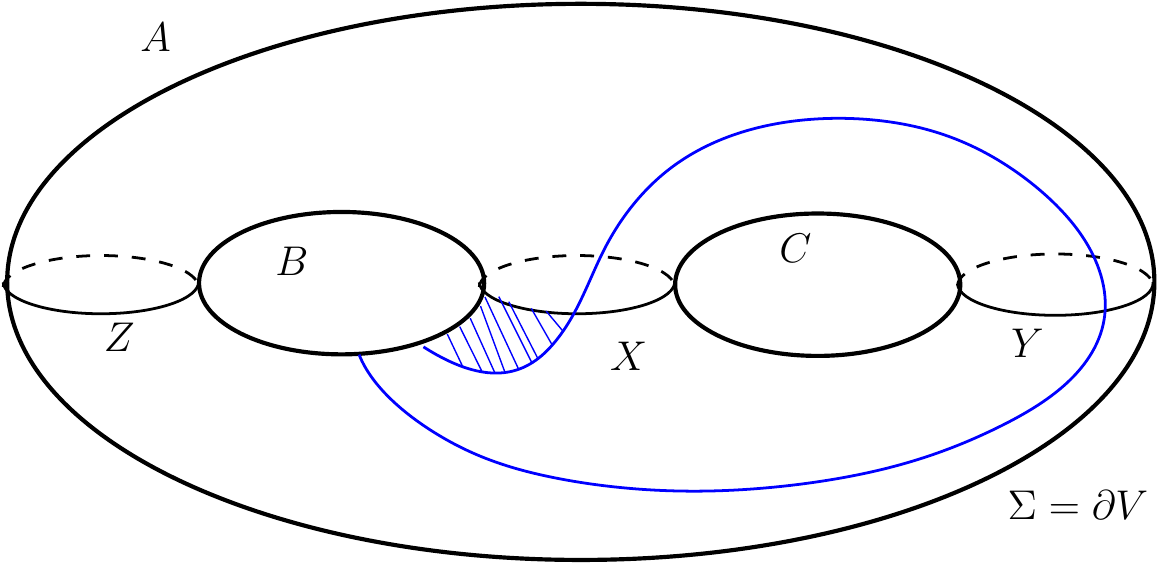}
	\caption{Wave on $\partial V$ with respect to $B$}
	\label{volodin-wave-w2}
	\end{subfigure}\hfill
	\begin{subfigure}{.49\linewidth}
	\centering
	\includegraphics[width=.65\linewidth]{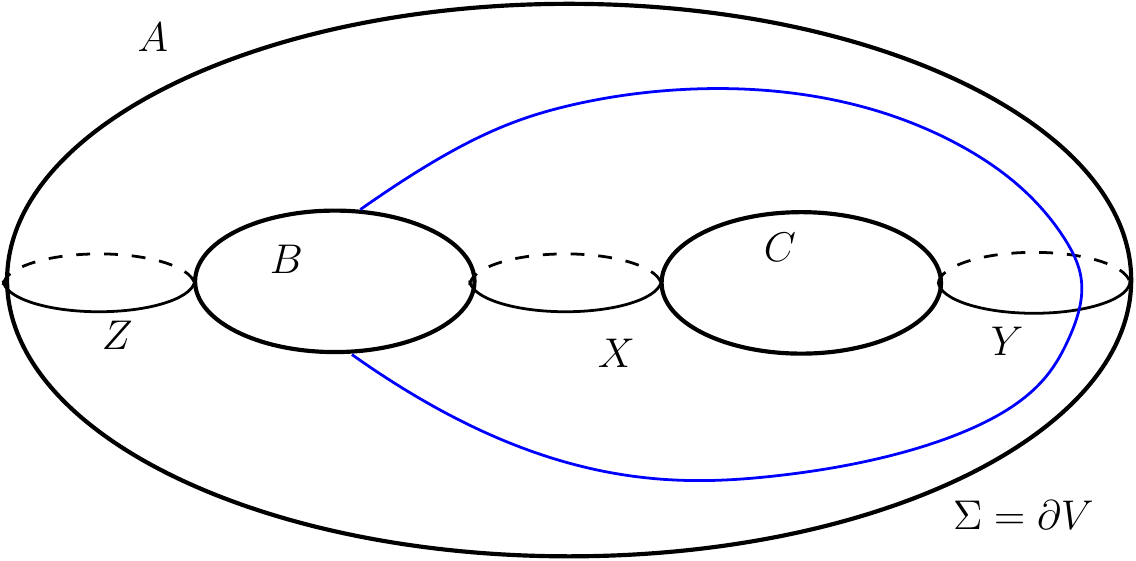}
	\caption{Wave w.r.t. $B$ after isotopy}
	\label{volodin-wave-w3}
	\end{subfigure}
	\caption{Wave of $\qsig$ w.r.t. $Y$ and $B$}
	\label{fig-volodin-wave}
	\end{figure}
       Similarly, the `wave' on $\partial W$ with respect to $C$ when viewed on $\partial V$, is a $CC$ arc on $\partial V$ since there are no $AA$ arcs on $\s$ or on $\t$.

	Therefore, when $a_Q < b_Q + c_Q$, by an isotopy, we can assume that $\qsig$ has a $BB$ or a $CC$ arc on one of $\s$ or $\t$. Without loss of generality let us assume that this wave is on $\s$. By Table \ref{arc-table}, $\qsig$ can have only one of $BB$ or $CC$ arc types on $\s$.
	
	First suppose that $\qsig$ has a $BB$ arc on $\s$. In this case, by Table \ref{arc-table} $\qsig$ cannot have $AA$, $CC$ or $AC$ arcs on $\s$. Then, since, $|BB|_1 > 0$,
	\[
	b_Q = 2|BB|_1+|BC|_1+|AB|_1 > |BC|_1 + |AB|_1= c_Q + a_Q.
	\]
	Likewise if $\qsig$ has a $CC$ arc on $\s$ then by Table \ref{arc-table} $\qsig$ cannot have $AA$, $BB$ or $AB$ arcs on $\s$. Then, since, $|CC|_1 > 0$,
	\[
	c_Q = 2|CC|_1+|BC|_1+|AC|_1 > |BC|_1 + |AC|_1= b_Q + a_Q.
	\]	
	\end{proof}
\begin{rem}\label{remark:aa_wave}
We note that by the work of Volodin et al. \citep{volodin}, when the inequality $a_Q > b_Q + c_Q$ holds, there is a `wave' on $\partial W$ with respect to $A$, which when viewed on $\partial V$ is an $AA$ type arc which either surrounds $B$ on $\s$ or on $\t$ or surrounds $C$ on $\s$ or on $\t$. Such an arc is shown in Figure \ref{fig:ayxa_azxa_arcs}.
\end{rem}

\begin{figure}\centering
\begin{subfigure}{.45\linewidth}
		\includegraphics[width=\linewidth]{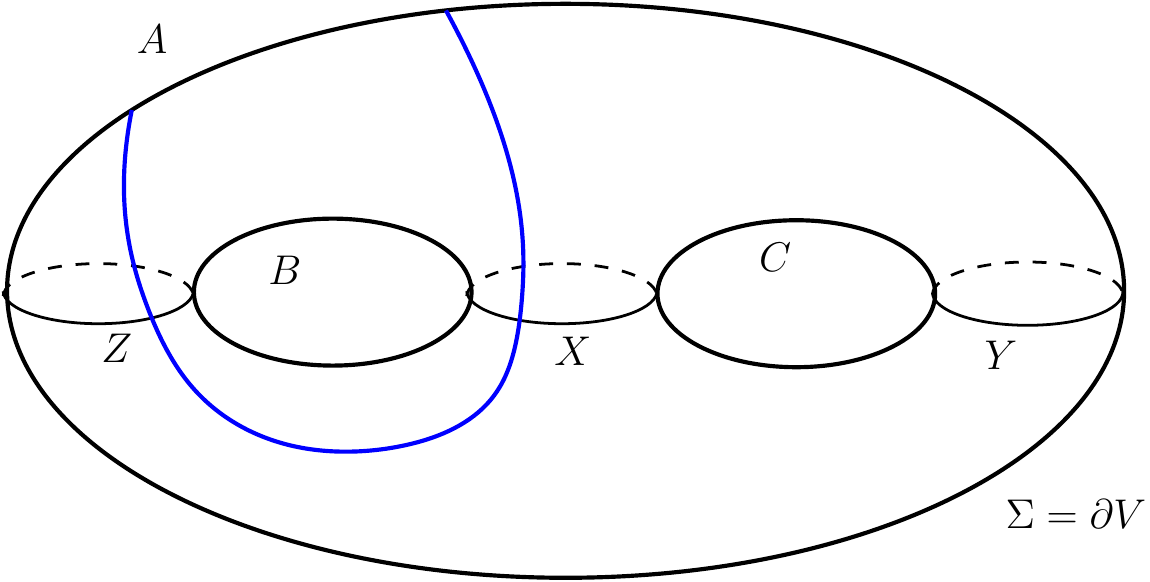}
		\subcaption*{AZXA arc}
	\end{subfigure}\hfill
	\begin{subfigure}{.45\linewidth}
		\includegraphics[width=\linewidth]{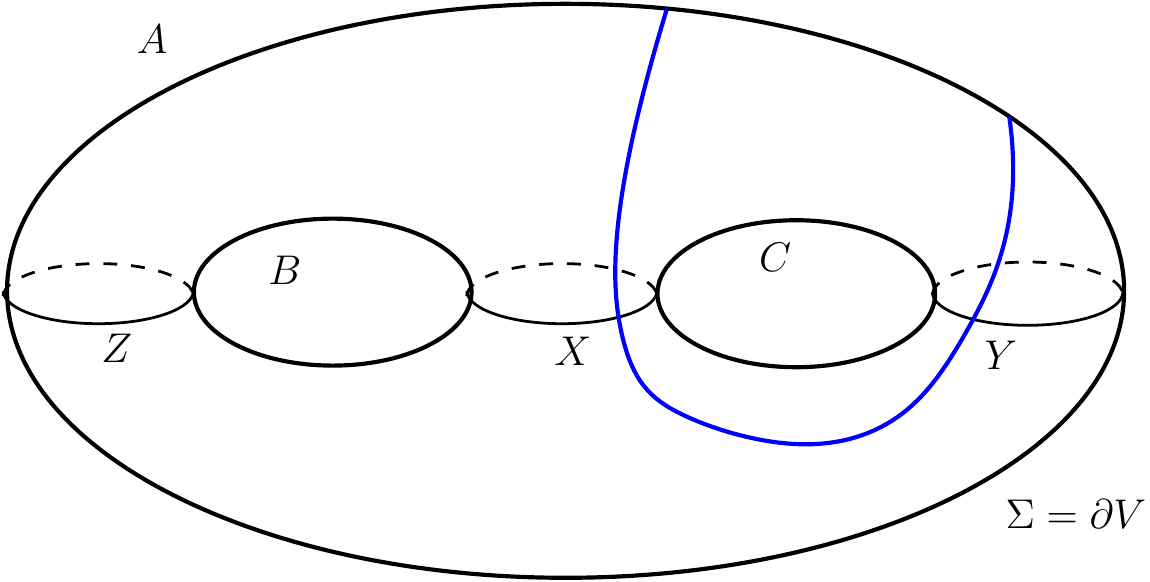}
		\subcaption*{AYXA arc}
	\end{subfigure}\hfill
\caption{A wave on $\partial W$ w.r.t. $A$ viewed on $\partial V$}\label{fig:ayxa_azxa_arcs}
\end{figure}

The proof of Theorem \ref{arc-lemma} also shows the following: 
\begin{cor}\label{cor:cor_to_arc_lemma}
(i) If a reducing sphere $Q$ satisfies $a_Q > b_Q + c_Q$, then $Q$ contains only $AA, AB$ and $AC$ type arcs on $\s$ and on $\t$.\\
(ii) If a reducing sphere $Q$ satisfies $b_Q > a_Q + c_Q$, then $Q$ contains only $BB, AB$ and $BC$ type arcs on $\s$ and on $\t$.\\
(iii) If a reducing sphere $Q$ satisfies $c_Q > b_Q + a_Q$, then $Q$ contains only $CC, AC$ and $BC$ type arcs on $\s$ and on $\t$.\\
In any of the above three cases, the number of arcs of any type are equal on $\s$ and on $\t$.
\end{cor}

\section{Arcs of a reducing curve}\label{section:arcs_pi1}

We assume the setup of the previous section. By Theorem \ref{arc-lemma}, and by invoking the order three symmetry, $\delta$, of $\mH$, in order to describe an arbitrary reducing sphere, $Q$, of $\mH$, it is enough to describe a reducing sphere in the case $a_Q > b_Q + c_Q$. Also, by Theorem \ref{arc-lemma}, the inequality $a_Q > b_Q + c_Q$ implies the existence of an $AA$ arc on $\s$ and an $AA$ arc on $\t$. Further, by Remark \ref{remark:aa_wave}, there must be an $AA$ arc on $\s$ or on $\t$ which is as shown in Figure \ref{fig:ayxa_azxa_arcs}. 

So, throughout this section, we assume that a reducing sphere $Q$ satisfies the inequality $a_Q > b_Q + c_Q$ and describe how the $AA$, $AB$ and $AC$ arcs of the corresponding reducing curve $\qsig$ can be presented on $\s$ and on $\t$.

\begin{figure}[h!]
	\centering
		\includegraphics[height=2in]{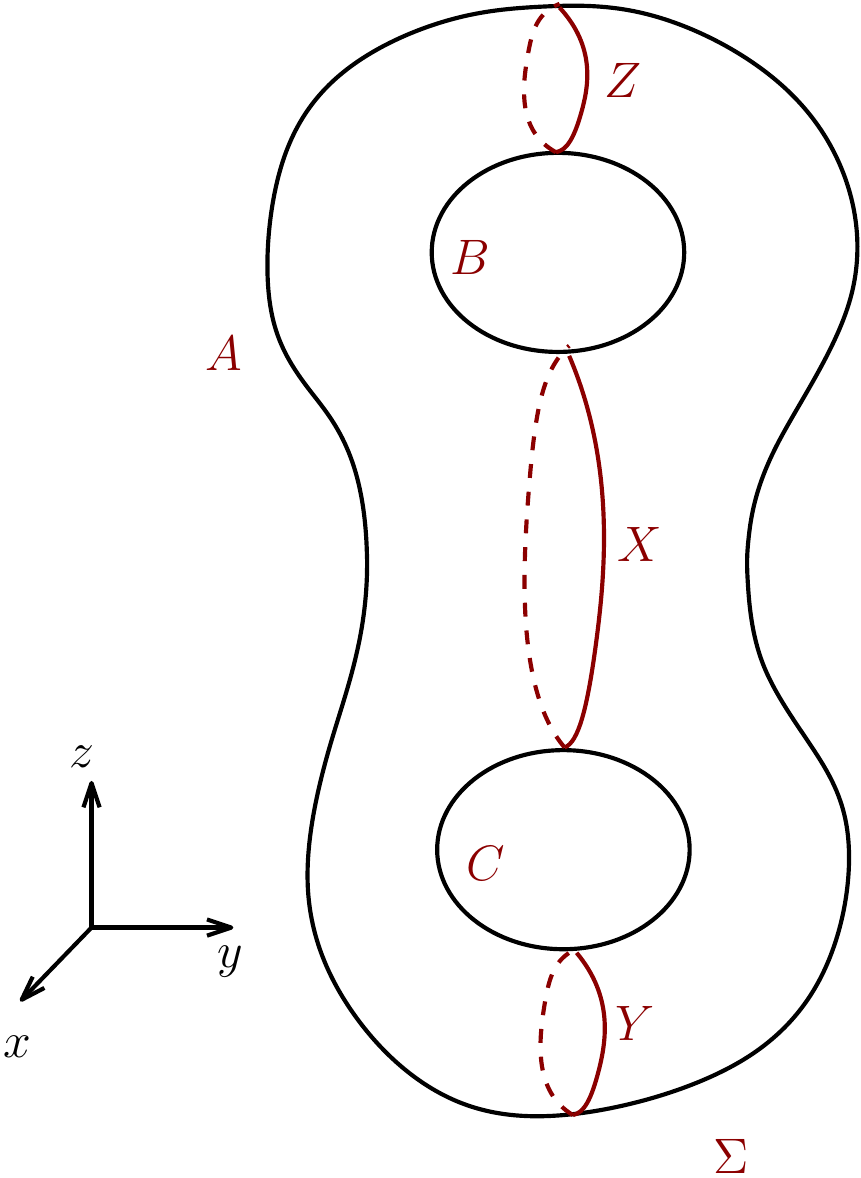}
	\caption{Setup of $\sig$ with $z$ coordinate as height}\label{fig:height}
\end{figure}
Let the two distinct arcs of $A$ each of which has one end-point on $Y$ and another end-point on $Z$ be denoted by $A_1$ and $A_2$. 

We can arrange the standard embedded genus-2 surface, $\partial V$, in $\mathbb{R}^3$ such that: (i) $A$ is the intersection of the $y-z$ plane with $\sig$ (ii) the points $Y\cap A$ and $Z\cap A$ lie along the $z$-axis with the $z$-coordinate of points on $A$ increasing from $Y \cap A$ to $Z\cap A$, (iii) no two points of $A_1$  or no two points of $A_2$ have the same $z$-coordinate. We will refer to this arrangement throughout this section. With this arrangement the $z$-coordinate can be thought of as a height function on $A_1$ and likewise on $A_2$. See Figure \ref{fig:height}. Owing to this observation, if $x_1,x_2\in A_1$ and if the height of $x_1$ is greater than the height of $x_2$, then we say $x_1$ is above $x_2$ or that $x_2$ is below $x_1$. If $x_3\in A_1$ such that the height of $x_3$ is between the heights of $x_1$ and $x_2$ then we say that $x_3$ is in between $x_1$ and $x_2$.

By an isotopy of $\qsig$, we can assume that it does not intersect the points $Z \cap A, Z \cap B, Y \cap A$ or $Y \cap C$. Since $|\qsig \cap A|$ is finite, there are only finitely many $AA$ arcs on $\s$ or on $\t$. These arcs are essential on $\s$ or on $\t$ as inessential arcs on a thrice boundered sphere are boundary reducible and a boundary reducing disk gives an isotopy of $\qsig$ reducing the intersection number of $\qsig$ with $A$.  Also these $AA$ arcs do not intersect each other as $\qsig$ is a simple curve. Since every essential arc on a thrice-boundered sphere from a boundary component to itself is separating, these $AA$ arcs are separating arcs of $\s$ and $\t$.

The curve $\qsig$ is a finite collection of $2n$ arcs on $\s$ and $\t$ for some natural number $n$. Fix an orientation of $\qsig$ and number these arcs sequentially as $\rho_1,\rho_2,\ldots, \rho_{2n}$ following the orientation of $\qsig$ such that $\rho_{i} \cap \rho_{j} \in A \cup B \cup C$ if and only if $i, j \in \{ 1,2, ..., 2n\}$ and $i-j \equiv \pm 1 \; \textrm{mod}\; 2n$. Without loss of generality we assume that for each $k\in\N$, $\rho_{2k-1}$ is on $\s$ and $\rho_{2k}$ is on $\t$.

\begin{lem}\label{beta-rem1}
	For $i = 1, 2$, an $AA$ arc of $\qsig$ on $\sig_i$ intersects $X \cap \sig_i$ exactly once.
\end{lem}
\begin{proof}
	Let $X_1 := X \cap \s$ Cut $\s$ along $X_1$ to get an annulus $S$, which has two boundary components. One boundary component is the curve $A$ and the other boundary component is a union of $B$ and $C$ with two copies of $X_1$. For brevity, we denote these two boundary components of $S$ by $\partial_A S$ and $\partial_X S$ respectively.
	
	If an $AA$ arc does not intersect $X_1$, then it is an arc in the annulus $S$ with both its endpoints on the same boundary component of $S$ and hence is boundary-reducible. The boundary-reducing disk is a bigon as the $AA$ arc intersects $A$ only at its endpoints. Then by the bigon criterion, we have an isotopy of $\qsig$ reducing its intersection with $A$ contradicting the minimal intersection position of $\qsig$ with $A$. So an $AA$ arc on $\s$ has to intersect $X_1$ at least once.
	
	Suppose now that an $AA$ arc, call it $\lambda$, intersects $X_1$ more than once. Orient $\lambda$ from one end point on $A$ to another and let $x_1, x_2,..., x_n$ be the points of intersection of $\lambda$ with $X_1$ listed in order when following the orientation of $\lambda$. $S$ cuts $\lambda$ into its component arcs on $S$.  All such component arcs on the annulus $S$ run from $\partial_X S$ to itself and hence are boundary reducible. Two boundary-reducing disks of such component arcs of $\lambda$ are either disjoint or one is contained in the other, otherwise $\lambda$ will have self-intersections. Hence by following a chain of containment of these boundary reducing disks, we can get an innermost disk $E$ which does not contain any other disk. We will now show that $E$ is a bigon on the surface $\sig$ formed by $\lambda$ and $X_1$.
	
	Let $E$ be the boundary reducing disk for an arc $\lambda_1$ on $S$ which is a component arc of $\lambda$ on $S$.  Let $\lambda_1$ considered on $\s$ join $x_i$ to $x_j$.  $E$ on $S$ is bounded by $\lambda_1$ and a portion, call it $\rho$, of the boundary component $\partial_X S$ of $S$. If $\rho$ includes $B$ or $C$ or both, then following the arc $\lambda$, in the orientation of $\lambda$ or in the opposite orientation, beyond $\lambda_1$, we arrive at yet another component arc of $\lambda$, call it $\lambda_2$, on $S$ which enters $E$. See Figure \ref{fig:AA_int_X}. 
	\begin{figure}[h!]
		\centering\includegraphics[width=\linewidth]{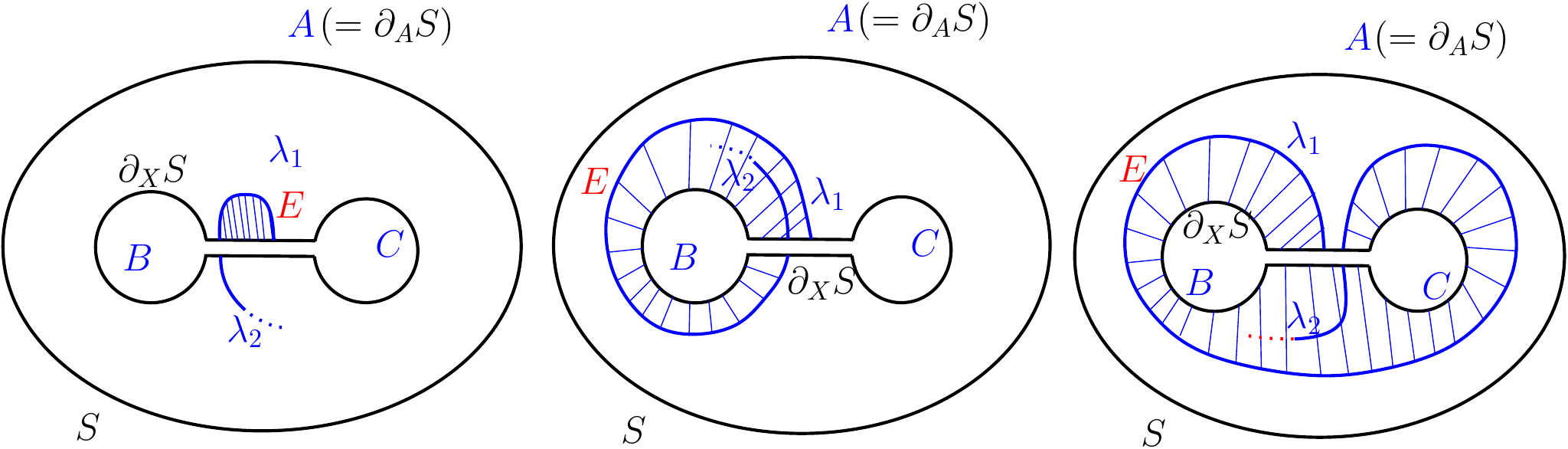}
		\caption{Any $AA$ arc intersects $X$ exactly once}\label{fig:AA_int_X}
	\end{figure}
	
	Since $\lambda_2$ enters $E$, it is completely contained in $E$. This is because, $\lambda_2$ cannot intersect $\lambda_1$, $B$ or $C$, so it has to intersect $X_1$ which then means that $\lambda_2$ is contained in $E$. Now, as remarked earlier, $\lambda_2$ is also boundary reducible and the boundary reduction disk of $\lambda_2$ to the boundary $\partial_X S$ is completely contained inside $E$, contradicting the assumption that $E$ is the innermost disk. This proves that $\rho$ cannot include $B$ or $C$ or both.
	
	This implies that $\rho$ contains only an arc from $X_1$. So $\lambda_1$ along with $\rho$ forms a bigon on the surface $\sig$. So by the bigon criterion there is an isotopy of $\qsig$ reducing the intersection with $X_1$. So if $\qsig$ is in minimal intersection position with $X$, then an $AA$ arc on $\s$ cannot intersect $X$ more than once.
	
	Hence, we conclude that any $AA$ arc of $\qsig$ on $\s$ intersects $X_1$ in exactly one point. 
	
	Likewise, we can repeat the same argument for an $AA$ arc on $\t$ and conclude that such an arc intersects $X \cap \t$ in exactly one point.
\end{proof}

We now use the fundamental group of the handlebody $V$ to get more information about the $AA$, $AB$ and $AC$ arcs of $\qsig$ when $a_Q > b_Q + c_Q$. In particular we will show that an $AA$ arc of $\qsig$ on $\s$ or on $\t$ can wind around the circles $B$ or $C$ at most once. 

Let $D_X, D_Y$ and $D_Z$ be the essential disks in $V$ bounded by the curves $X, Y$ and $Z$ respectively. Consider the core curves of the handlebody $V$, dual to the disks $D_Z$ and $D_Y$, based at a point $x_0$ lying in the interior of $D_X$ and denote them by $B_{x_0}$ and $C_{x_0}$ respectively. Note that $B_{x_0}$ and $C_{x_0}$ are freely isotopic in $V$ to $B$ and $C$ respectively. The element in $\pi_1(V, x_0)$ whose representative is a loop which traces $B_{x_0}$ once in the direction shown in Figure \ref{fig:pi_1_hbdy} will be denoted by $b$ and the element which traces $C_{x_0}$ once in the direction shown in Figure \ref{fig:pi_1_hbdy} will be denoted by $c$.
\begin{figure}[h!]
	\centering\includegraphics[height=2in]{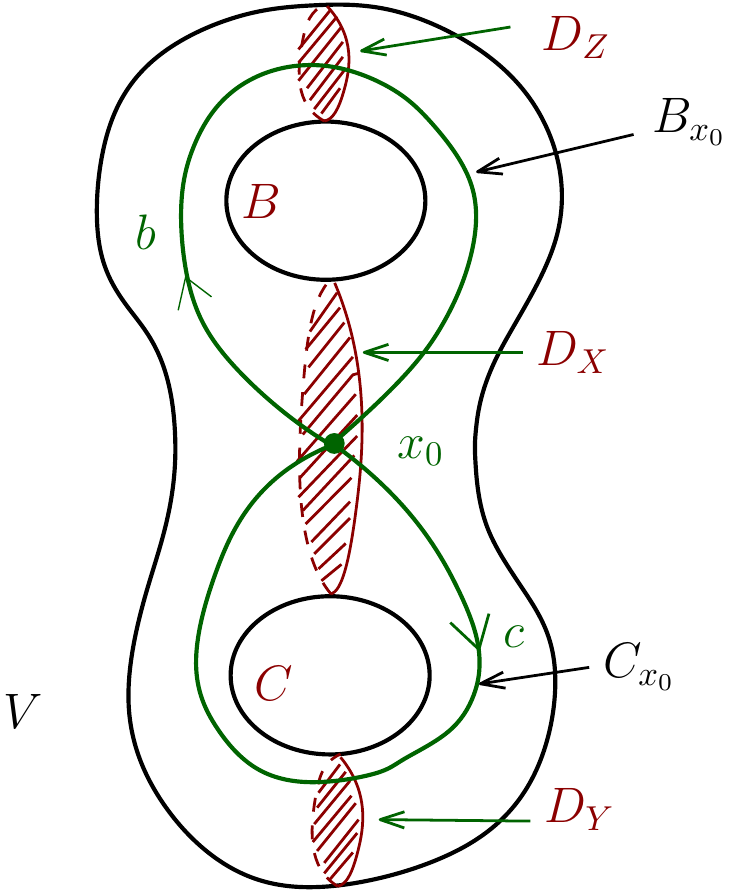}
	\caption{Handlebody $V$: core curves and dual disks}\label{fig:pi_1_hbdy}
\end{figure}

We consider the presentation $\langle b,c|-\rangle$ of the fundamental group of $V$ based at $x_0$ and henceforth refer to this presentation simply as $\pi_1(V)$. With a slight abuse of notation, we use the same letters to indicate the elements in $\pi_1(V)$ and the loops in $V$. By this abuse of notation the loop $B_{x_0}$ in Figure \ref{fig:pi_1_hbdy} is also denoted by $b$, for instance. 

Now, let $Q$ be a reducing sphere which is in minimal position with the curves $A, B, C, X, Y$ and $Z$ and which satisfies $a_Q > b_Q + c_Q$. Consider some orientation of $\qsig$ and view it as a union of oriented arcs $\rho_1, \rho_2, \ldots, \rho_{2n}$ as mentioned earlier, where the orientation on $\rho_i$, $1 \leq i \leq 2n$, is the one induced from the orientation of $\qsig$. Let $i$ denote an integer in the set $\{1, 2 ,..., 2n\}$ throughout the following discussion with addition and subtraction performed modulo $2n$ when indices are added or subtracted. Each $\rho_i$ is an $AA$, $AB$ or an $AC$ arc with a starting point indicated by its orientation. Since $\rho_i$ is an oriented arc we would distinguish between an $AC$ arc and a $CA$ arc, and likewise between an $AB$ and a $BA$ arc depending upon the orientation of $\rho_i$. Let $\delta_i$ be an arc contained in the interior of $V$ oriented from $x_0$ to the starting point of $\rho_i$ such that $\delta_i$ does not intersect $D_Y$ or $D_Z$. By $\overline{\delta_i}$, we mean $\delta_i$ with the reversed orientation. Now, $\sigma_i := \delta_i \cdot \rho_i \cdot \overline{\delta_{i+1}}$ is a loop in $V$ based at $x_0$. This loop, $\sigma_i$, represents a word in the in the generators $b$ and $c$ in $\pi_1(V)$, which we call the \textit{word of $\rho_i$}.

Also, $\qsig$ is freely homotopic to the product loop $\prod_{i=1}^{2n} \sigma_i$ and correspondingly represents a product word in the generators $b$ and $c$ in $\pi_1(V)$. Since $\qsig$ bounds a disk in $V$, the corresponding reduced word must be cyclically trivial. This last requirement imposes restrictions on the arcs $\rho_i$'s which we discuss in the following. 

The union of curves $A \cup B \cup C \cup X \cup Y \cup Z$ divides $\sig$ into four disks, which we call the \textit{component hexagons of $\sig$}, owing to the fact that the boundary of each of these disks contain exactly one arc each of the six curves $A, B, C, X, Y$ and $Z$. In this context, an arc of $\rho_i$ on a component hexagon of $\sig$ is of the type $DE$ if it starts on the curve $D$ and ends on the curve $E$ in the orientation of $\rho_i$, where $D, E \in \{A, B, C, X, Y, Z\}$. A sub-arc of $\rho_i$ is called its \textit{initial-arc} if it is completely contained in one of the component hexagons of $\sig$ and is of the type $DE$ where $D \in \{A, B, C\}$ and $E \in \{X, Y, Z, A, B, C\}$. Similarly, a sub-arc of $\rho_i$ is called its \textit{terminal-arc} if it is completely contained in one of the component hexagons of $\sig$ and is of the type $DE$ where $D \in \{X, Y, Z, A, B, C\}$ and $E \in \{A, B, C\}$. Note that neither the terminal nor the initial arc of $\rho_i$ of type $DE$ can have $D = E$.

\begin{lem}\label{lem1-ZAZ-bigon}
Consider any two component hexagons $D_1$ and $D_2$ of $\sig$ both of which have a side $A_1$ in common or both of which have a side $A_2$ in common. If $D_1$ intersects $\qsig$ in an $AZ$ or a $ZA$ arc, then $D_2$ cannot intersect $\qsig$ in an $AZ$ or a $ZA$ arc. Likewise, both $D_1$ and $D_2$ cannot simultaneously intersect $\qsig$ in a $YA$ or an $AY$ arc.
\end{lem}
\begin{proof}
	\begin{figure}[h!]
		\begin{subfigure}{.43\linewidth}
			\centering\includegraphics[height=2in]{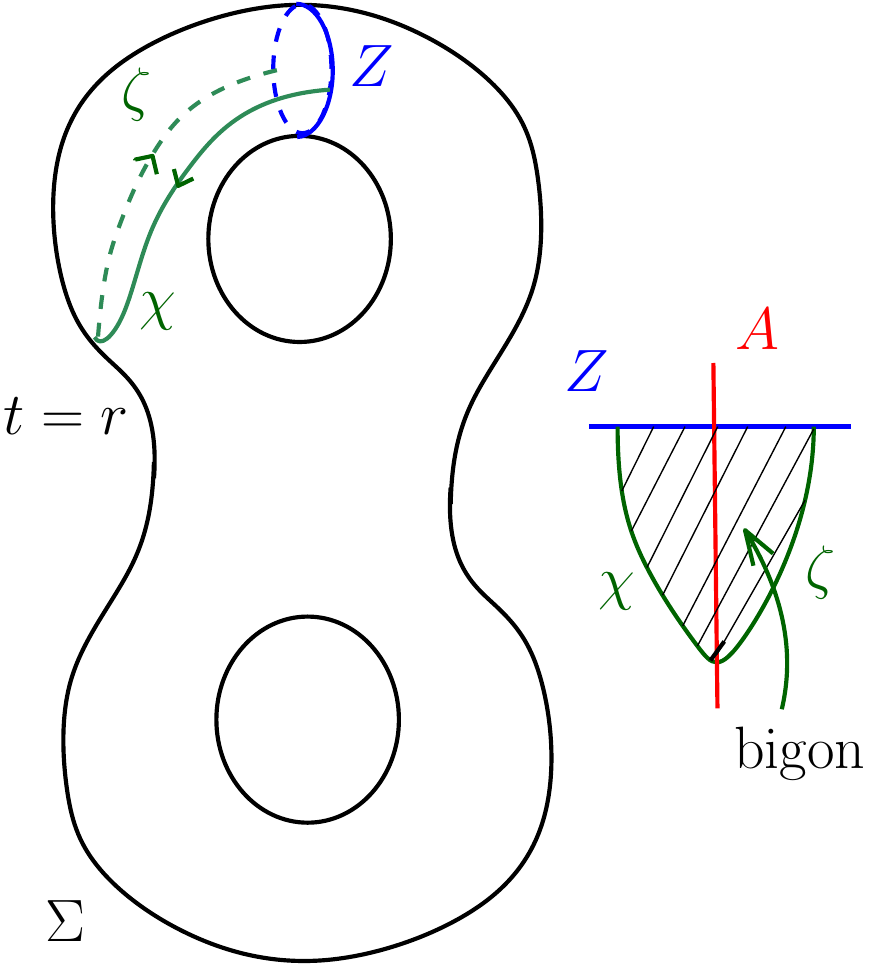}
			\caption{$t$ and $r$ coincide}\label{fig:ZAZ-bigon-case1}
		\end{subfigure}
		\begin{subfigure}{.28\linewidth}
			\centering\includegraphics[height=2in]{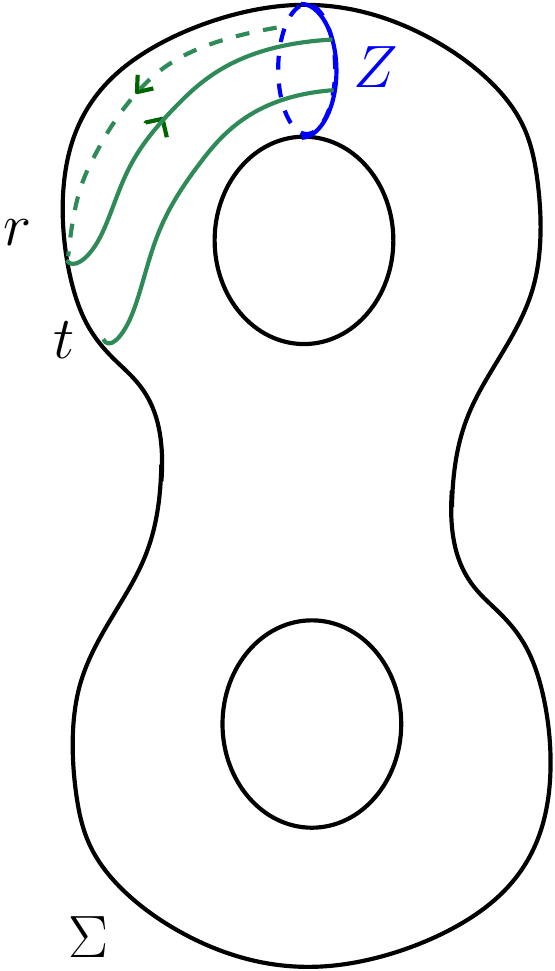}
			\caption{$t$ is below $r$}
		\end{subfigure}
		\begin{subfigure}{.27\linewidth}
			\centering\includegraphics[height=2in]{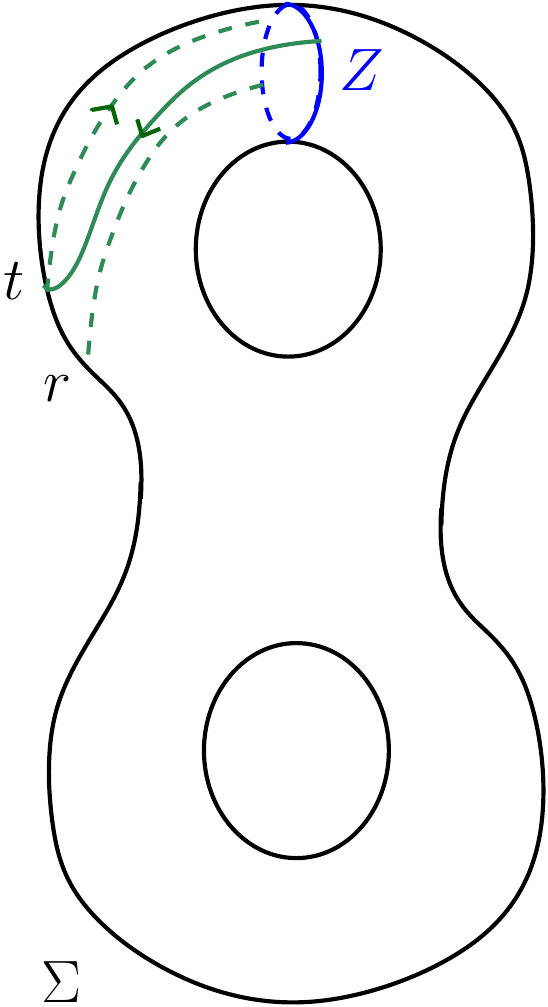}
			\caption{$t$ is above $r$}	
		\end{subfigure}
		\caption{A Z-A-Z bigon}\label{fig:ZAZ-bigon}
	\end{figure}
	
	Let $\chi$ be an arc of $\qsig$ in $D_1$ with one endpoint on $Z$ and the other on $A_1$. Let $t := \chi \cap A_1$. If possible let $\zeta$ be an arc on $D_2$ with end-points on $Z$ and $A_1$ and $r := \zeta \cap A_1$. 
	If $t = r$, then $\chi$ and $\zeta$ bound a bigon with $Z$ on $\sig$ contradicting the minimal intersection position of $\qsig$ with $Z$. 
	If $t$ is below $r$ on $A_1$, then the continuation arc of $\zeta$ on $D_1$ is contained in the disk cut out of $D_1$ by $\chi$, $Z$ and $A_1$ and hence is another arc $\zeta_1$ with one endpoint as $r$ and the other on $Z$. But then, $\zeta$ and $\zeta_1$ bound a bigon with $Z$ on $\sig$, once again contradicting the minimal intersection position of $\qsig$ with $Z$.
	If $r$ is below $t$, then the continuation arc of $\chi$ on $D_2$ and $\chi$ together bound a bigon with $Z$ leading us to a similar contradiction. See Figure \ref{fig:ZAZ-bigon}.
	
	So $r$ cannot be above or below $t$ on $A_1$ nor can it coincide with $t$ and hence $\zeta$ does not exist. The proof for the case where $\chi$ is an arc of $\qsig$ in $D_1$ with one endpoint on $Y$ and the other on $A_1$ is similar.
\end{proof}

\begin{cor}\label{bigon-cor1}
	When the word of a $\rho_i$ which is an $AA$, $BA$ or a $CA$ arc of $\qsig$ is concatenated with word of $\rho_{i+1}$, which is an $AA$, $AB$ or an $AC$ arc the concatenated word cannot have trivial relators of $\pi_1(V)$.
\end{cor}

\begin{proof}
	This is so because a trivial relator at the word interface contradicts Lemma \ref{lem1-ZAZ-bigon}. 
	
	Suppose that $\rho_i$ is an $AA$, $BA$ or a $CA$ arc on $\s$ whose word ends with a letter $b^{-1}$, and that $\rho_{i+1}$ is an $AA$, $AB$ or an $AC$ arc on $\t$ whose word starts with a letter $b$. See Figure \ref{fig:ZAZ-bigon-case1}. The terminal-arc of $\rho_i$ on $\s$ is a $ZA$ arc and the initial-arc of $\rho_{i+1}$ on $\t$ is an $AZ$ arc, contradicting Lemma \ref{lem1-ZAZ-bigon}.
	
	Similarly, if $\rho_i$ is an arc on $\s$ with a word which ends with any letter $x$, where $x \in \{b, b^{-1}, c, c^{-1}\}$ and the word of $\rho_{i+1}$ begins with $x^{-1}$ then $\qsig$ forms a bigon with $Y$ or $Z$ on $\sig$ and this will lead to a contradiction with the minimal position of $\qsig$ with $Y$ and $Z$. Since the roles of $\s$ and $\t$ can be reversed in the above argument, the corollary is proved.
\end{proof}

\begin{cor}\label{bigon-cor2}
	If the word of a $\rho_i$, which is an $AB$ arc (on $\s$ or $\t$), ends in $b$ or $b^{-1}$ then the word of $\rho_{i+1}$, which has to be a $BA$ arc (on $\t$ or $\s$, respectively) cannot start with $b^{-1}$ or $b$, respectively. Similarly, if the word of a $\rho_i$, which is an $AC$ arc (on $\s$ or $\t$), ends in $c$ or $c^{-1}$ then the word of $\rho_{i+1}$ which has to be a $CA$ arc (on $\t$ or $\s$, respectively) cannot start with $c^{-1}$ or $c$, respectively.
\end{cor}

\begin{proof}
	As in the proof of Corollary \ref{bigon-cor1} if $\rho_i$ is an $AB$ arc on $\s$ (or on $\t$) whose word ends with a letter $b$ or $b^{-1}$, then $\rho_{i+1}$ will be a $BA$ arc on $\t$ (or on $\s$), and its word cannot begin with $b^{-1}$ or $b$ respectively, because if it does, then the terminal-arc of $\rho_i$ and the initial-arc of $\rho_{i+1}$ form a bigon with $Z$ on $\sig$ leading to a contradiction with the assumption that $\qsig$ is in minimal position with $Z$. See Figure \ref{fig:ZBZ-bigon}.
	
	\begin{figure}[h!]
		\centering\includegraphics[height=2in]{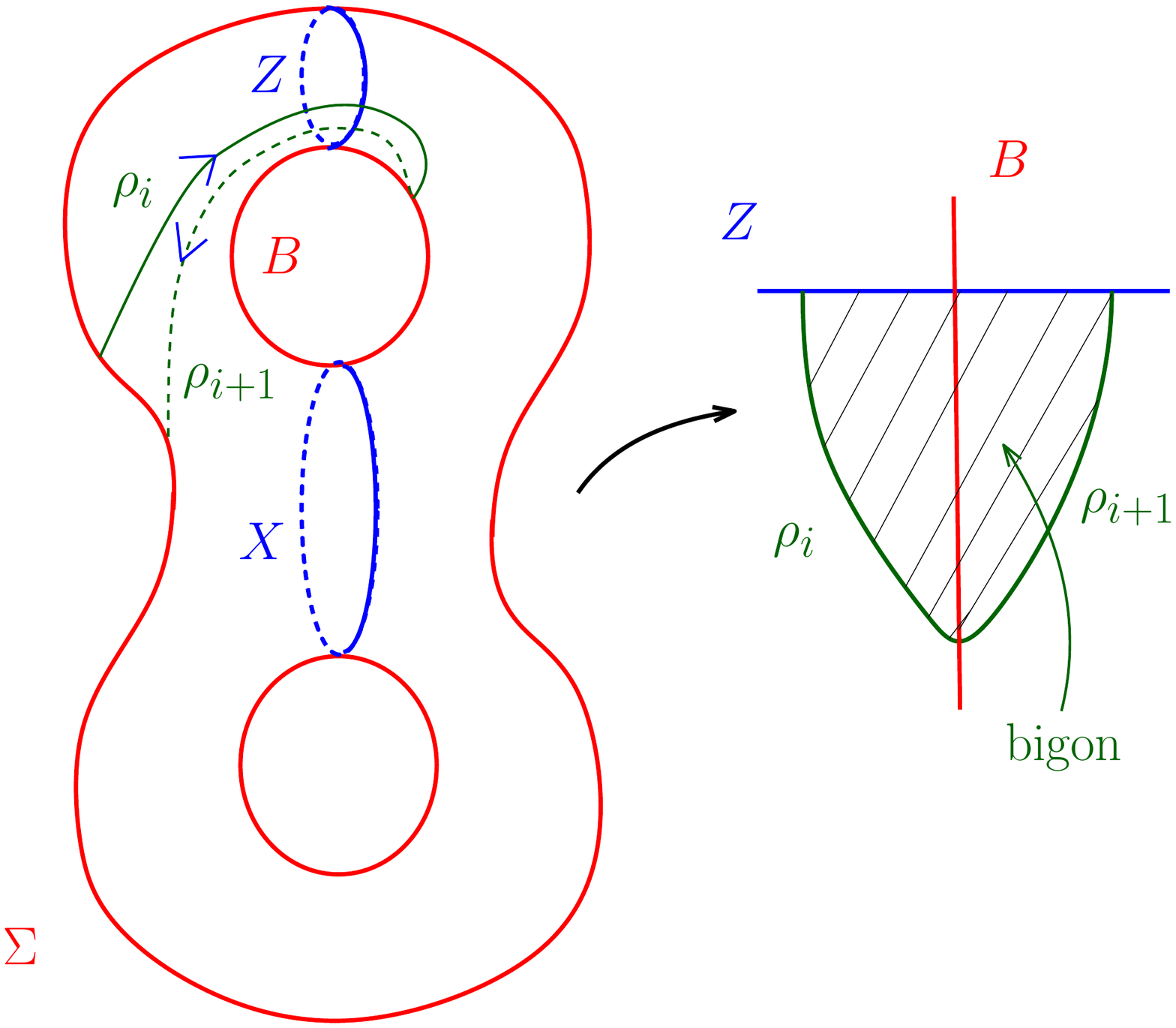}
		\caption{$AB, BA$ concatenation cannot have trivial relators }\label{fig:ZBZ-bigon}
	\end{figure}
	
	Similarly, if $\rho_i$ is an $AC$ arc on either $\s$ (or $\t$), whose word ends with a letter $c$ or $c^{-1}$, then $\rho_{i+1}$ will be a $CA$ arc on $\t$ (or $\s$), and its word cannot begin with $c^{-1}$ or $c$ respectively, because if it does, then the terminal-arc of $\rho_{i}$ and the initial-arc of $\rho_{i+1}$ form a bigon with $Y$ leading to a contradiction with the assumption that $\qsig$ is in minimal position with $Y$.
\end{proof}

\begin{lem}\label{lem2-ZABZ-bigon} 
	A concatenated word of a sequence $\rho_{i-1}, \rho_{i}$ and $\rho_{i+1}$ of arcs of $\qsig$ where $\rho_{i}$ is an $AB$ arc cannot contain the trivial relators $bb^{-1}$ or $b^{-1}b$. Similarly, a concatenated word of a sequence $\rho_{i-1}, \rho_{i}$ and $\rho_{i+1}$ of arcs of $\qsig$ where $\rho_{i}$ is an $AC$ arc cannot contain the trivial relators $cc^{-1}$ or $c^{-1}c$.
\end{lem}
\begin{proof}
	Corollary \ref{bigon-cor1} implies that a trivial relator $bb^{-1}$ or $b^{-1}b$ cannot occur when the word of $\rho_{i-1}$ is concatenated with the word of a $\rho_{i}$ when $\rho_{i-1}$ is an $AA, BA$ or a $CA$ arc and $\rho_{i}$ is an $AB$ arc. Corollary \ref{bigon-cor2} implies that a trivial relator $bb^{-1}$ or $b^{-1}b$ cannot occur when the word of $\rho_i$ is concatenated with the word of a $\rho_{i+1}$ when $\rho_i$ is an $AB$ arc and $\rho_{i+1}$ is a $BA$ arc. So we need to prove that when $\rho_{i}$ is an $AB$ arc whose word is empty word, $\{\}$, the concatenated words of $\rho_{i-1}, \rho_i$ and $\rho_{i+1}$ arc sequence cannot be of the form $w_1b\{\}b^{-1}w_2$ or of the form $w_1b^{-1}\{\}bw_2$ where $w_1$ and $w_2$ are some words in $b, b^{-1}, c$ and $c^{-1}$. But if this happens, then as in the proof of Corollary \ref{bigon-cor2}, the terminal arc of the $\rho_{i-1}$, the whole of $\rho_i$ and the initial-arc of $\rho_{i+1}$ form a bigon with $Z$ as shown in figure \ref{fig:3arc-bigon}. 
	\begin{figure}[h!]
		\centering\includegraphics[width=\linewidth]{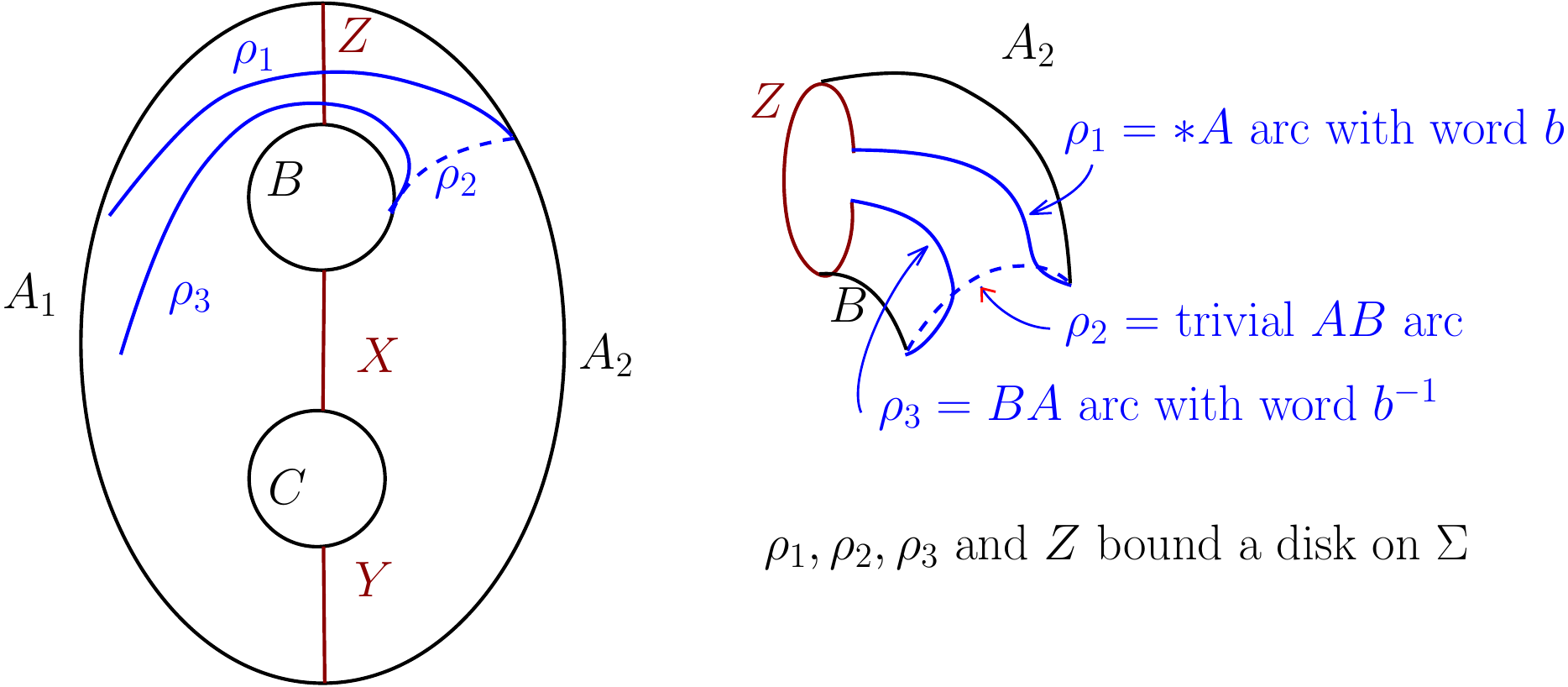}
		\caption{Concatenation of $*A, AB, BA$ arcs bounds bigon on $\sig$ with $Z$}\label{fig:3arc-bigon}
	\end{figure}
	This contradicts the assumption that $\qsig$ is in minimal position with $Z$. This proves the first statement of the lemma. The proof of the second statement of this lemma is similar.
\end{proof}

Now, we describe the words of $AA, AB, BA, AC$ and $CA$ arcs. For $j = 1, 2$, let $X_j = X \cap \sig_j, Y_j = Y \cap \sig_j$ and $Z_j = Z \cap \sig_j$.

We will first describe the words of $AA$ arcs on $\s$. The same argument applies to words of $AA$ arcs on $\t$. An $AA$ arc on $\s$ can be described completely upto isotopy by a sequence of arc types in the component hexagons of $\s$ because these hexagons are disks. Further, it is possible to give an exhaustive description of these arc sequences owing to the conditions on an $AA$ arc imposed by Lemma \ref{beta-rem1}. Let $m$ be a non-negative integer, $\epsilon_1, \epsilon_2 \in \{0,1\}$ and $k$ be any integer for the following discussion.

Table \ref{tab:pi1_words} gives all possible arc sequences of an $AA$ arc on $\s$ which starts on $A_j$ and ends on $A_{j'}$, where $j,j' \in \{1,2\}$, and gives its corresponding word in $\pi_1(V)$. Note that the arc types alternate on the two component hexagons of $\s$. The first, third, etc. arc-types are on the component hexagon which has $A_1$ as an edge and the second, fourth etc. arc-types are on the component hexagon which has $A_2$ as an edge.
\begin{table}
\centering
\begin{tabular}{c c c c}
Start & End & Arc-type sequence and Word\\
\hline\\
$A_1$ & $A_2$ & $AX, XA$\\ & & empty\\
\hline\\
$A_1$ & $A_2$ & $AZ, (ZY, YZ)^m, ZX, XY, (YZ, ZY)^m, YA$\\ & & $b(cb)^{m}(c^{-1}b^{-1})^{m}c^{-1}$\\
\hline\\
$A_1$ & $A_2$ & $AZ, (ZY, YZ)^m, ZY, YX, XZ, (ZY, YZ)^m, ZA$\\ & & $b(cb)^{m}cb^{-1}(c^{-1}b^{-1})^{m}c^{-1}$\\ 
\hline\\
$A_1$ & $A_2$ & $AY, (YZ, ZY)^m, YX, XZ, (ZY, YZ)^m, ZA$\\ & & $c^{-1}(b^{-1}c^{-1})^m(bc)^{m}b$\\
\hline\\
$A_1$ & $A_2$ & $AY, (YZ, ZY)^m, YZ, ZX, XY, (YZ, ZY)^m, YA$\\ & & $c^{-1}(b^{-1}c^{-1})^mb^{-1}c(bc)^{m}b$ \\ 
\hline\\
\end{tabular}
\caption{Arc-type sequences and Words in $\pi_1(V)$}\label{tab:pi1_words}
\end{table}

Figure \ref{fig:AA-arcs} shows these arc sequences for a sample.

\begin{figure}[h!]
	\begin{subfigure}{.3\linewidth}
		\centering
		\includegraphics[height=2in]{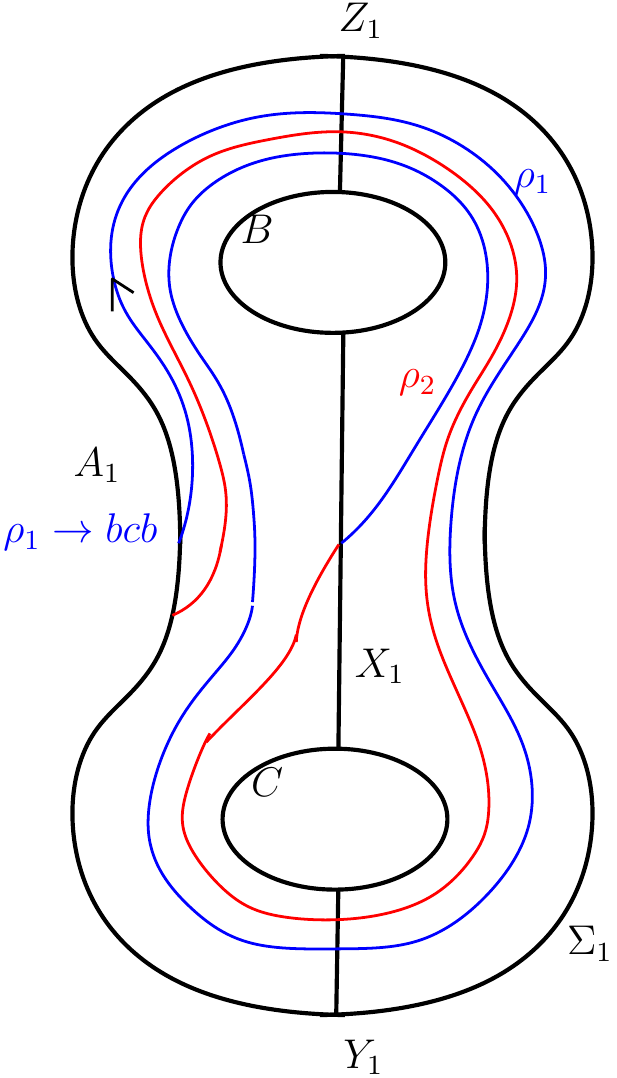}
		\caption{$A_1A_1$ arc}
	\end{subfigure}
	\begin{subfigure}{.3\linewidth}
		\centering
		\includegraphics[height=2in]{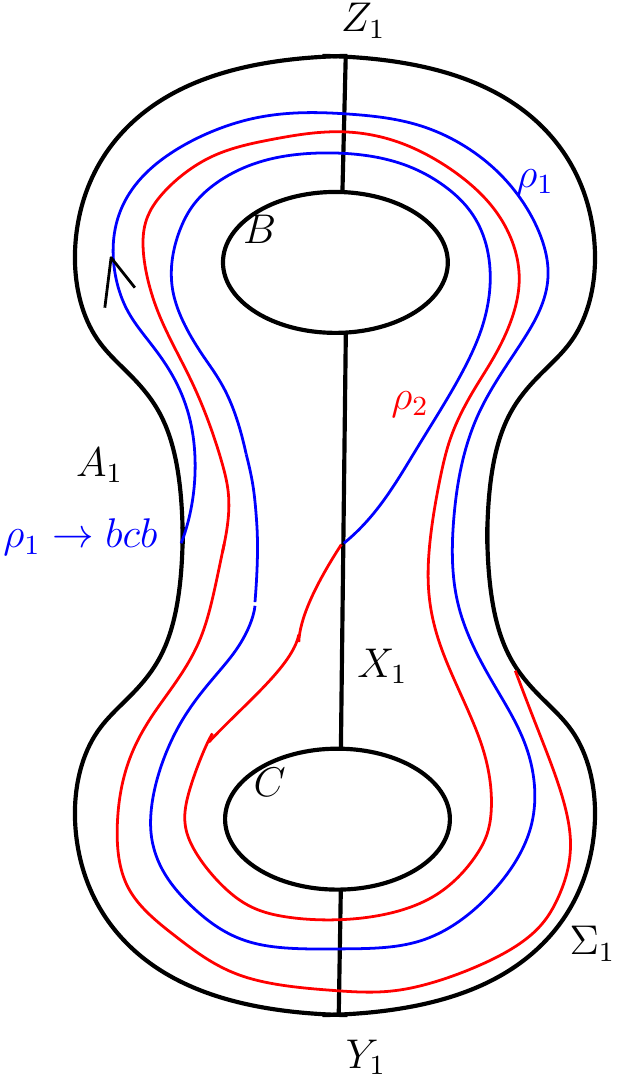}
		\caption{$A_1A_2$ arc}
	\end{subfigure}
	\begin{subfigure}{.3\linewidth}
		\centering
		\includegraphics[height=2in]{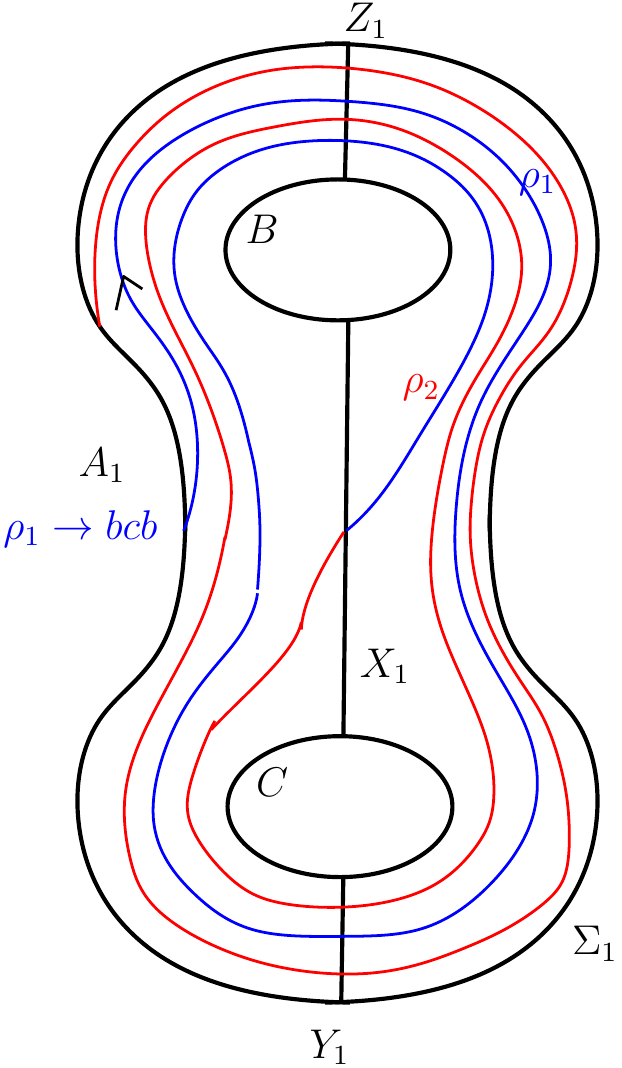}
		\caption{$A_1A_1$ arc}
	\end{subfigure}
	\caption{Arcs for $AA$ words}\label{fig:AA-arcs}
\end{figure}

	Using this technique, we summarize the reduced word forms in $\pi_1(V)$ of various arc-types as follows. Let $w$ represent some reduced word in $\pi_1(V)$. An $AA$ arc of $\qsig$ on $\s$ or on $\t$ starting on $A_1$ and ending on $A_2$ has words with three possible reduced forms, viz: empty, $bwc^{-1}$, $c^{-1}wb$. An $AA$ arc starting on $A_2$ and ending on $A_1$ has words with three possible reduced forms, viz: empty, $cwb^{-1}$, $b^{-1}wc$. An $AA$ arc starting on $A_1$ and ending on $A_1$ has words with the following possible reduced forms: $b$, $b^{-1}$, $c$, $c^{-1}$, $bwb^{-1}$, $c^{-1}wc$. An $AA$ arc starting on $A_2$ and ending on $A_2$ has words with following possible reduced forms: $b$, $b^{-1}$, $c$,  $c^{-1}$, $cwc^{-1}$, $b^{-1}wb$. The $AB, BA, AC$ and $CA$ arcs of $\qsig$ on $\s$ or on $\t$ have reduced word forms based on the word form of its $AA$ arc, which necessarily exists when $a_Q>b_Q + c_Q$. Depending on the beginning letter of the word described by an $AA$ arc, an $AB$ arc starting on $A_1$ describes a reduced word of the following forms: $b^k$, $bwcb^k$, $c^{-1}b^k$ or $c^{-1}wc^{-1}b^k$. Depending on the beginning letter of the word described by an $AA$ arc, an $AB$ arc starting on $A_2$ describes a reduced word of the following forms: $b^k$, $cb^k$, $cwcb^k$, $b^{-1}wc^{-1}b^k$. Likewise, an $AC$ arc starting on $A_1$ describes a reduced word of the following forms: $c^k$, $bc^k$, $bwbc^k$, $c^{-1}wb^{-1}c^k$. An $AC$ arc starting on $A_2$ describes a reduced word of the following forms: $c^k$, $cwbc^k$, $b^{-1}c^k$, $b^{-1}wb^{-1}c^k$. The word forms of a $BA$ arc are the inverses of the word forms of an $AB$ arc and the word forms of a $CA$ arc are the inverses of the word forms of an $AC$ arc of the corresponding kind. 
	
	Since $\qsig$ represents the conjugacy class of the trivial reduced word, the following Theorem shows that the variable word $w$ in the word-form summary, presented above, has to be trivial. 

\begin{thm}\label{thm:aa-atmost-2}
	The word in $\pi_1(V)$ of an $AA$ arc of $\qsig$ can have a length at most $2$. Correspondingly, the word of an $AB, BA, AC$  or $CA$ arc of $\qsig$ can have a word of the form $wb^k$ and $wc^k$ respectively, for some integer $k$, where $w$ can have a length at most $1$.
\end{thm}

\begin{proof}
	 	 
	When we concatenate the words of the arcs $\rho_i$ of $\qsig$ on $\s$ and $\t$ from $i = 1$ to $i=2n$, then the words should reduce to an empty word. Since the words described by various arcs, as summarized above are already in their reduced form, at least one trivial relator, a relator of the form $xx^{-1}$, where $x$ could be any of $b, c, b^{-1}, c^{-1}$, has to occur at the join of some $\rho_i$ and $\rho_{i+1}$. Referring to the summary of word-forms of arcs, if the word $w$ in the word of an arc which is of the form $xwy$ is non-trivial, where $x, y \in \{b, c, b^{-1}, c^{-1}\}$ then a trivial relator does not arise when words of arcs which form the word of $\qsig$ are concatenated, leaving us with a conclusion that $\qsig$ cannot have arcs with non-trivial words $w$.
		
	Now consider an arc $\rho_{i}$ of $\qsig$. Suppose that $\rho_{i}$ is a $BA$ or a $CA$ arc, whose word is of any of the types described above where $w$ is not the empty word.
	\begin{itemize}
		\item[\textbf{Case 1:}] {\bf $\rho_{i+1}$ and $\rho_{i+2}$ is an $AB, BA$ or an $AC, CA$ arc combination.} Suppose first that $\rho_{i+1}$ is an $AB$ arc and $\rho_{i+2}$ is a $BA$ arc. In this case, by Corollary \ref{bigon-cor1} and Lemma \ref{lem2-ZABZ-bigon}, for any trivial relator to arise, if at all it does, when words in $\pi_1(V)$ are concatenated, the words of both $\rho_{i+1}$ and $\rho_{i+2}$ must be empty. But this is impossible as both $\rho_{i+2}$ and $\rho_{i}$ are on $\s$ or on $\t$ and $\rho_{i+2}$ is a $BA$ arc which cannot describe an empty word when $w$ is not the empty word. So no trivial relator arises when concatenating the words of $\rho_{i}, \rho_{i+1}$ and $\rho_{i+2}$. A similar argument shows that no trivial relator arises when concatenating the words of $\rho_{i}, \rho_{i+1}$ and $\rho_{i+2}$ when $\rho_{i+1}, \rho_{i+2}$ is an $AC, CA$ arc combination.
		
		\item[\textbf{Case 2:}] {\bf $\rho_{i+1}$ is an $AA$ arc.} Even in this case, by Corollary \ref{bigon-cor1}, there cannot be a trivial relator when the words of $\rho_{i}$ and $\rho_{i+1}$ are concatenated. However, the word of $\rho_{i+1}$ could be empty, in which case it should be an arc from $A_1$ to $A_2$ or an arc from $A_2$ to $A_1$.
		
		\begin{enumerate}
			\item[\textbf{Subcase 1:}] \textbf{$\rho_{i}$ ends on $A_2$:} In this sub-case, the word of $\rho_{i}$ ends with $b$ or $c^{-1}$. Here we refer to the summary of words of arcs. So if the word of $\rho_{i+1}$ is empty it must be an arc from $A_2$ to $A_1$, and so $\rho_{i+2}$ has to start on $A_1$ and hence has a word starting with either $b$ or $c^{-1}$. As a result there is no trivial relator while concatenating these three words.
			
			\item[\textbf{Subcase 2:}] \textbf{$\rho_{i}$ ends on $A_1$:} In this sub-case, the word of $\rho_{i}$ ends with $b^{-1}$ or $c$. So if the word of $\rho_{i+1}$ is empty it must be an arc from $A_1$ to $A_2$, and so $\rho_{i+2}$ has to start on $A_2$ and hence has a word starting with either $b^{-1}$ or $c$. As a result there is no trivial relator while concatenating these three words.
		\end{enumerate}
	\end{itemize}
	If $\rho_{i}$ is an $AB$ or an $AC$ arc whose word is of any of the types described above where $w$ is not the empty word, then the corresponding $AB, BA$ arc combination or an $AC, CA$ arc combination can never be the empty word by Corollary \ref{bigon-cor2}.
	
	The arcs $\rho_i$ for $i=1$ to $i=2n$ of $\qsig$ are a sequence of $AA$ arcs, $AB, BA$ arc combinations and $AC, CA$ arc combinations. When the corresponding words in $\pi_1(V)$ are concatenated by following the order of $\rho_i$'s, at-least one of these pieces gives a non-trivial reduced word when $w$ is not the empty word and none of the concatenation results in a trivial relator. This shows that the word described by $\qsig$ cannot be the empty word resulting in a contradiction to the assumption that $\qsig$ bounds a disk in $V$. This shows that the word of any $AA$ arc on $\s$ or $\t$ can only be of the form $xwy$ with an empty $w$ and where $x,y\in \{b,c,b^{-1},c^{-1},\{\}\}$. Hence the words of $AA$ arc on $\s$ or $\t$ can contain at most two letters. Correspondingly, we infer from the summary of arcs above that the word of an $AB$ or of an $AC$ arc of $\qsig$ is of the form $wb^k$ and $wc^k$ respectively, for some integer $k$, where $w$ can have a length at most $1$.
\end{proof}

\begin{rem} An $AA$ arc with empty word on $\s$ (or on $\t$) will intersect an $AA$ arc with two-letter word on $\s$ (or on $\t$). So $\qsig$ cannot contain two $AA$ arcs on $\s$ (or on $\t$), one with an empty word and another with a two-letter word.
\end{rem}

\begin{rem} 
By the word forms arrived at, preceding the Theorem \ref{thm:aa-atmost-2}, we can further refine the words of $AB, BA, AC$ and $CA$ arcs based on whether their $A$-end is on $A_1$ or $A_2$ as follows.
\begin{itemize}
\item An $AB$  arc which starts on $A_1$ can only have a word of the form $b^k$ or $c^{-1}b^k$ and correspondingly a $BA$ arc which ends on $A_1$ can only have the inverse words $b^k$ or $b^kc$.
\item An $AB$  arc which starts on $A_2$ can only have a word of the form $b^k$, $cb^k$ and correspondingly a $BA$ arc which ends on $A_2$ can only have the inverse words $b^k$ or $b^kc^{-1}$.
\item An $AC$  arc which starts on $A_1$ can only have a word of the form $c^k$, $bc^k$ and correspondingly a $CA$ arc which ends on $A_1$ can only have the inverse words $c^k$ or $c^kb^{-1}$.
\item An $AC$  arc which starts on $A_2$ can only have a word of the form $c^k$, $b^{-1}c^k$ and correspondingly a $CA$ arc which ends on $A_2$ can only have the inverse words $c^k$ or $c^kb$.
\end{itemize}
\end{rem}

When $a_Q > b_Q + c_Q$ for a reducing sphere $Q$, the following proposition gives one more condition which applies to word lengths of $AA$ arcs.


\begin{prop}\label{prop:two_side_two_letter_aa} 
	Suppose that $a_Q > b_Q + c_Q$ and that $\qsig$ has no $AA$ arc on $\s$ or on $\t$ describing a single letter word in $\pi_1(V)$. Then $\s$ and $\t$, cannot simultaneously contain $AA$ arcs each of which describe a two letter word in $\pi_1(V)$. 
\end{prop}

\begin{proof} Assume that $\qsig$ has no $AA$ arc on $\s$ or on $\t$ describing a single letter word in $\pi_1(V)$ and assume to the contrary that both $\s$ and $\t$ contain an $AA$ arc which describe a two letter word in $\pi_1(V)$. Without loss of generality, assume that $\s$ contains an $AA$ arc with end points on $A_1$ and $A_2$ which describes the word $bc^{-1}$ in the orientation of $\qsig$. In this case $\t$ cannot contain an $AA$ arc describing the same word $bc^{-1}$ or the inverse word $cb^{-1}$ owing to Lemma \ref{lem1-ZAZ-bigon}. So, the $AA$ arc on $\t$ will have to describe the word $b^{-1}c$ or the word $c^{-1}b$ depending on whether it starts on $A_2$ or $A_1$, respectively. We will now show that even this cannot occur. For the purpose of this argument, we note that an $AB$ arc of $\sig$ has to continue as a $BA$ arc. So we tie up such $AB, BA$ arc combinations and call them $ABA$ arcs. We concatenate the words described by the $AB$ and the $BA$ arcs involved and call the concatenated and reduced word in $\pi_1(V)$ as the word of the $ABA$ arc. Likewise we tie up $AC, CA$ arc combinations and call them $ACA$ arcs. We likewise define the words of $ACA$ arcs. Now, $\qsig$ can be seen as a sequence of $AA$, $ABA$ and $ACA$ arcs. Further, the concatenation of the words described by these arcs should reduce to the empty word. Now, we gather the possible words described by these arcs when they are traced from $A_i$ to $A_j$ for $i,j \in \{1,2\}$ in Table \ref{tab:two_letters}. 
	
	In this table \ref{tab:two_letters} an arc starts on $A_1'$ means that an arc starts on $A_1$ and on $\s$, whereas an arc starts on $A_1''$ means that it starts on $A_1$ and on $\t$. An arc ends on $A_1'$ means that the arc ends on $A_1$ and on $\s$ whereas it ends on $A_1''$ means that the arc ends on $A_1$ and on $\t$.
	
	\begin{table}[h!]
		\begin{tabular}{|c|c|c|c|c|}
			\hline
			Ends on $\longrightarrow$ & $A_1'$ & $A_1''$ & $A_2'$ & $A_2''$ \\
			Starts on $\downarrow$ & & & & \\ \hline
			$A_1'$ &  & $c, bc, b, bc$ & $bc^{-1}$ & $\phi, b$\\ \hline
			$A_1''$ & $c^{-1}, b^{-1}, c^{-1}b^{-1}$ &  & $\phi, c^{-1}$ & $c^{-1}b$ \\ \hline
			$A_2'$ & $cb^{-1}$ & $\phi, c$ &  & $c, b, cb$ \\ \hline
			$A_2''$ & $\phi, b^{-1}$ & $b^{-1}c$ & $c^{-1}, b^{-1}c^{-1}, b^{-1}$ &  \\ \hline
		\end{tabular}
		\caption{Words described by the $AA, ABA$, and $ACA$ arcs traced from $A_i$ to $A_j$, $i,j\in\{1,2\}$}\label{tab:two_letters}
	\end{table} 
	
	Figure \ref{fig:arcs-AA=bc-} shows the possible arcs on $\s$ along with the $AA$ arc whose word is $bc^{-1}$. Figure \ref{fig:arcs-AA=cb-} shows the possible arcs on $\t$ along with the $AA$ arc whose word is $b^{-1}c$ or $c^{-1}b$ depending on how it is traced.
	
	\begin{figure}[h!]
	\begin{subfigure}{.49\linewidth}
		\centering
		\includegraphics[height=1.5in]{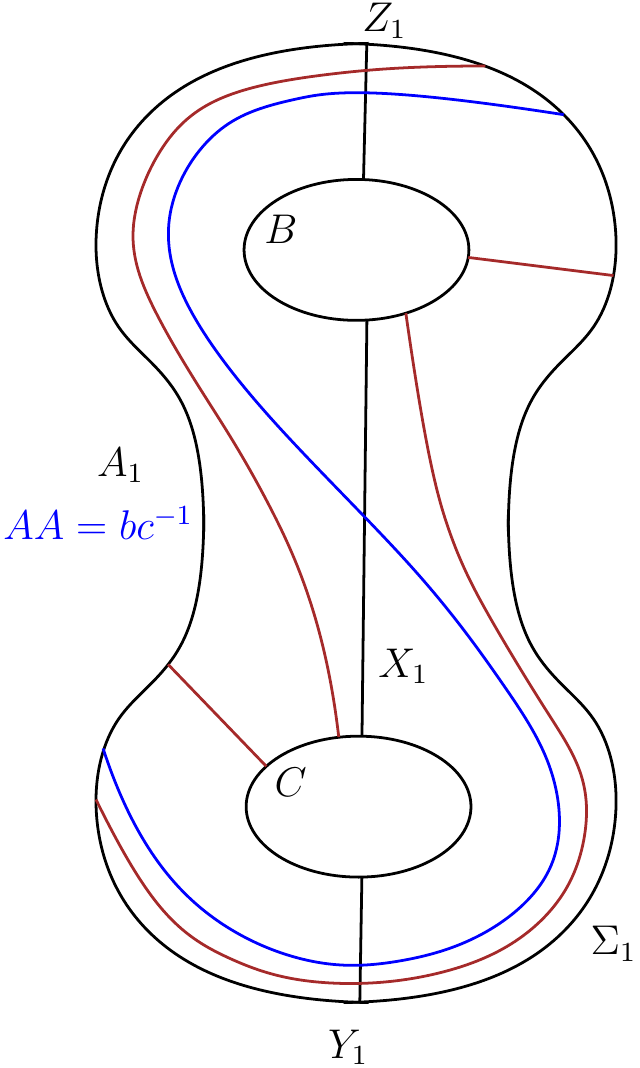}
		\caption{When $AA$ describes $bc^{-1}$}\label{fig:arcs-AA=bc-}
	\end{subfigure}
	\begin{subfigure}{.49\linewidth}
		\centering
		\includegraphics[height=1.5in]{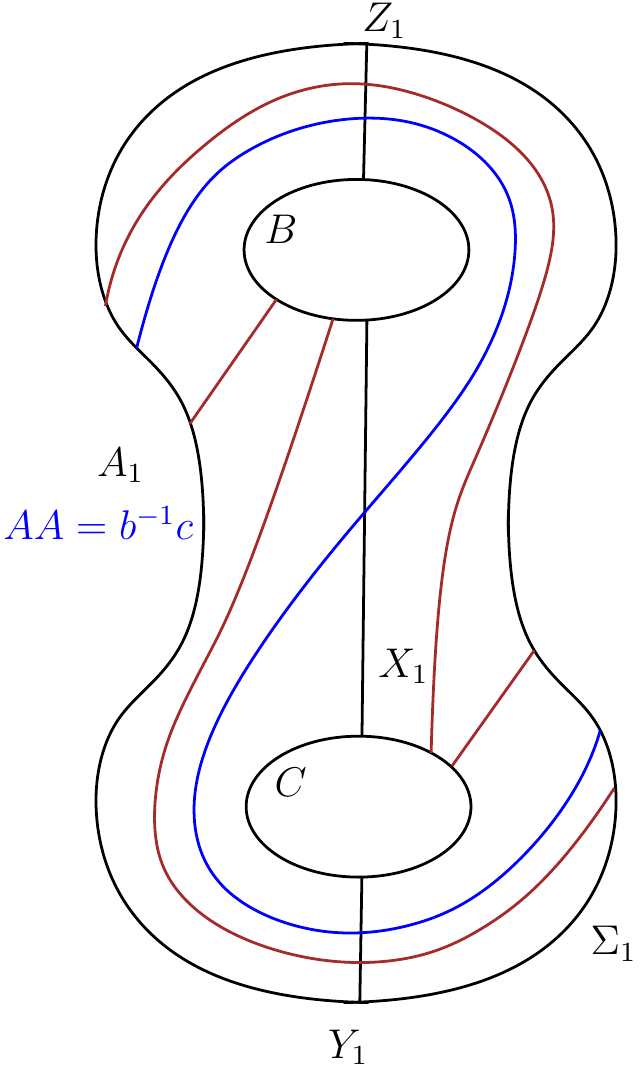}
		\caption{When $AA$ describes $b^{-1}c$}\label{fig:arcs-AA=cb-}
	\end{subfigure}
	\caption{Possible arcs on $\Sigma_1$ with two-lettered $AA$ arcs}\label{fig:arcs-with-AA}
\end{figure}
	
	Let $\kappa_1$ be the $AA$ arc of $\qsig$ whose word is $bc^{-1}$. Starting at $\kappa_1$ and following the orientation of $\kappa_1$, we list the $AA$, $ABA$ and $ACA$ arcs of $\qsig$ as $\kappa_i$'s $1 \leq i \leq p$ for some positive integer $p$. The end point of $\kappa_p$ is the starting point of $\kappa_1$. We first define $w_1$ to be the word of $\kappa_1$. At the $i^{th}$ step, for $1 \leq i \leq p-1$, we define $w_{i+1}$ to be the reduced word corresponding to concatenation of $w_i$ with the word of $\kappa_{i+1}$. By induction on $i$, we show that the length of $w_i$, written $|w_i|$, is always more than one. Since $w_p$ represents the word of $\qsig$ and $w_p$ should have length zero, we have a contradiction. Initially $|w_1| = |bc^{-1}| \geq 2$. $\kappa_1$ ends on $A_2$. So $\kappa_2$ has its initial arc on $\t$ and has the possible words from the row against $A_2''$ in Table \ref{tab:two_letters}. Note that none of the words start with $c$. So the concatenation of $w_1$ with the word of $\kappa_2$ is already reduced and $|w_2| \geq |w_1| \geq 2$. This can be taken as the base case. Now suppose that, for some integer $k, 1 \leq k \leq p-1$, the reduced word $w_k$ is formed and $|w_l| \geq 2$ for all $l, 1 \leq l \leq k$. The arc $\kappa_k$ ends on any of $A_1', A_1'', A_2'$ or $A_2''$ and so the $\kappa_{k+1}$ will start on $A_1'', A_1', A_2''$ or $A_2'$ respectively. In the above table, for any of these combinations, we notice that there is no reduction upon concatenation of the words, except possibly when the word of $\kappa_k$ is empty. So, except possibly in the case that the word of $\kappa_k$ is empty, $|w_{k+1}| \geq 2$. Let us now consider the remaining case when the word of $\kappa_k$ is empty and suppose that $\kappa_{k+1}$ starts on $A_1'$. Then $\kappa_k$ should have ended on $A_1''$. From the Table \ref{tab:two_letters}, there is only one possibility for an arc describing the empty word and ending on $A_1''$, namely an $ACA$ arc which starts on $A_2'$, which in turn means that the $\kappa_{k-1}$ must have ended on $A_2''$. For all arcs ending on $A_2''$ we notice that the ending letter of such arcs is $b$ or $c$, unless it is the empty word. So the possible words for $\kappa_{k-1}$, when $\kappa_k$ has an empty word, have no reduction when concatenated with the word of $\kappa_{k+1}$, unless the word of $\kappa_{k-1}$ is trivial. So unless the $\kappa_{k-1}$ has an empty word, $|w_{k+1}| \geq 2$. If $\kappa_{k-1}$ has an empty word then from the Table \ref{tab:two_letters}, we infer that it must have been an $ABA$ arc starting on $A_1'$. Continuing thus, we see that unless $\kappa_{k+1}$ is preceded by a finite sequence of $ACA$ and $ABA$ arcs which end on $A_1''$ and $A_2''$, respectively, contributing empty words to $w_{k+1}$, we have that $|w_{k+1}| \geq 2$. Even in the case when there are a sequence of $ACA$ and $ABA$ arcs preceeding $\kappa_{k+1}$ which end on $A_1''$ and $A_2''$, respectively, contributing empty words to $w_{k+1}$, this sequence has to end in finite number of steps as the index starts reducing and when such a sequence ends, the arc immediately preceding this sequence of arcs will contribute a word which ends with $b$ or $c$. So there is no reduction on concatenation with the word of $\kappa_{k+1}$. Hence $|w_{k+1}| \geq 2$. The argument when $\kappa_{k+1}$ starts on $ A_2', A_1''$ or $A_2''$ is similar to the case when $\kappa_{k+1}$ starts on $A_1'$. So, even in the case that the word of $\kappa_k$ is empty, we see that $|w_{k+1}| \geq 2$ and so by induction, we conclude that $|w_{l}| \geq 2$ for all $l \in \mathbb{N}$ when the indexing is taken modulo $p$. This contradicts the fact that $\qsig$ must have a trivial word and so $|w_p| = 0$. This completes the proof of this theorem.
\end{proof}

Using Proposition \ref{prop:two_side_two_letter_aa}, we conclude that when $a_Q > b_Q + c_Q$, for a reducing sphere $Q$, and when $Q \neq P$, the word length of $AA$ arcs on $\s$ and $\t$ can have the following possibilities: 
\begin{itemize}
\item[Case (i):] one of $\s$ or $\t$ have an $AA$-arc describing a single letter word
\item[Case (ii):] neither $\s$ nor $\t$ have an $AA$-arc describing a single letter word and both $\s$ and $\t$ have an $AA$-arc with an empty word, and 
\item[Case (iii):] neither $\s$ nor $\t$ have an $AA$-arc describing a single letter word and one of $\s$ or $\t$ has an $AA$-arc describing a two-letter word and the other has an $AA$-arc describing an empty word.
\end{itemize}


For Case (ii) we can choose an isotopy so that we can convert this case into Case (i) as follows: $\qsig$ has to intersect $Y$ or $Z$, else by cutting along $Y \cup Z$, we notice that $\qsig$ is a separating curve on a four-boundered sphere and either piece separated by $\qsig$ cannot be a torus with one boundary component, which is a contradiction to the fact that $Q$ is a reducing sphere. So suppose $\qsig$ intersects $Z$ without loss of generality. Since $AA$ arcs have empty words, an $AB$ or a $BA$ arc must intersect $Z$ in an essential way. A finite sequence of $AB, BA$ arcs connect this intersecting $AB$ or a $BA$ arc to an $AA$ arc. Sliding the curve $\qsig$ past $Z$ so that the intersection is shifted onto the first $AA$ arc gives the required isotopy. See Figure \ref{fig:2to1-isotopy}.

\begin{figure}[h!]
	\centering
	\includegraphics[height=2in]{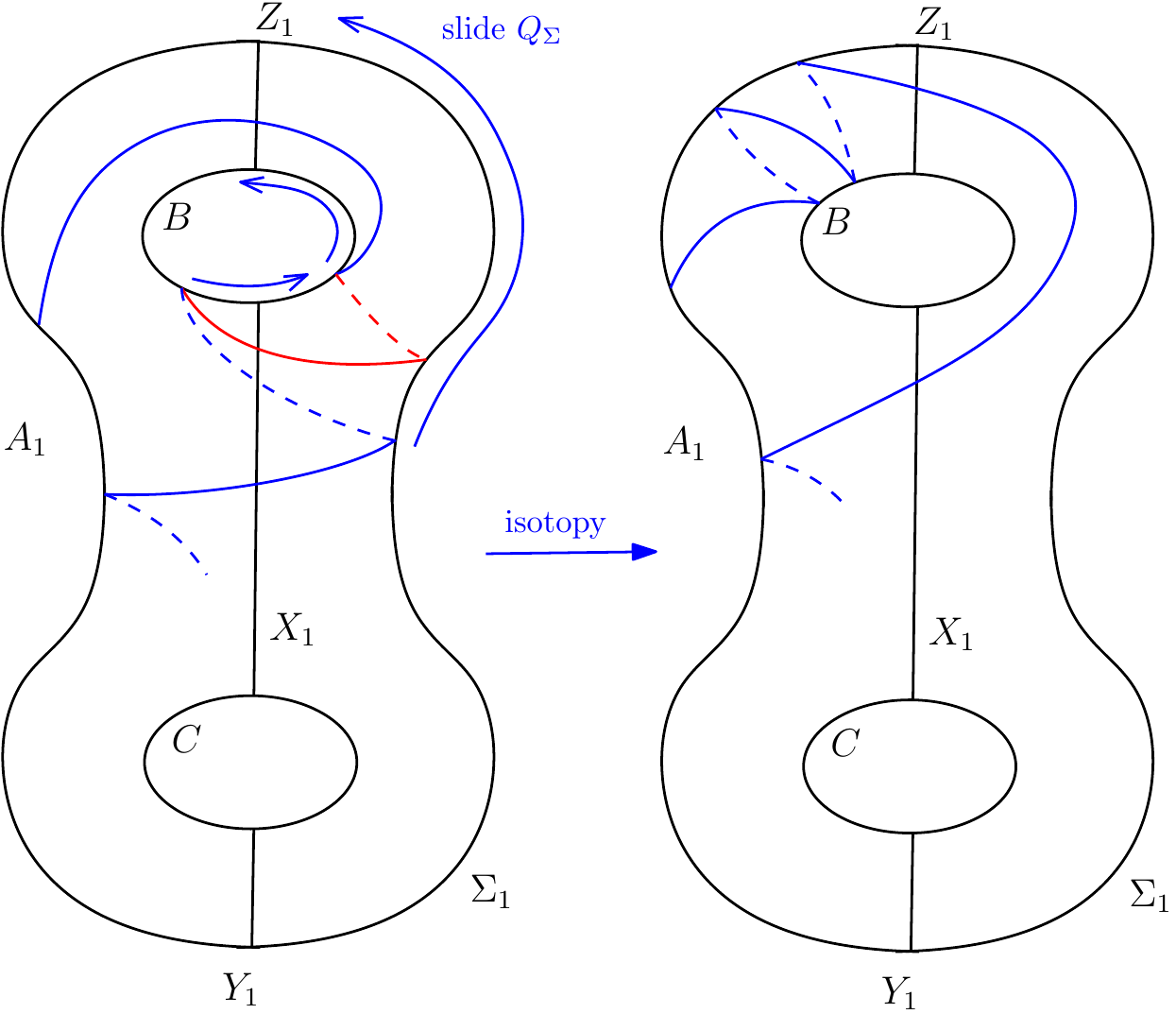}
	\caption{Sliding $\qsig$ past $Z$ making empty-word $AA$ to single lettered $AA$}\label{fig:2to1-isotopy}
\end{figure}

Even for Case (iii) we can choose an isotopy to convert this case into Case (i). Without loss of generality, suppose that $\s$ has an $AA$ arc describing the word $bc^{-1}$. Choose such an $AA$ arc, $\kappa$, which has the topmost $z$ coordinate on $A_1$. Now any arc of $\qsig$ on $\s$ having a $z$-coordinate above the $z$-coordinate of $\kappa$, cannot end on $A_2$. This can be inferred by looking at the component hexagon of $\s$ containing $A_2$ as shown in Figure \ref{fig:3to1-isotopy}. So for any such arc, including $\kappa$, the order of intersection on $\s$ of first with $Z$ and then with $A_1$ can be swapped as intersection with $A_2$ and then an intersection with $Z$ on $\t$ as shown in Figure \ref{fig:3to1-isotopy}.

\begin{figure}[h!]
	\centering
	\includegraphics[height=2in]{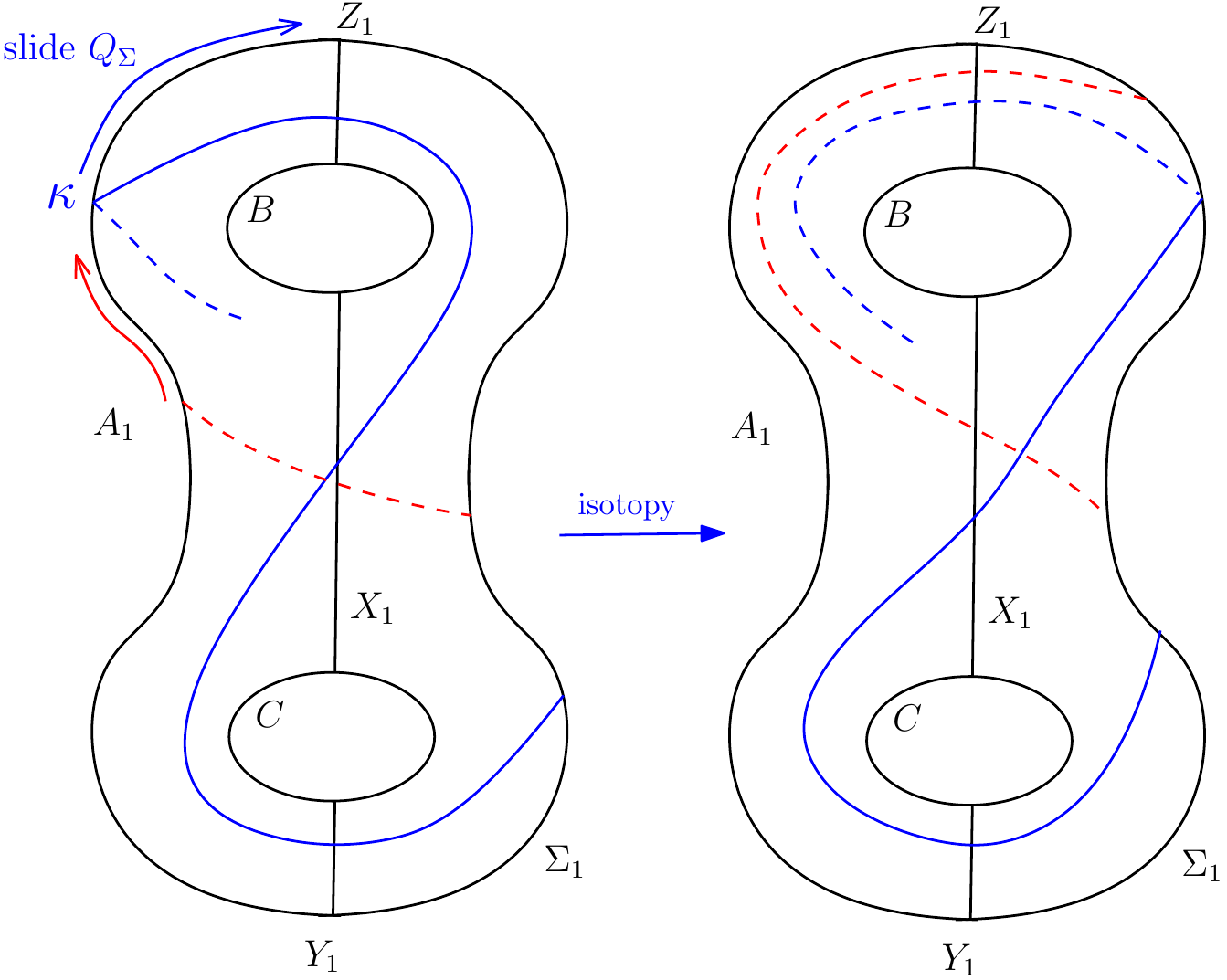}
	\caption{Isotopy taking two lettered $AA$ to single lettered $AA$}\label{fig:3to1-isotopy}
\end{figure}

Using these isotopies, whenever $a_Q > b_Q + c_Q$ for a reducing sphere $Q \neq P$, we can always assume that there is an $AA$ arc of $\qsig$ on $\s$ or on $\t$ which describes a single letter word in $\pi_1(V)$.


\section{Lowering Intersections Using $\beta$}\label{section:beta_lowering}

In this section we study the action of the automorphism $\beta$, which is as described in the section \ref{section:intro}. If $Q$ is a reducing sphere with $Q \neq P$ and $a_Q > b_Q + c_Q$, we show that an application of one of $\beta^{-1}$ or $\beta$ to $Q$ decreases its geometric intersection number with $A$.  We do so by showing that the geometric intersection number of $\qsig$ with one of the curves $\eta_1 := \beta(A)$ or $\eta_2 := \beta^{-1}(A)$ shown in Figure \ref{fig:eta1-eta2} is lesser than that with $A$.

\begin{figure}[h!]
	\centering
	\includegraphics[height=2in]{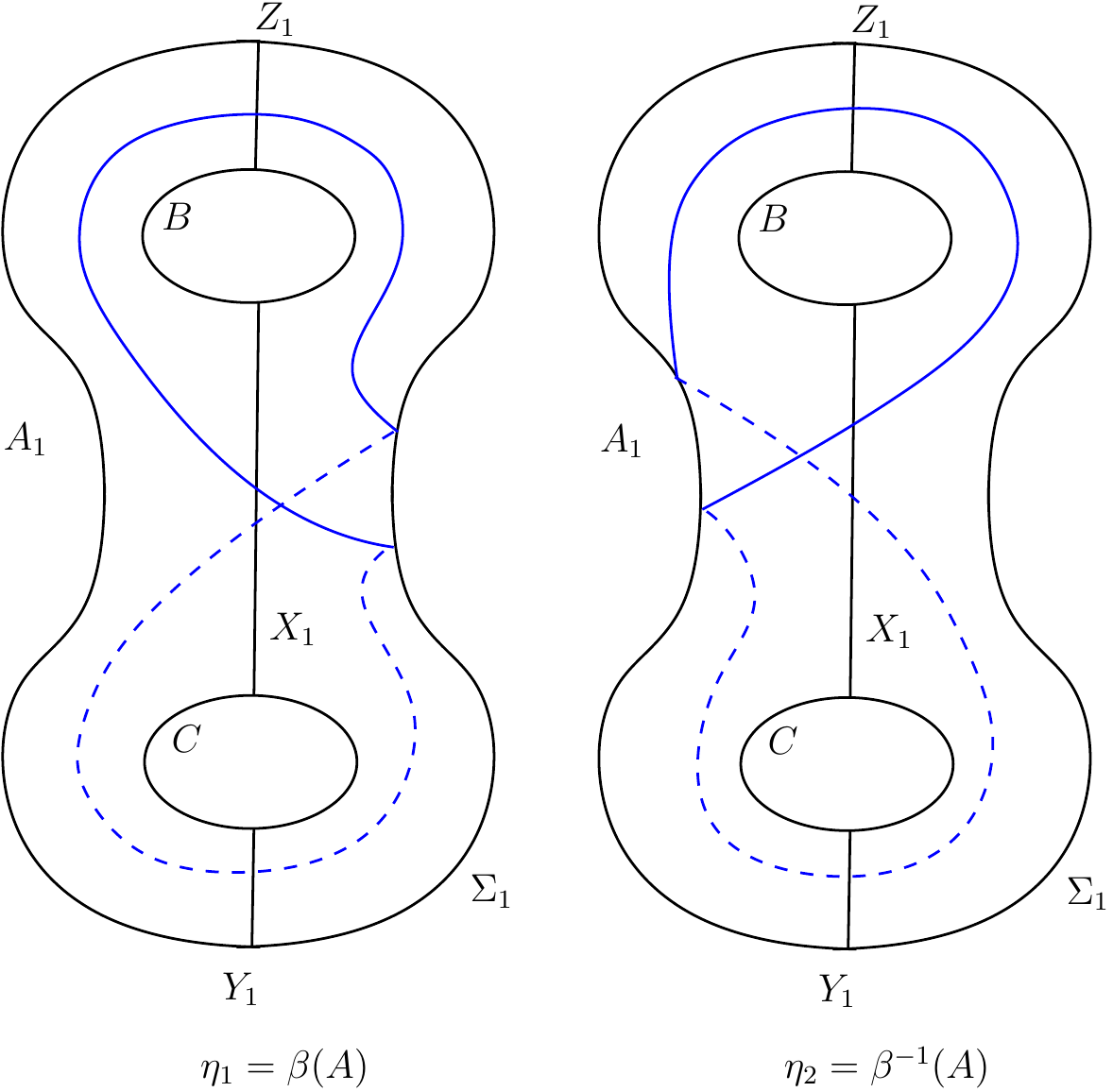}
	\caption{Curves $\eta_1=\beta(A)$ and $\eta_2=\beta^{-1}(A)$}\label{fig:eta1-eta2}
\end{figure}

\begin{thm}\label{thm:beta_single_letter}
	If $a_Q > b_Q + c_Q$  for a reducing sphere $Q$ then either $|\beta(\qsig)\cap A|<|\qsig\cap A|$ or $|\beta^{-1}(\qsig)\cap A|<|\qsig \cap A|$.
\end{thm}

Without loss of generality, as remarked at the end of the previous section, we can assume that $\qsig$ contains an $AA$-arc on $\s$ or on $\t$ which describes a single letter word $b, b^{-1}, c$ or $c^{-1}$ in $\pi_1(V)$. Before we prove this Theorem, we need a few results which give some conditions on $AA, AB, BA, AC$ and $CA$ arcs.

Let us assume that $\qsig$ has an $AA$ arc on $\s$, call it $\chi$, which has both its endpoints on $A_1$ and its word is $b$ or $b^{-1}$ based on its orientation. Every $AA$ arc on $\s$ which satisfies the description of $\chi$, \textit{i.e.} which starts and ends on $A_1$ and has the word $b$ or $b^{-1}$ will be called an $AA$ arc \textit{parallel to $\chi$}.

Let $\omega_1$ and $\omega_2$ be distinct $AA$ arcs of $\qsig$ on $\s$ parallel to $\chi$. Let the $z$-coordinates of the endpoints of $\omega_1$ be $z_{11}$ and $z_{12}$ with $z_{11} < z_{12}$ and the $z$-coordinates of the endpoints of $\omega_2$ be $z_{21}$ and $z_{22}$ with $z_{21} < z_{22}$.  Then either $z_{11} < z_{21} < z_{22} < z_{12}$ or $z_{21} < z_{11} < z_{12} < z_{22}$. This is because, if the endpoints of $\omega_1$ and $\omega_2$ alternate on $A_1$, then $\omega_1$ would intersect $\omega_2$ on $\s$. If $z_{11} < z_{21} < z_{22} < z_{12}$ then we say that $\omega_2$ is nested inside $\omega_1$. Or if $z_{21} < z_{11} < z_{12} < z_{22}$, then we say that $\omega_1$ is nested inside $\omega_2$.

Since $\chi$ cuts an annulus $S_{\chi}$ out of $\s$ containing the circle $B$, all the $A$-ends of the $AB, BA$ arcs of $\qsig$ on $\s$ have to lie on $A_1$ and between the endpoints of $\chi$. With the $z$-coordinate as height we classify the points of $\qsig \cap A_1$ into five stacks. The \textit{first stack of points} with the largest $z$-coordinates consists of $A$-ends of arc segments of $\qsig$ which lie outside $S_{\chi}$, if any, which connect $Z_1$ to $A_1$ and which are not end points of $AA$ arcs on $\s$ parallel to $\chi$. Every $AA$ arc on $\s$ which is parallel to $\chi$ has two ends, one with a higher $z$-coordinate and one with a lower $z$-coordinate. The \textit{second stack of points} consists of those $A$-ends of $AA$ arcs on $\s$ parallel to $\chi$, which have higher $z$-coordinate than their counterparts. The \textit{fourth stack of points} consists of those $A$-ends of $AA$ arcs on $\s$ parallel to $\chi$ which have lower $z$-coordinate than their counterparts. The \textit{third stack of points} consists of the $A$-ends of $AB, BA$ arcs of $\qsig$, if any on $\s$. The \textit{fifth stack of points} consists of the $A$-ends of the arc segments of $\qsig$ on $\s$, if any, which have a $z$-coordinate lower than that of any point in the fourth stack of points and are not the end-points of $AA$ arcs parallel to $\chi$. See figure \ref{fig:A1-stacks}.

\begin{figure}[h!]
	\centering
	\includegraphics[width=.5\linewidth]{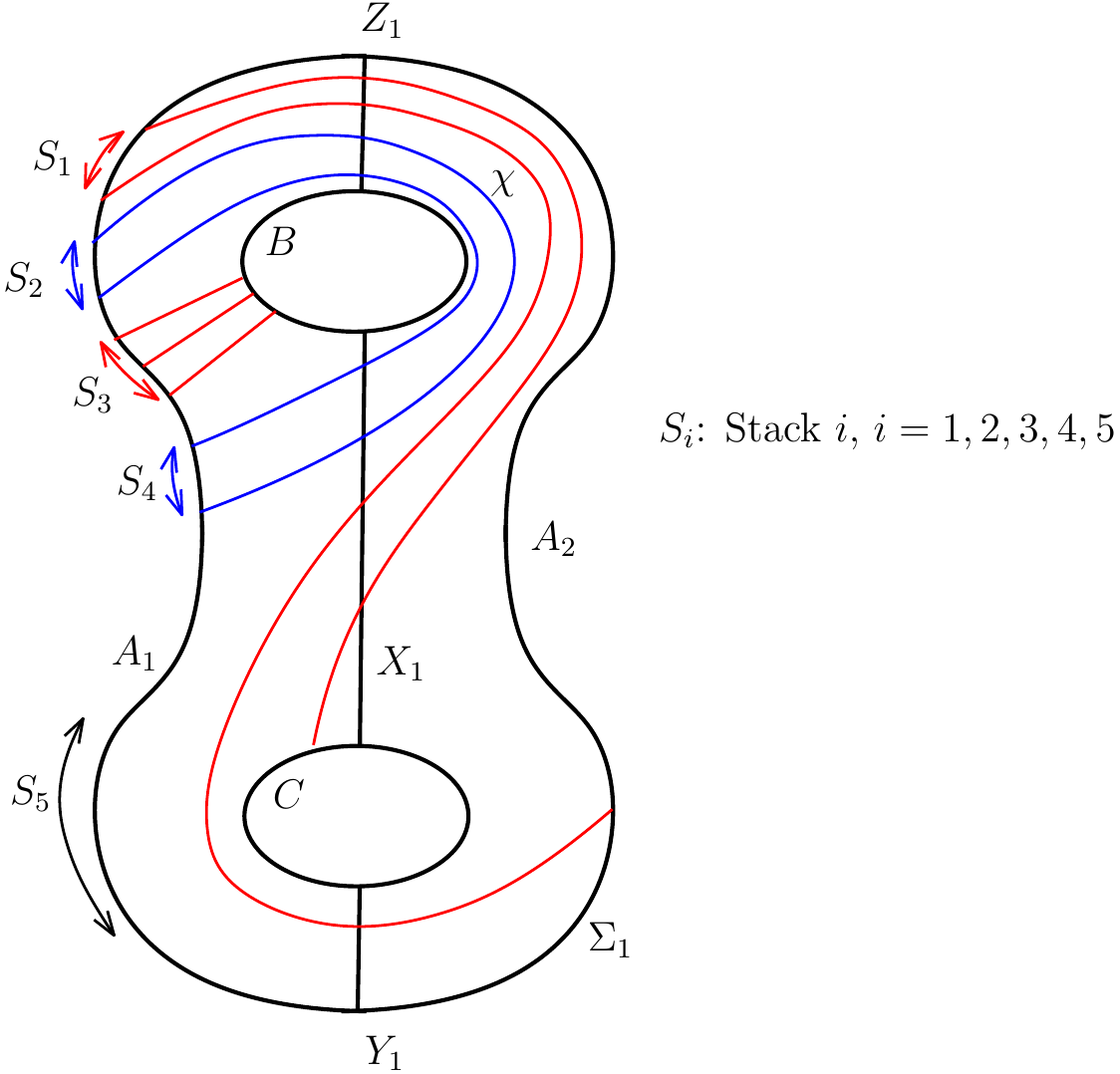}
	\caption{Stacks of points from $\qsig\cap A$ on $A_1$}\label{fig:A1-stacks}
\end{figure}

Now without loss of generality, we can assume that $\chi$ is the innermost among all these parallel $AA$-arcs on $\s$ having word $b$ or $b^{-1}$.

Under the above circumstances, we note the following conditions on $AA$, $AB, BA, AC$ and $CA$ arcs of $\qsig$ on $\t$.

\begin{prop}\label{prop:aa_arc}
	$\qsig$ cannot have an $AA$-arc on $\t$ with both ends on $A_1$, whose word is either $b$ or $b^{-1}$. Also $\qsig$ cannot have an $AA$-arc on $\t$ with word $bc^{-1}$ or $cb^{-1}$. 
\end{prop}

\begin{proof}
	The presence of such an arc (see figure \ref{fig:impossible-AA}) on $\t$ along with $\chi$ contradicts Lemma \ref{lem1-ZAZ-bigon}.
\end{proof}

\begin{figure}[h!]
	\begin{subfigure}{.49\linewidth}
		\centering
		\includegraphics[height=2in]{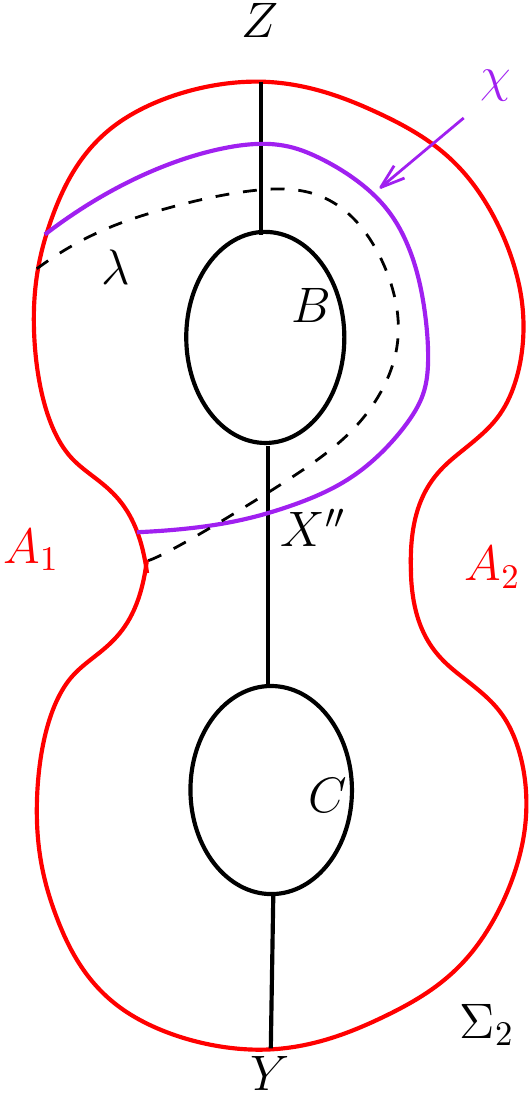}
		\caption{}
	\end{subfigure}
	\begin{subfigure}{.49\linewidth}
		\centering
		\includegraphics[height=2in]{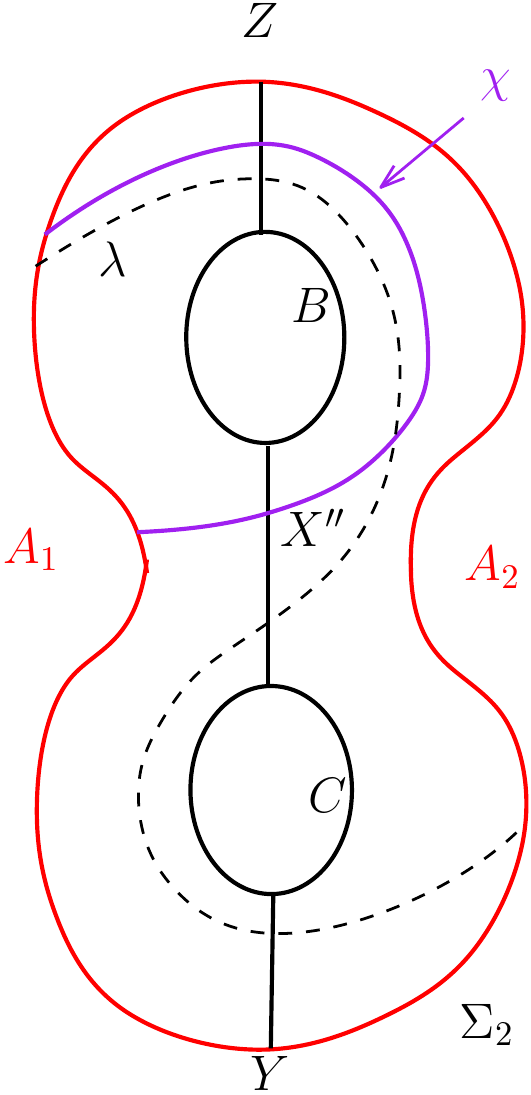}
		\caption{}
	\end{subfigure}
	\caption{Impossible $AA$ arcs on $\Sigma_2$ with $\chi$ on $A_1$}\label{fig:impossible-AA}
\end{figure}

\begin{prop}\label{prop:ab_b_power_k}
	On $\t$, if $\qsig$ has an arc with word $b^k, k \in \mathbb{Z}$ such that one of its endpoints is on $B$ and another is on $A_1$, then its $A_1$-end should lie in stack one, two or three barring the point in stack three with the least $z$-coordinate.
\end{prop}

\begin{proof}
	If possible, suppose $\qsig$ has an arc $\lambda$ on $\t$ with word $b^k, k \in \mathbb{Z}$ such that one of its endpoints is on $B$ and another is on $A_1$, where its $A_1$-end in stack four or five or is the point in stack three with the least $z$-coordinate. Then all arcs on $\t$ having their $A_1$ ends on stacks one, two and three must end on $B$. Including the $A_1$- end point of $\chi$ in stack one and the $A_1$- end point of $\lambda$, there are at least $b_Q +1$ points in stack one, two and three on $A_1$ which are to be connected via arcs on $\t$ to $b_Q$ points on $B$, which is impossible. See Figure \ref{fig:impossible-AB1}. So, such an arc $\lambda$ cannot exist.

\begin{figure}[h!]
	\begin{subfigure}{0.32\linewidth}
		\centering
		\includegraphics[height=2in]{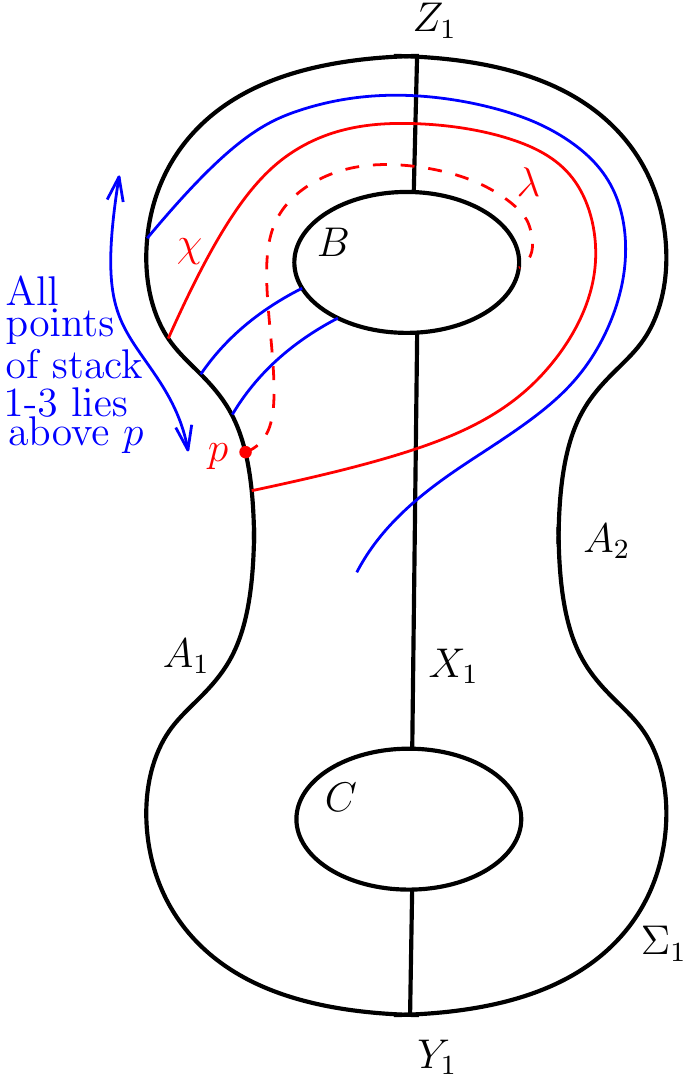}
		\caption{$\lambda$ from $A_1$ with $b^k$}\label{fig:impossible-AB1}
	\end{subfigure}
	\begin{subfigure}{0.32\linewidth}
		\centering
		\includegraphics[height=2in]{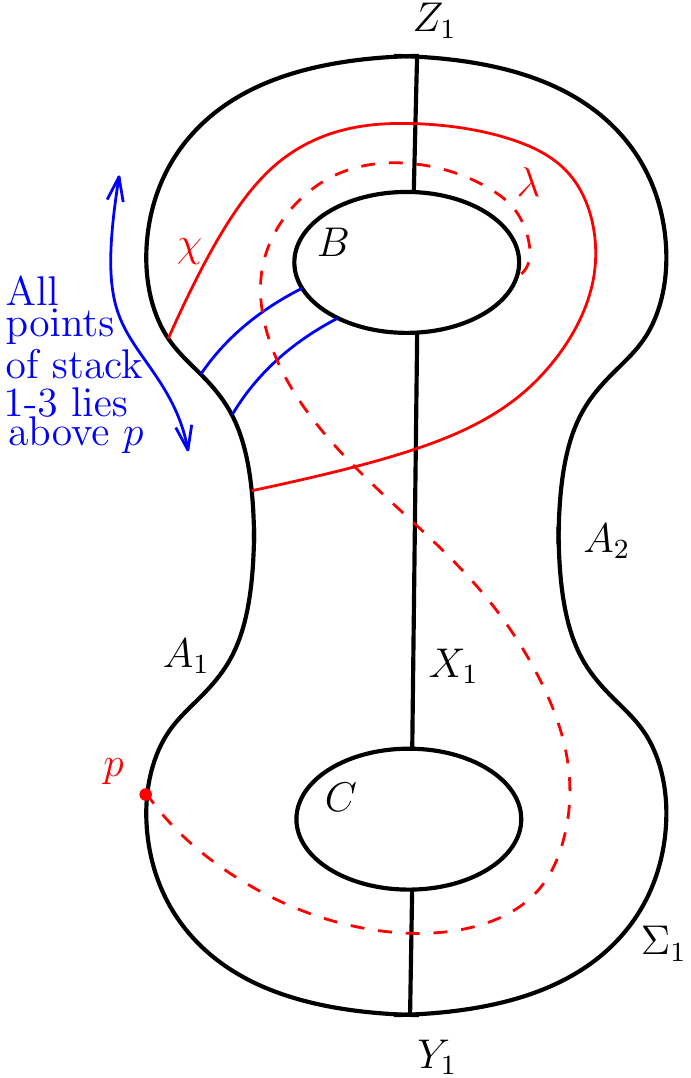}
		\caption{$\lambda$ from $A_1$ with $c^{-1}b^k$}\label{fig:impossible-AB2}
	\end{subfigure}
	\begin{subfigure}{0.32\linewidth}
		\centering
		\includegraphics[height=2in]{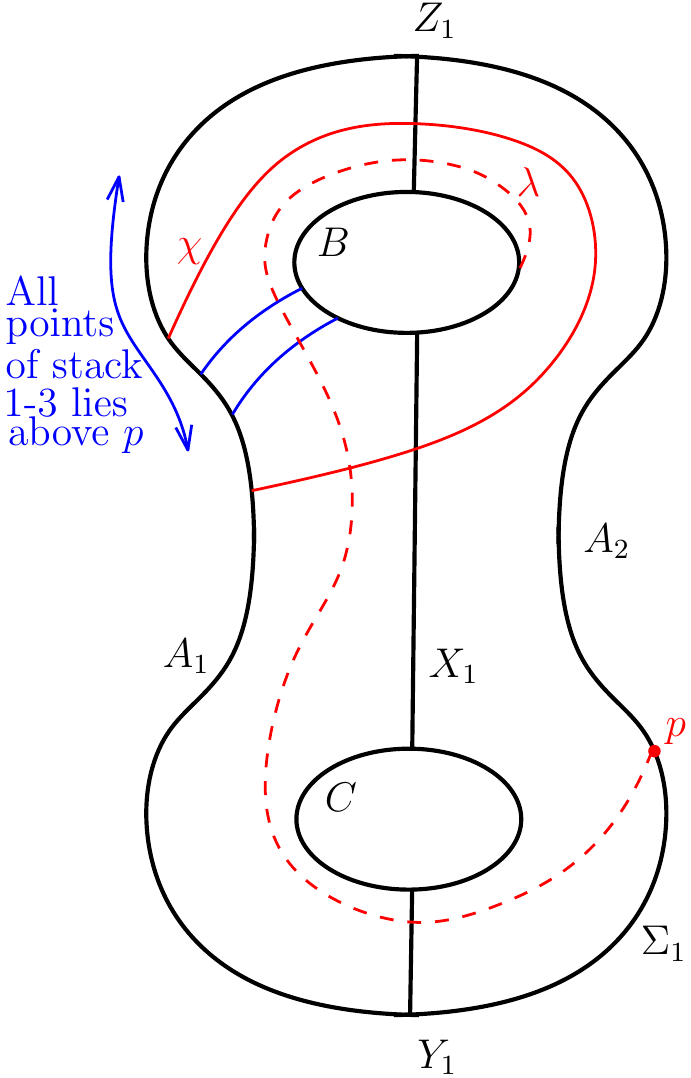}
		\caption{$\lambda$ from $A_2$ with $cb^k$}\label{fig:impossible-AB3}
	\end{subfigure}
	\caption{Impossible $AA$ arcs and $AB$ arcs on $\Sigma_2$ with $\chi$ on $A_1$ (Arc continuation of the top end of $\chi$ must end on $B$. Also all arcs on $\Sigma_2$ with one end in stack 3 must have their other end on $B$.)}\label{fig:impossible-AB}
\end{figure}

\end{proof}

\begin{prop}\label{prop:ab_cinv_b_power_k} If $\t$ contains an arc of $\qsig$ with endpoints on $A_1$ and $B$ with the word $c^{-1}b^k$ then such an $A_1$-end should be in stack three, four or five, barring the topmost point in the third stack with the largest $z$-coordinate.
\end{prop}
\begin{proof}
	Suppose to the contrary that $\lambda$ is such an $AB$-arc on $\t$ with an endpoint on $A_1$ with word $c^{-1}b^k$ and with its endpoint in stack one or two of $A_1$ or the topmost point in stack three on $A_1$. See Figure \ref{fig:impossible-AB2}. Then $\lambda$ will force every arc of $c_Q$ on $\t$ starting from stacks three, four and five to end on $B$ due to Theorem \ref{thm:aa-atmost-2}. Including the $A_1$- end point of $\chi$ on stack four and the $A_1$-  end point of $\lambda$, there are at least $b_Q +1$ points in stack three four and five on $A_1$ which are to be connected via arcs on $\t$ to $b_Q$ points on $B$, which is impossible. Hence such a $\lambda$ cannot exist. 
\end{proof}

\begin{prop}\label{prop:no_ab_cb_power_k} If an arc of $\qsig$ on $\t$ has its end points on $A_2$ and $B$, then its word has to be of the form $b^k$ for some integer $k$. 
\end{prop}
\begin{proof}
	If possible suppose $\lambda$ is an arc on $\t$ with word $cb^k$ and having its end points on $A_2$ and $B$. See Figure \ref{fig:impossible-AB3}. Note that by Theorem \ref{thm:aa-atmost-2}, this is the only other option for such an arc apart from having a word of the form $b^k$ for some integer $k$.
	
	Due to Lemma \ref{lem1-ZAZ-bigon}, such a $\lambda$ will force all arcs on $\t$ having one end-point on $A_1$ to start with a word $c^{-1}$ or end with $c$, as the case maybe, and they must have the other end-point on $B$. But, including the $A_1$-end-points of $\chi$, there are at-least $b_Q+2$ points on $A_1$ whereas only $b_Q$ points on $B$. Therefore, such a $\lambda$ cannot exist.
\end{proof}

\begin{prop}\label{prop:no_ac_bc_power_k} $\qsig$ cannot have an $AC$-arc on $\t$ with word $bc^k$, where $k \in \mathbb{Z}$
\end{prop}
\begin{proof}
	The presence of such an arc on $\t$ along with $\chi$ contradicts Lemma \ref{lem1-ZAZ-bigon}.
\end{proof}

We now define $AA$-arcs, $\eta'$ and $\zeta'$ on $\s$ which are not arcs of $\qsig$ and are disjoint from arcs of $\qsig$ on $\s$ and are parallel to $\chi$ as follows. Define $\eta'$ to be an arc whose ends on $A_1$ are nested inside the ends of $\chi$ so that $A$-ends of all $AB, BA$-arcs on $\s$, if any, are between the ends of $\eta'$. Likewise, define $\zeta'$ to be an arc whose ends on $A_1$ are such that (i) every arc of $\qsig$ on $\s$ parallel to $\chi$ is nested inside $\zeta'$, (ii) the $A$-ends of every $AB, BA$-arc on $\s$ is also between the ends of $\zeta'$ and (iii) the $A$-ends of no other arc of $\qsig$ on $\s$ is between the ends of $\zeta'$. Figure \ref{fig:eta-zeta} shows $\eta'$ and $\zeta'$.

\begin{figure}
	\centering
	\includegraphics[width=.32\linewidth]{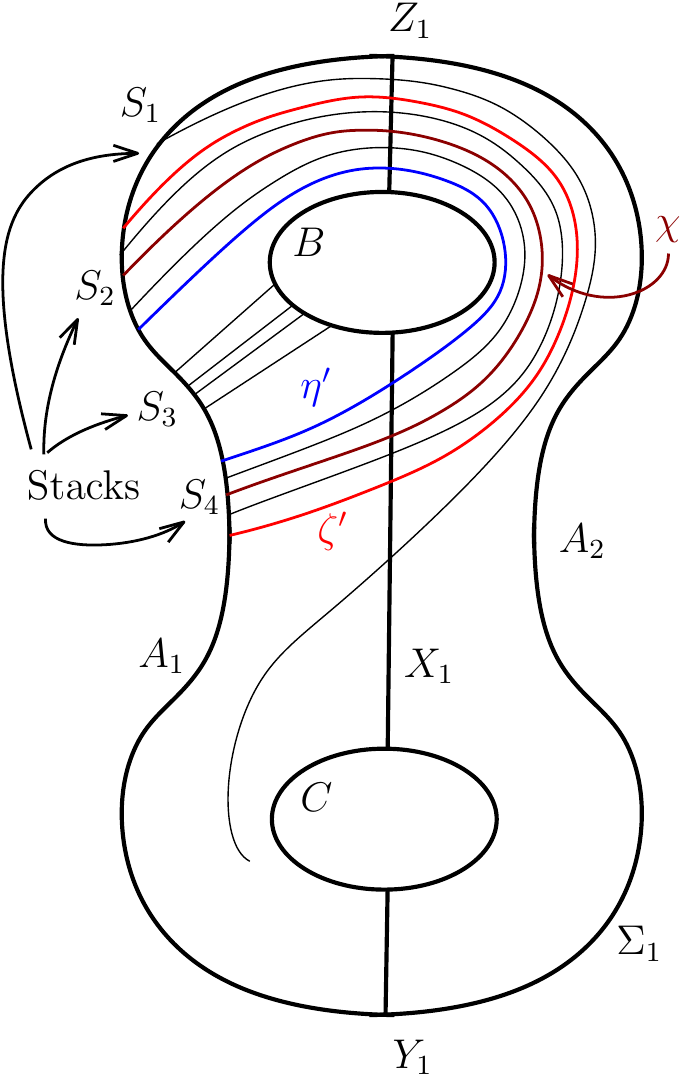}
	\caption{Construction of $\eta'$ and $\zeta'$}\label{fig:eta-zeta}
\end{figure}

Define $\eta''$ and $\zeta''$ to be $AA$ arcs (also not of $\qsig$) on $\t$ such that: (i) $\partial\eta''=\partial\eta'$, $\partial\zeta''=\partial\zeta'$, (ii) the word of $\eta''$ or $\zeta''$ is $c^{-1}$ or $c$ and (iii) $\eta''$ and $\zeta''$ intersect $c_Q$ minimally on $\t$ where the end points of all arcs on $\t$ are rigid \textit{i.e.}cannot be moved while considering their intersection numbers. The simple closed curves $\eta'\cup\eta''$ or $\zeta'\cup\zeta''$ are isotopic to each other on $\sig$ and in turn are isotopic to $A' := \beta^{-1}(A)$. We now show that $|c_Q\cap A'|<|c_Q\cap A|$. By construction, $|c_Q\cap\eta'|=|c_Q\cap\zeta'|=0$. So it is enough to show that either $|c_Q\cap A|>|c_Q\cap\eta''|$ or $|c_Q\cap A|>|c_Q\cap\zeta''|$. We will first show that $|c_Q\cap A|\ge|c_Q\cap\eta''|$ and $|c_Q\cap A|\ge|c_Q\cap\zeta''|$. 

Taking the restrictions on the possible $AA$ arcs of $\qsig$ on $\t$ into account, the intersection of such arcs with $\eta''$ and $\zeta''$ are as shown in Figure \ref{fig:aa-eta-zeta}. Likewise, the intersection of possible $AB, BA$ arcs of $\qsig$ on $\t$ with $\eta''$ and $\zeta''$ are as shown in Figure \ref{fig:ab-eta-zeta}. The intersection of possible $AC, CA$ arcs of $\qsig$ on $\t$ with $\eta''$ and $\zeta''$ are as shown in Figure \ref{fig:ac-eta-zeta}.

\begin{figure}[h!]
	\begin{subfigure}{.24\linewidth}
		\centering
		\includegraphics[height=1.5in]{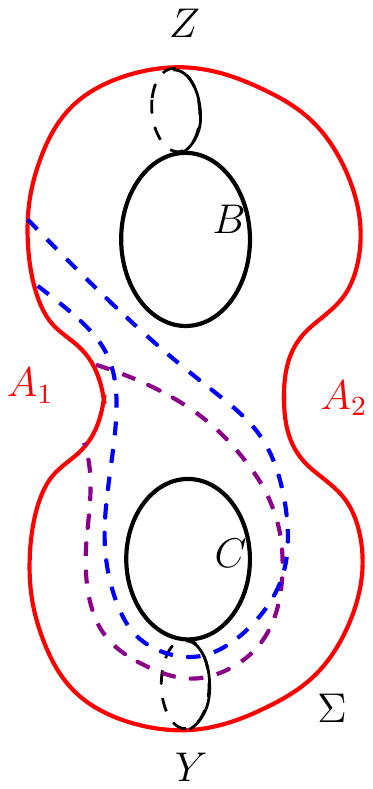}
		\caption{}
	\end{subfigure}
	\begin{subfigure}{.24\linewidth}
		\centering
		\includegraphics[height=1.5in]{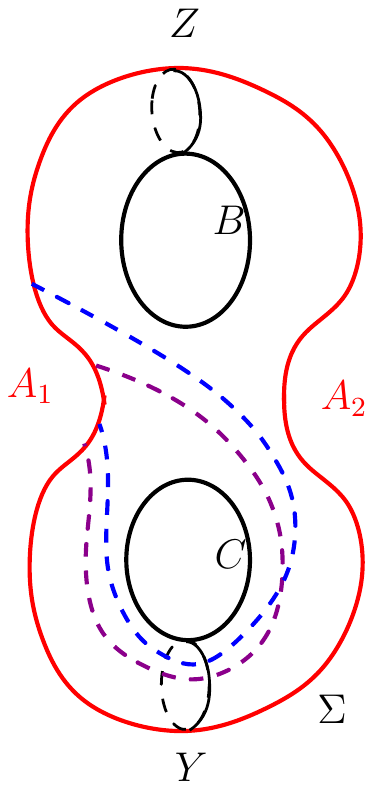}
		\caption{}
	\end{subfigure}
	\begin{subfigure}{.24\linewidth}
		\centering
		\includegraphics[height=1.5in]{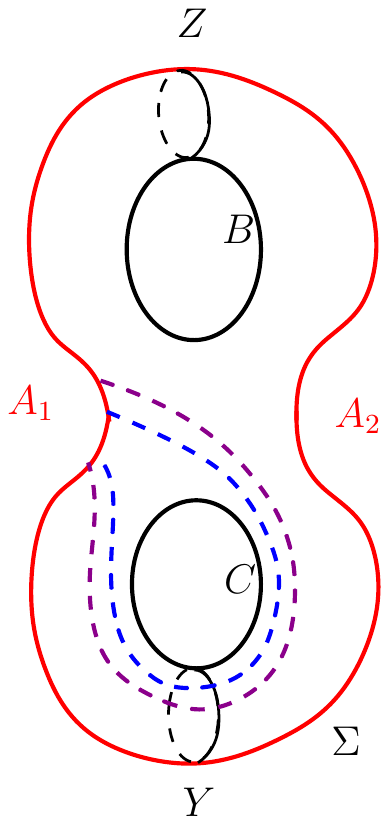}
		\caption{}
	\end{subfigure}
	\begin{subfigure}{.24\linewidth}
		\centering
		\includegraphics[height=1.5in]{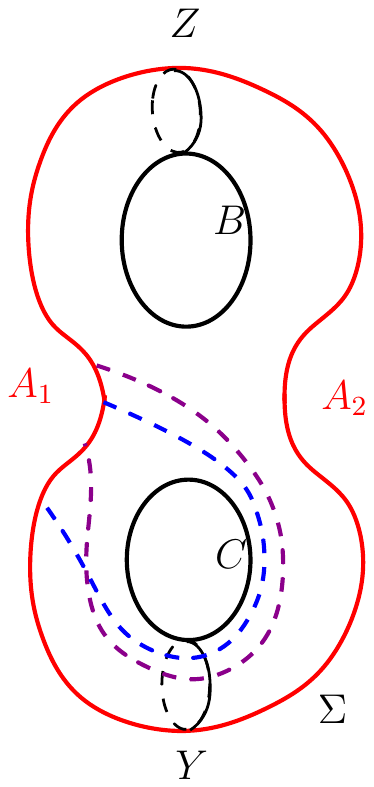}
		\caption{}
	\end{subfigure}

	\begin{subfigure}{.24\linewidth}
		\centering
		\includegraphics[height=1.5in]{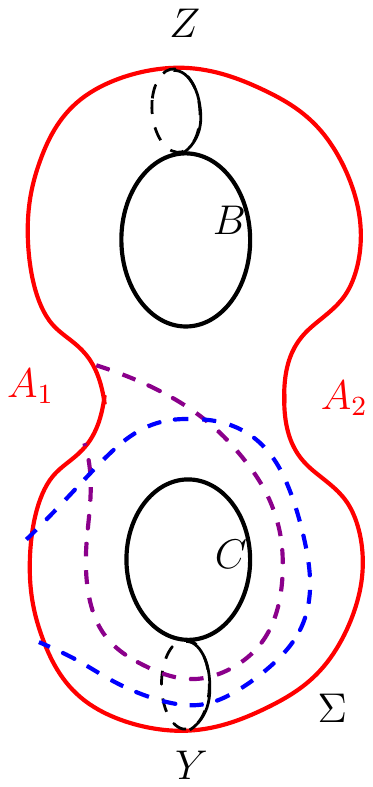}
		\caption{}
	\end{subfigure}
	\begin{subfigure}{.24\linewidth}
		\centering
		\includegraphics[height=1.5in]{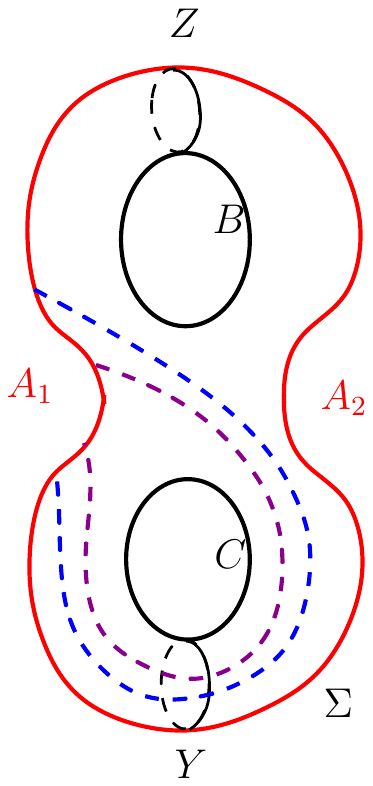}
		\caption{}
	\end{subfigure}
	\begin{subfigure}{.24\linewidth}
		\centering
		\includegraphics[height=1.5in]{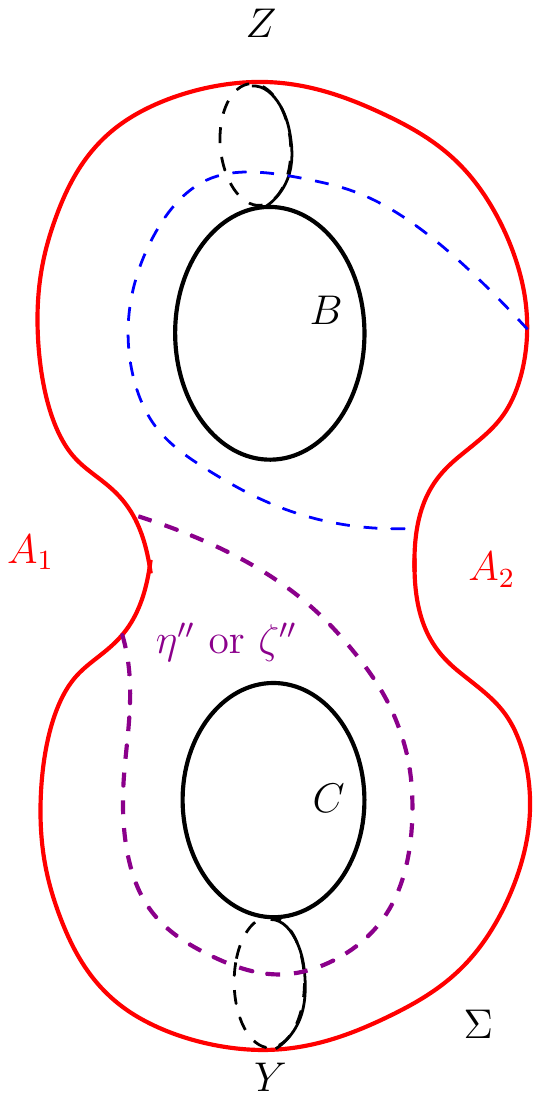}
		\caption{}
	\end{subfigure}
	\begin{subfigure}{.24\linewidth}
		\centering
		\includegraphics[height=1.5in]{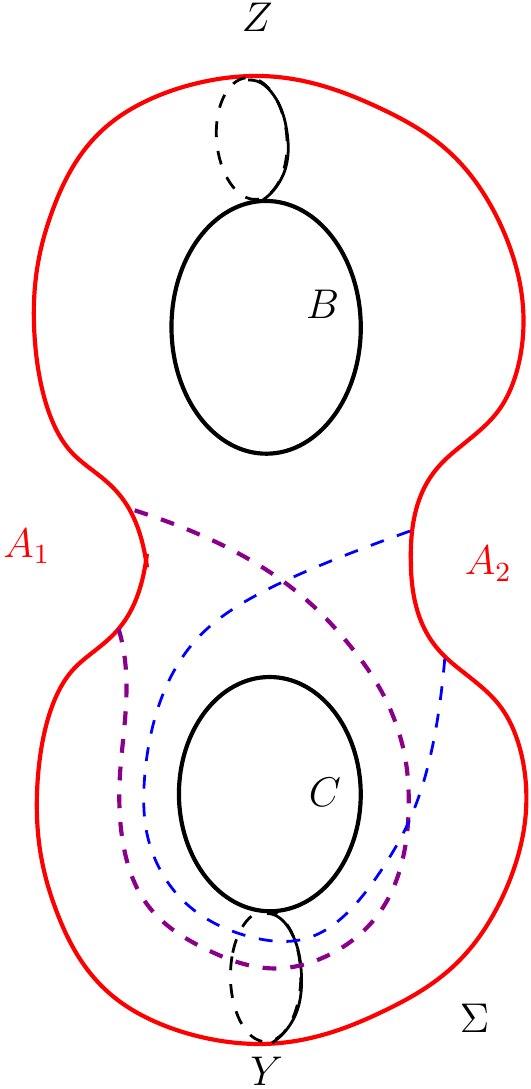}
		\caption{}
	\end{subfigure}

	\caption{Intersection of $AA$ with $\eta''$ and $\zeta''$}\label{fig:aa-eta-zeta}
\end{figure}

\begin{figure}[h!]
		\begin{subfigure}{.3\linewidth}
			\centering
			\includegraphics[height=1.5in]{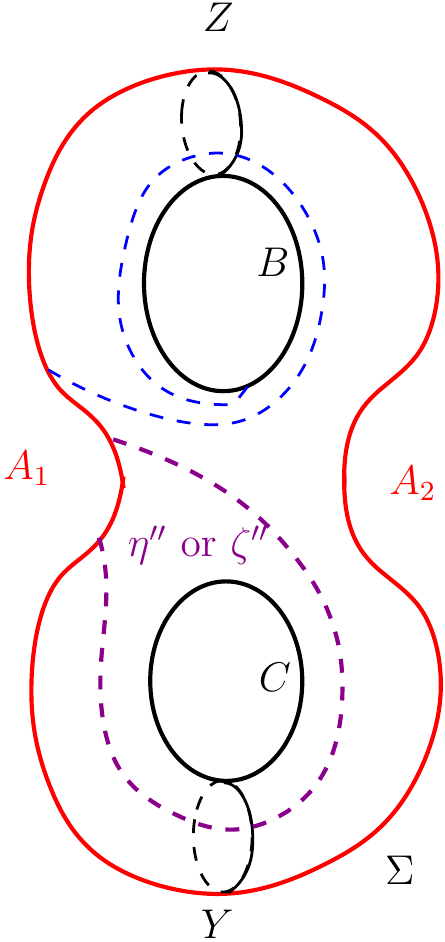}
			\caption{}
		\end{subfigure}
		\begin{subfigure}{.3\linewidth}
			\centering
			\includegraphics[height=1.5in]{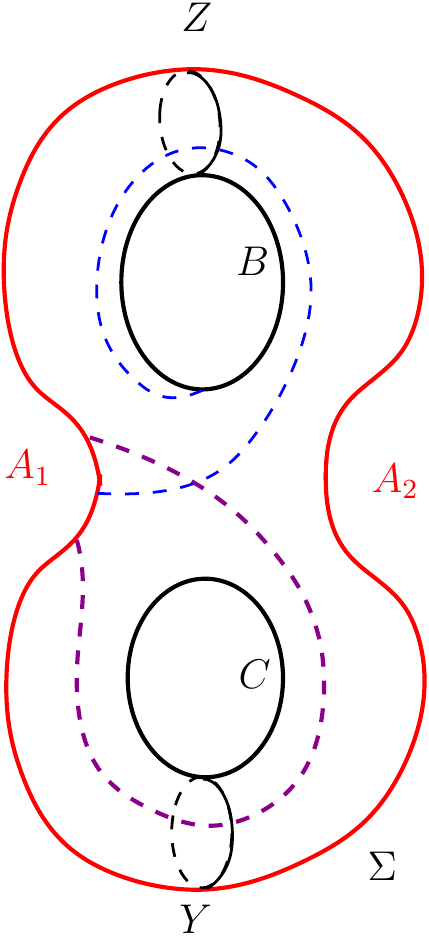}
			\caption{}
		\end{subfigure}
		\begin{subfigure}{.3\linewidth}
			\centering
			\includegraphics[height=1.5in]{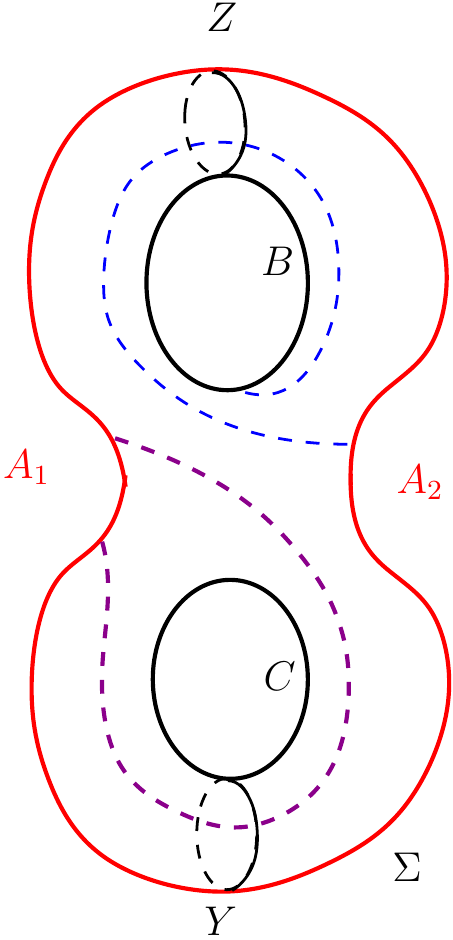}
			\caption{}
		\end{subfigure}
	\caption{Intersection of possible $AB$, $BA$ with $\eta''$ and $\zeta''$} \label{fig:ab-eta-zeta}
\end{figure}

\begin{figure}[h!]
	\begin{subfigure}{.3\linewidth}
		\centering
		\includegraphics[height=1.5in]{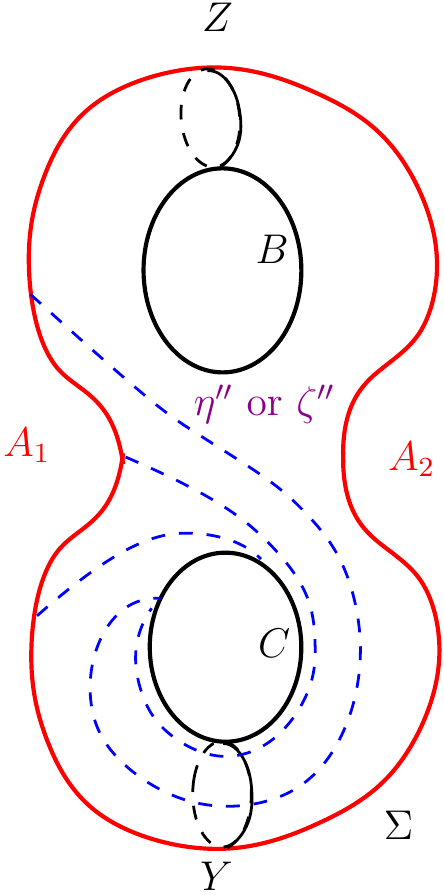}
		\caption{}
	\end{subfigure}
	\begin{subfigure}{.3\linewidth}
		\centering
		\includegraphics[height=1.5in]{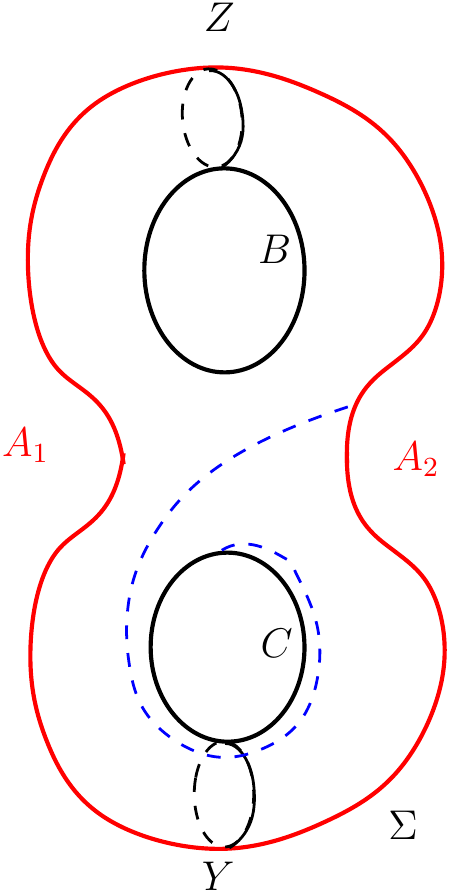}
		\caption{}
	\end{subfigure}
	\begin{subfigure}{.3\linewidth}
		\centering
		\includegraphics[height=1.5in]{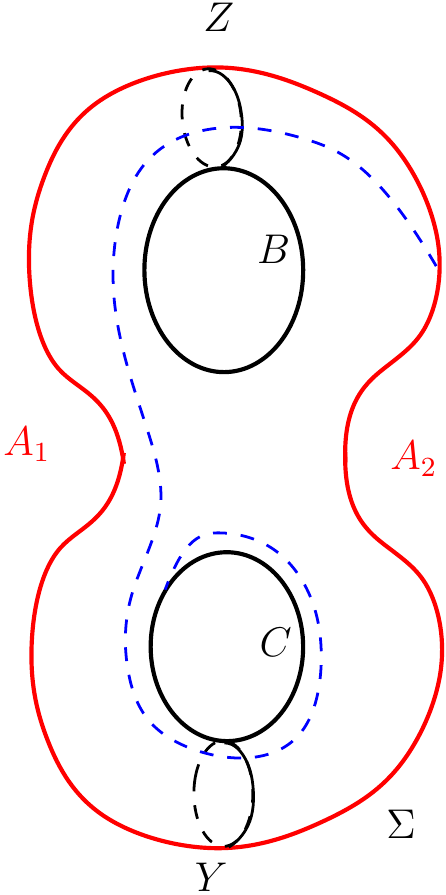}
		\caption{}
	\end{subfigure}
	\caption{Intersection of possible $AC$, $CA$ with $\eta''$ and $\zeta''$} \label{fig:ac-eta-zeta}
\end{figure}


From the above discussion and Figures \ref{fig:aa-eta-zeta}, \ref{fig:ab-eta-zeta} and \ref{fig:ac-eta-zeta}, we conclude the following.

\begin{prop}
	Every arc $\kappa$ of $\qsig$ on $\t$ satisfies $|\kappa \cap \eta''| \leq |\kappa \cap A|$ and $|\kappa \cap \zeta''| \leq |\kappa \cap A|$. 
\end{prop}

Next we show that at least one of the two strict inequalities: $|\qsig \cap A|>|\qsig \cap \eta''|$ or $|\qsig \cap A|>|\qsig \cap \zeta''|$ holds, by showing the following.

\begin{lem}
	There exists an arc $\lambda$ of $\qsig$ on $\t$ such that $|\lambda \cap \eta''|<|\lambda \cap A|$ or $|\lambda \cap \zeta''|<|\lambda \cap A|$.
\end{lem}
\begin{proof}
	Consider the points $\chi \cap A_1$. Let the arc of $\qsig$ on $\t$ with the point in stack two of $\chi \cap A_1$ as an endpoint be $\theta_1$ and the arc with the point in stack four of $\chi \cap A_1$ as an endpoint be $\theta_2$. By proposition \ref{prop:ab_cinv_b_power_k}, $\theta_1$ cannot be an $AB$ arc with the word $c^{-1}b^k$ or a $BA$ arc with the word $b^kc$. We have the following cases.
	\begin{itemize}
		\item[Case 1] Suppose that $\theta_1$ is an $AA$-arc with endpoints on $A_1$ and $A_2$ with the empty word. In this case, $\theta_1$ is disjoint from $\eta''$ and is the required $\lambda$. See Figure \ref{fig:theta1-case1}. 

\begin{figure}[h!]
	\begin{subfigure}{.24\linewidth}
		\centering
		\includegraphics[height=1.5in]{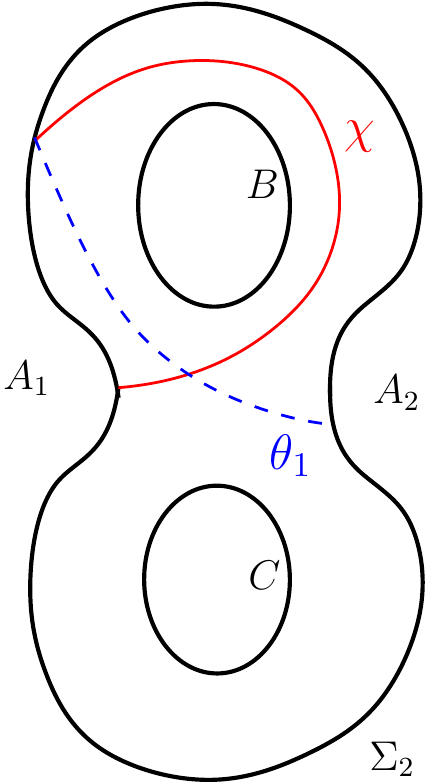}
		\caption{}\label{fig:theta1-case1}
	\end{subfigure}
	\begin{subfigure}{.72\linewidth}
		\centering
		\includegraphics[height=1.5in]{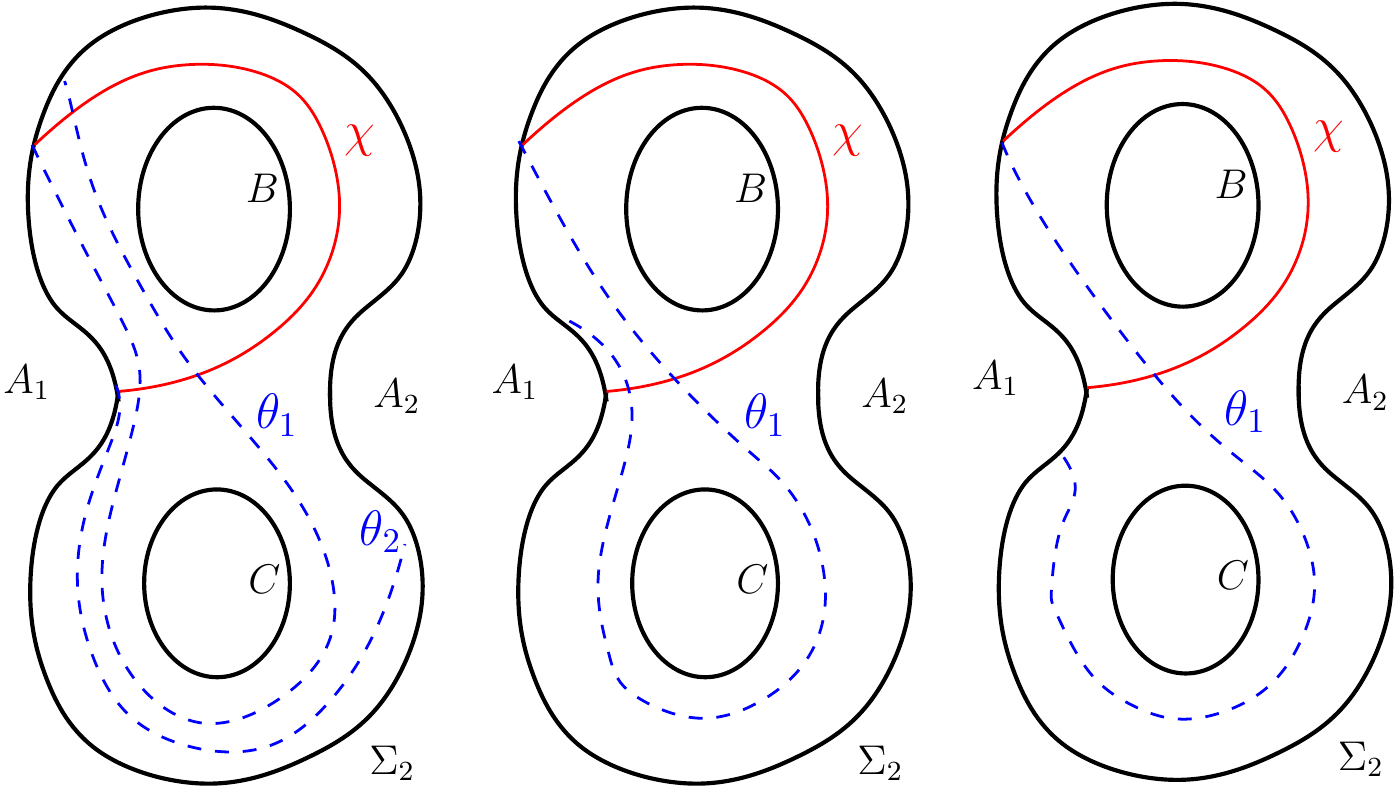}
		\caption{}\label{fig:theta1-case2}	
	\end{subfigure}
	\begin{subfigure}{.24\linewidth}
		\centering
		\includegraphics[height=1.5in]{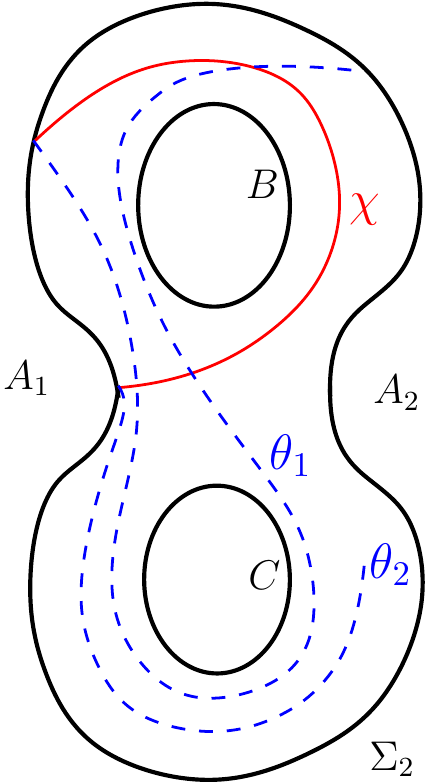}
		\caption{}\label{fig:theta1-case3}	
	\end{subfigure}
	\begin{subfigure}{.24\linewidth}
		\centering
		\includegraphics[height=1.5in]{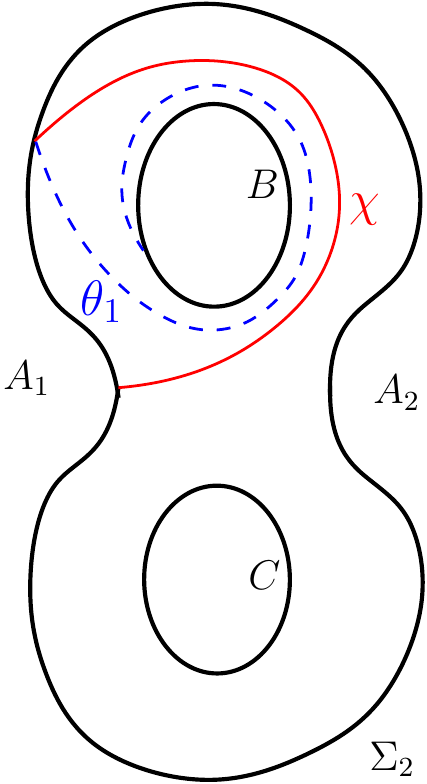}
		\caption{}\label{fig:theta1-case4}	
	\end{subfigure}
	\begin{subfigure}{.24\linewidth}
		\centering
		\includegraphics[height=1.5in]{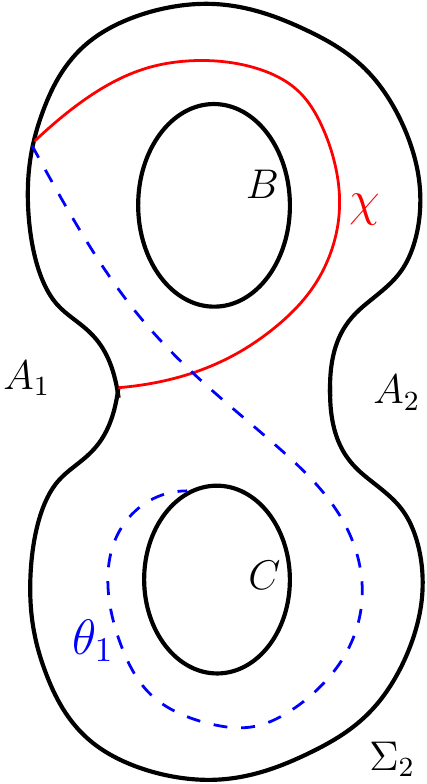}
		\caption{}\label{fig:theta1-case5}	
	\end{subfigure}
	\caption{The arc $\lambda$ of $Q_\sigma$ on $\Sigma_2$}\label{fig:theta1}
\end{figure}

		\item[Case 2] Suppose $\theta_1$ is an $AA$-arc with both its endpoints on $A_1$ with the word $c$ or $c^{-1}$. If the other end point of $\theta_1$ is in stack one or two, then $\theta_1$ cuts off an annulus from $\t$ containing $C$ as one boundary component and so $\theta_2$, cannot be an $AC$ arc. So, if the other end point of $\theta_1$ is in stack one or two then $\theta_2$ is an $AA$ arc with the word $c^{-1}b$, $c$ or $c^{-1}$ or an $AB$ arc with the word $c^{-1}b^k$, $k \in\mathbb{Z}$ (or a $BA$ arc with the inverse word), and in any case, $\theta_2$ is disjoint from $\eta''$ and is the required $\lambda$. If the other end point of $\theta_1$ is in stack three then it intersects $\eta''$ once, whereas $\theta_1$ intersects $A$ twice. So $\theta_1$ is the required $\lambda$.  Finally, if the other end point of $\theta_1$ is in stack four or five then it is disjoint from $\eta''$ and $\theta_1$ is the required $\lambda$. See Figure \ref{fig:theta1-case2}.
		\item[Case 3] If $\theta_1$ is an $AA$ arc with endpoints on $A_1$ and $A_2$ with the word $c^{-1}b$, then $\theta_2$ is either an $AA$ arc with the word $c^{-1}b$ or is an $AB$ arc with the word $c^{-1}b^k$ (or a $BA$ arc with the inverse word). In any case, $\theta_2$ is disjoint from $\eta''$ and is the required $\lambda$. See Figure \ref{fig:theta1-case3}.
		\item[Case 4] If $\theta_1$ is an $AB$ arc with the word $b^k$, $k \in\mathbb{Z}$, then it is disjoint from $\eta''$ and is the required $\lambda$. See Figure \ref{fig:theta1-case4}.
		\item[Case 5] Suppose that $\theta_1$ is an $AC$-arc or a $CA$-arc with the word $c^k$, $k \in \mathbb{Z}$. Then, $\theta_1$ is disjoint from $\zeta''$ and is the required $\lambda$. See Figure \ref{fig:theta1-case5}.
	\end{itemize}
\end{proof}
Owing to the symmetry of $\sig$, all of the above discussion holds even if we replace $\chi$ by any $AA$ arc on $\sig$ with both endpoints on $A_1$ or both endpoints on $A_2$ describing a single letter word in $\pi_1(V)$. For comprehensiveness, Figure \ref{fig:singleAA} shows the eight possibilities for such an $AA$ arc. We note that, in some cases we might have to use $A'' = \beta(A)$ instead of $A' = \beta^{-1}(A)$ throughout the above discussion for the inequalities to hold.

\begin{figure}\centering
\includegraphics[height=3in]{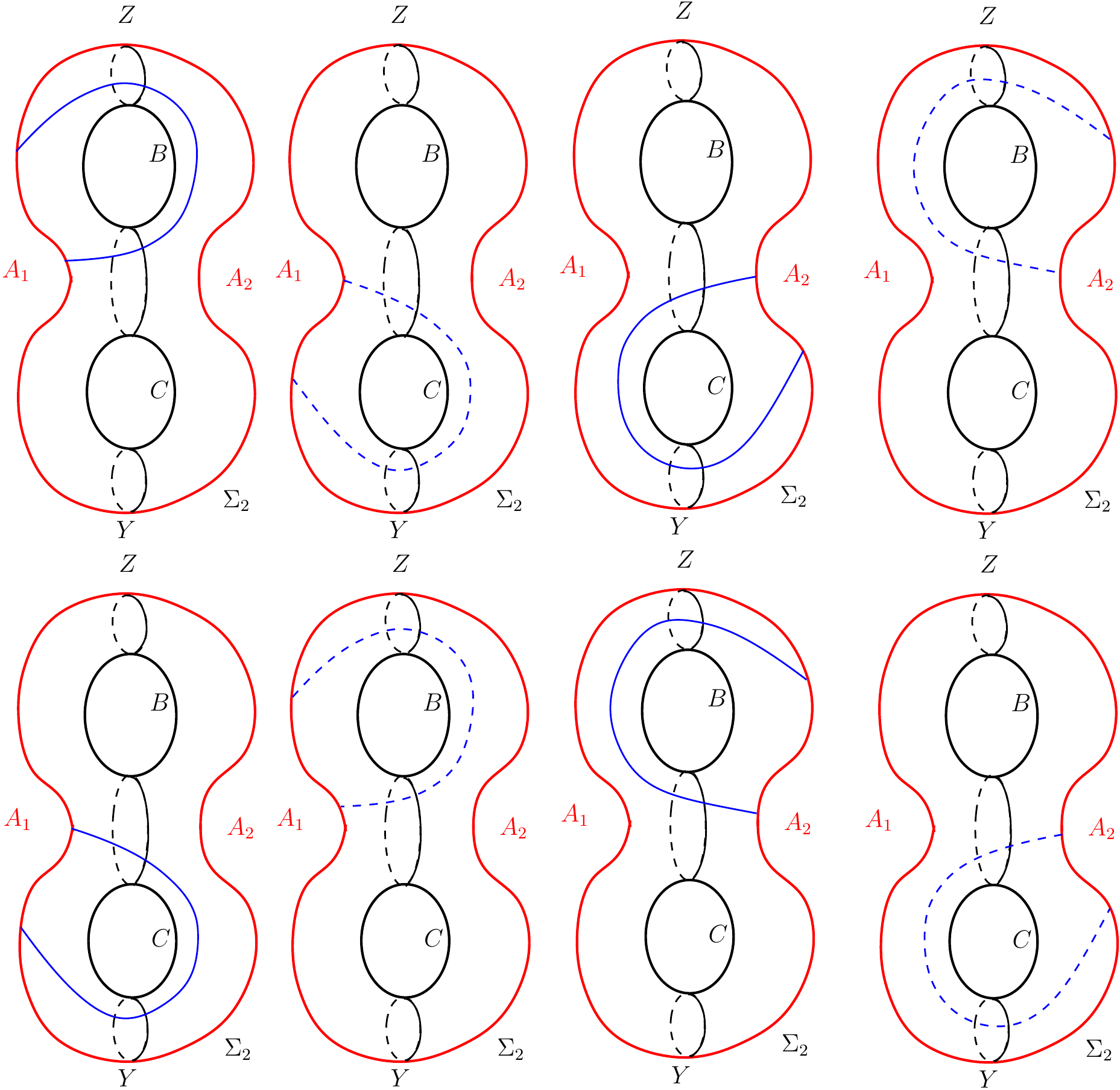}
\caption{All eight possibilities for a single word $AA$ arc}\label{fig:singleAA}
\end{figure}

So the above discussion proves Theorem \ref{thm:beta_single_letter}.

\section{A Reduction Algorithm}\label{section:algorithm}

We have the following:

\begin{lem}\label{lem:reduce_a_Q}
Suppose that $a_Q > b_Q + c_Q$ for a reducing sphere $Q$ of $\mH$. Then $b_Q \neq c_Q$, unless $b_Q = c_Q = 0$ in which case $a_Q$ is automatically $2$.
\end{lem}
\begin{proof}
First we note that all three of $a_Q, b_Q$ and $c_Q$ cannot be simultaneously zero because then $\qsig$ will be a trivial curve on $\partial V$ implying that $Q$ is not a reducing sphere.

If $b_Q = c_Q = 0$, we can apply $\beta^k$ for some integer $k$ to $\mH$ to reduce $a_Q$ to $2$ at which stage we get, $P$, the standard reducing sphere. But then $\beta^{-k}(P) = Q$ and $P$ is invariant under $\beta$ implying $P = Q$, which means $a_Q = 2$.

If possible, suppose that $b_Q = c_Q \neq 0$. Since, $a_Q > b_Q + c_Q$, we can use Theorem \ref{thm:beta_single_letter} to lower $a_Q$. Suppose that $k$ is an integer of smallest absolute value such that $\beta^k(Q) = R$ and $a_R < b_R + c_R$. Since $\beta$ does not change $b_Q$ and $c_Q$, $b_R = c_R = b_Q = c_Q \neq 0$. So, neither $b_R > a_R + c_R$ nor $c_R > b_R + a_R$, contradicting Theorem \ref{arc-lemma} as applied to $R$. So, we cannot have $b_Q = c_Q \neq 0$. 
\end{proof}

Among the four automorphisms of $\mH$ mentioned in section \ref{section:red_curves}, we note that the non-zero integral powers of $\beta$ are the only ones which change the value of $a_Q + b_Q + c_Q$. Lemma \ref{lem:reduce_a_Q} shows that the lowest possible value of $a_Q + b_Q + c_Q$ is $2$. Three possibilities for $Q$ with $a_Q + b_Q + c_Q = 2$ are:  (i) $a_Q = 2, b_Q = c_Q = 0$ giving the standard reducing sphere $P$,  (ii) $b_Q = 2, a_Q = c_Q = 0$ giving the simple reducing sphere $P'$ or (iii) $c_Q = 2, b_Q = a_Q = 0$ giving the simple reducing sphere $P''$. The reducing curves of $P, P'$ and $P''$ are shown in Figure \ref{fig:std-red-curves}. We note that $\delta^{-1}(P') = P$ and $\delta(P'') = P$, where $\delta$ is the order three automorphism of $\mH$ described in section \ref{section:red_curves}.

\begin{figure}[h!]
	\begin{subfigure}{.32\linewidth}
		\centering
		\includegraphics[width=\linewidth]{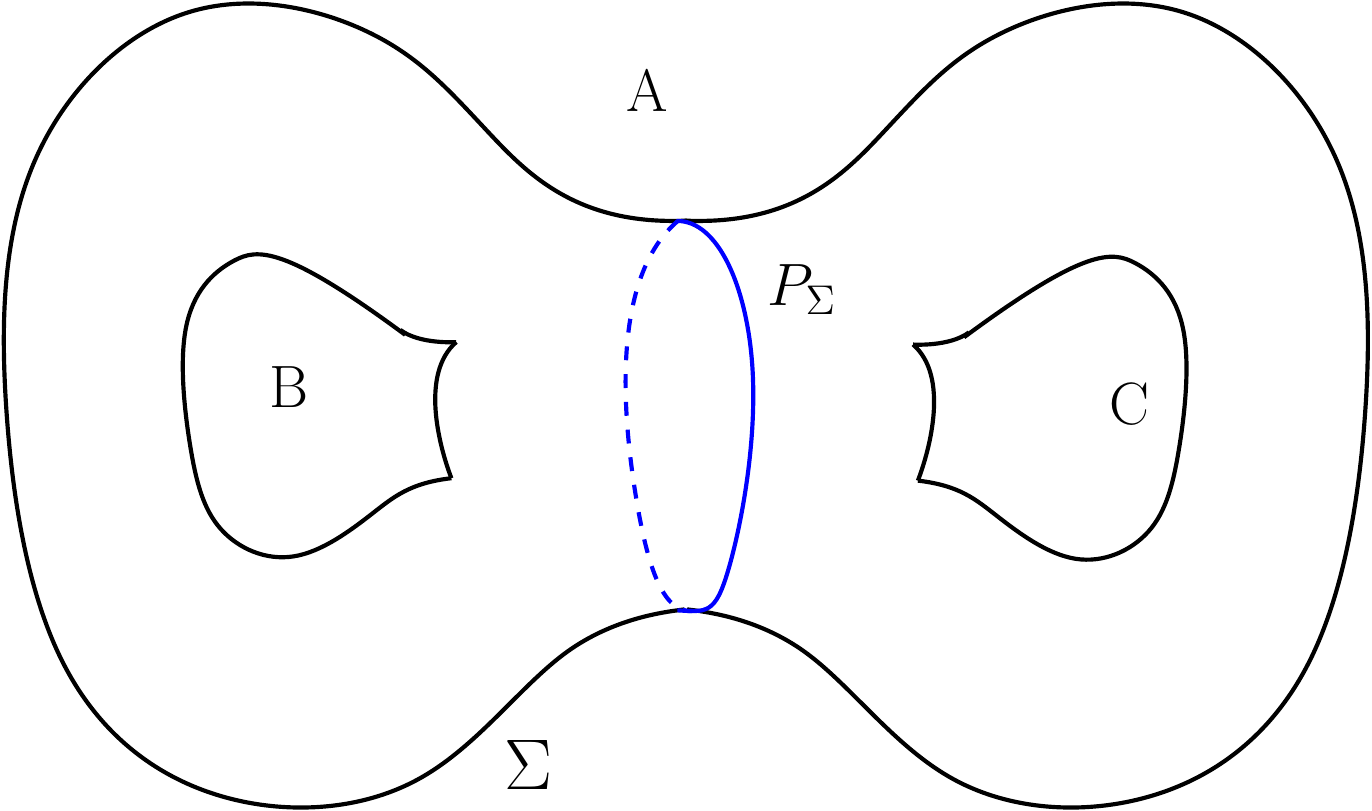}
		\caption{Reducing curve of $P$}
	\end{subfigure}
	\begin{subfigure}{.32\linewidth}
		\centering
		\includegraphics[width=\linewidth]{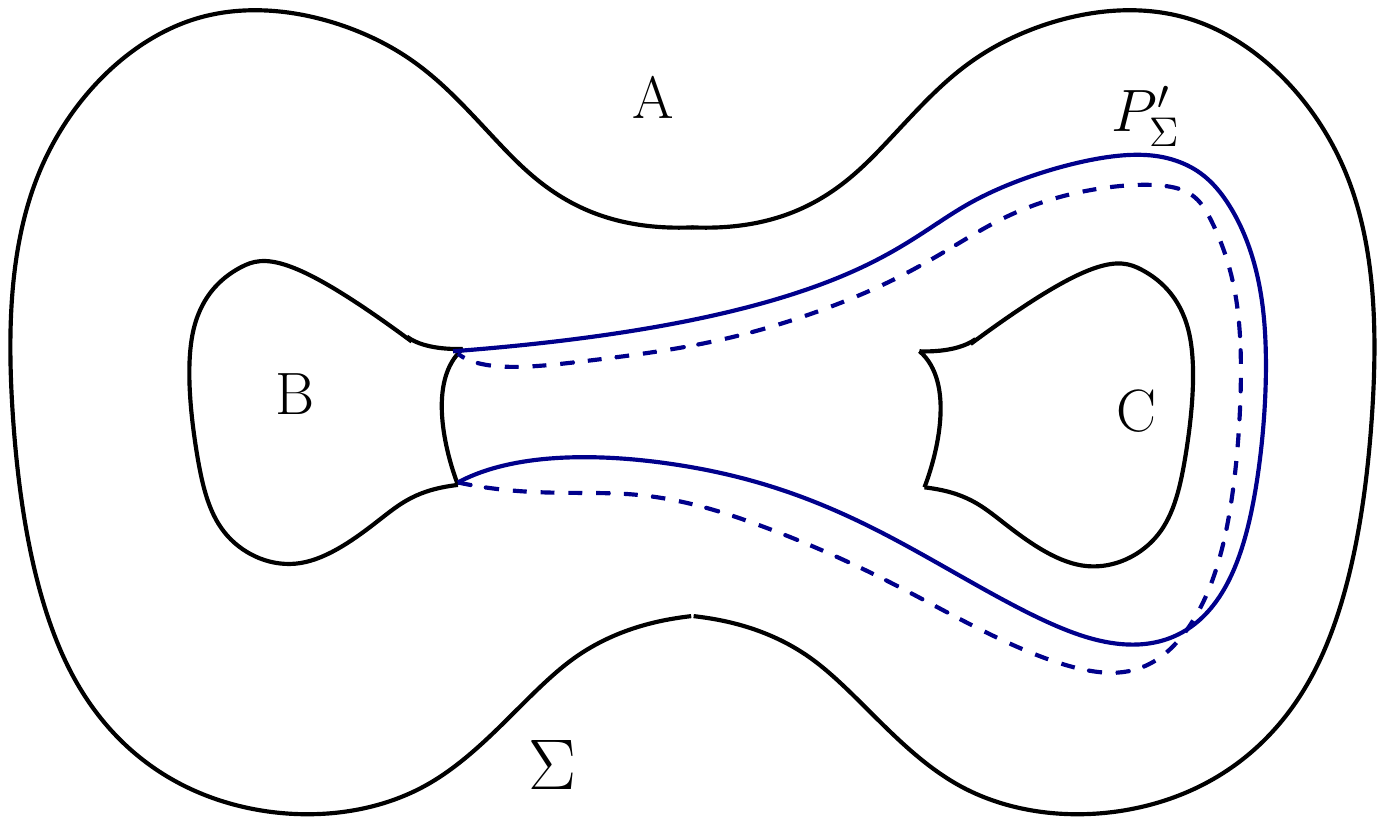}
		\caption{Reducing curve of $P'$}
	\end{subfigure}
	\begin{subfigure}{.32\linewidth}
		\centering
		\includegraphics[width=\linewidth]{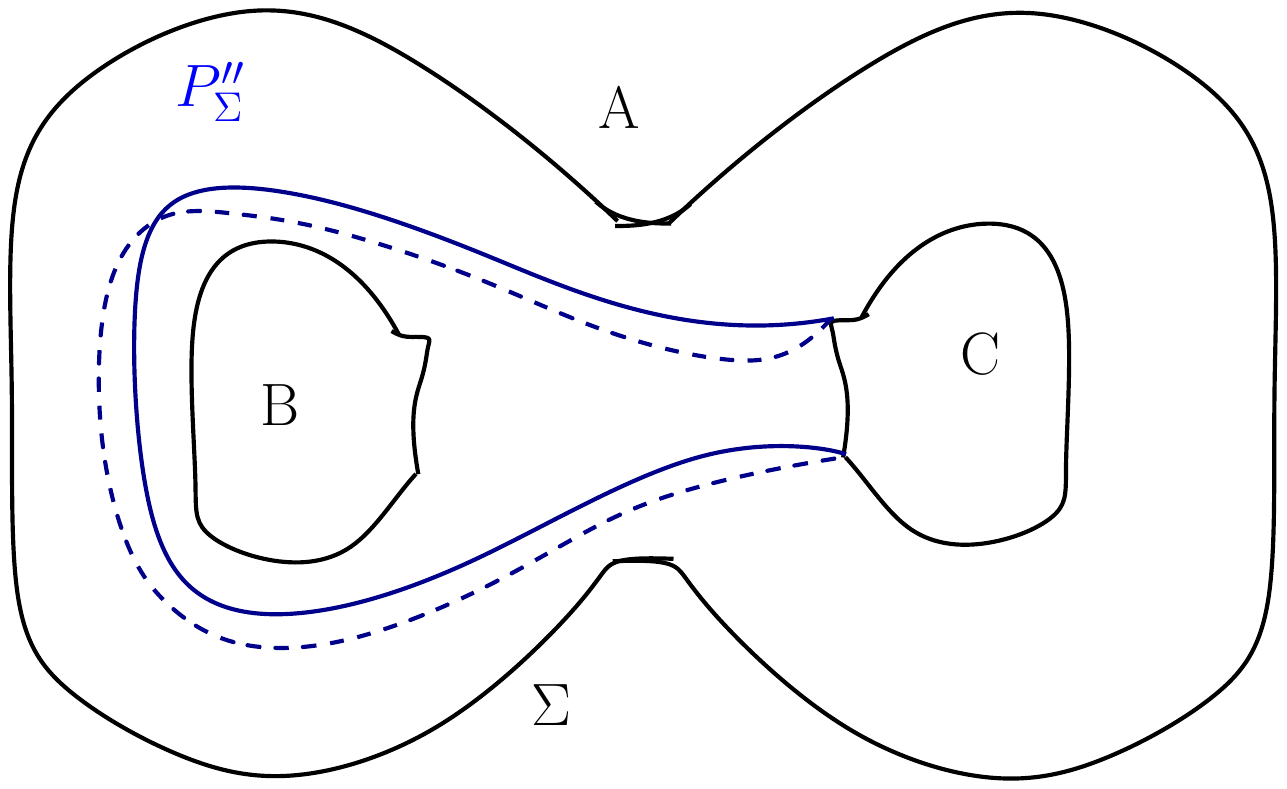}
		\caption{Reducing curve of $P''$}
	\end{subfigure}
	\caption{The reducing curves $P_\Sigma, P_\Sigma'$, and $P_\Sigma''$}\label{fig:std-red-curves}
\end{figure}

Let $Q$ be a reducing sphere of $\mH$. Using automorphisms $\beta$ and $\delta$, we have the following algorithm which transforms $Q$ into $P$:

\begin{algorithm} Reduction Algorithm: \\
\begin{enumerate}
\item[Step 1:] If $a_Q + b_Q + c_Q = 2$ then $Q = P, P'$ or $P''$.\\ 
						If $Q = P$, Exit. \\
						Else: apply $\delta$ or $\delta^{-1}$ to $P''$ or $P'$, respectively, to get $P$. Exit.
						
\item[Step 2:] While $a_Q + b_Q + c_Q > 2$, Do
			\begin{itemize}
					\item[Step (i):] If $a_Q > b_Q + c_Q$, go to Step (ii), else: let $\delta(Q) = R'$, $\delta^{-1}(Q) = R''$. If $a_{R'} > b_{R'} + c_{R'}$, set $Q$ to be $R'$. If $a_{R''} > b_{R''} + c_{R''}$, set $Q$ to be $R''$. Go to Step (ii).
					\item[Step (ii):] Let $\beta(Q) = S'$, $\beta^{-1}(Q) = S''$. If $a_{S'} < a_Q$, set $Q$ to be $S'$. If $a_{S''} < a_Q$, set $Q$ to be $S''$. 
			\end{itemize}
\item[Step 3:] Now, $a_Q + b_Q + c_Q = 2$ and $Q = P, P'$ or $P''$.\\
					 If $Q = P$, Exit. \\
					 Else: apply $\delta$ or $\delta^{-1}$ to $P''$ or $P'$, respectively, to get $P$. Exit.	
\end{enumerate}
\end{algorithm}

Note that in Step (i) within the while-loop in the reduction algorithm for $Q$, if $a_Q < b_Q + c_Q$, Theorem \ref{arc-lemma} asserts that one of $a_{R'} > b_{R'} + c_{R'}$ or $a_{R''} > b_{R''} + c_{R''}$ will hold. Also, in Step (ii) within the while-loop, Theorem \ref{thm:beta_single_letter} asserts that one of $a_{S'} < a_Q$ or $a_{S''} < a_Q$ holds. So each iteration of the while loop indeed updates $Q$ to a new reducing sphere $Q'$ such that $a_{Q'} + b_{Q'} + c_{Q'} < a_Q + b_Q + c_Q$.

By applying the reduction algorithm to the given reducing sphere $Q$, we can obtain $P$. So by applying the automorphisms $\beta, \beta^{-1}, \delta$ or $\delta^{-1}$ in the inverse sequence, as dictated by the reduction algorithm, we can obtain $Q$ from $P$. Note that this inverse process applies a reduced word, $w$, of automorphisms $\beta, \beta^{-1}, \delta$ or $\delta^{-1}$ to $P$ to obtain $Q$. 

Let $h \in G_2$, then $h(P) = Q$ for some reducing sphere $Q$. By using the reducing algorithm, we get a word $w$ in the automorphisms $\beta, \beta^{-1}, \delta$ or $\delta^{-1}$ such that $w(Q) = P$, so $w\circ h(P) = w(Q) = P$. This implies that $w \circ h \in Stab(P)$. If we assume that $Stab(P)$ is finitely generated by the automorphisms $\alpha, \beta$ and $\gamma$, then since $w$ is also a finite word generated by $\beta$ and $\delta$, we have that $h$, an arbitrary element in $G_2$, can also be  written as a word in $\alpha, \beta, \gamma, \delta$ and their inverses. Hence we can conclude that $G_2$ is finitely generated by $\{\alpha, \beta, \gamma, \delta\}$.  

We conclude this article by asking a question. Let a triple $(x,y,z)$ of non-negative integers be called a \textit{non-triangular triple} if it satisfies one of the inequalities $x>y+z, y>z+x$ or $z>x+y$. From Corollary \ref{cor:cor_to_arc_lemma}, the set of triples $(a_Q, b_Q, c_Q)$ of intersection numbers with the curves $A, B, C$ as in Figure \ref{fig:standard_curves} of a given reducing sphere $Q$ for $\mH$ is a specific subset $\Lambda$ of the set of non-triangular triples of integers.

\begin{question} Describe the set $\Lambda$ which is in one-to-one correspondence with the set of reducing spheres for the genus-2 Heegaard splitting $\mH$ of $S^3$.
\end{question} 

\bibliographystyle{plainnat}
\bibliography{bib}	

\begin{thebibliography}{10}
\providecommand{\natexlab}[1]{#1}
\providecommand{\url}[1]{\texttt{#1}}
\expandafter\ifx\csname urlstyle\endcsname\relax
  \providecommand{\doi}[1]{doi: #1}\else
  \providecommand{\doi}{doi: \begingroup \urlstyle{rm}\Url}\fi

\bibitem[Akbas(2008)]{akbas}
Erol Akbas.
\newblock A presentation for the automorphisms of the 3-sphere that preserve a
  genus two {H}eegaard splitting.
\newblock \emph{Pacific J. Math.}, 236\penalty0 (2):\penalty0 201--222, 2008.

\bibitem[Cho(2008)]{cho}
Sangbum Cho.
\newblock Homeomorphisms of the 3-sphere that preserve a {H}eegaard splitting
  of genus two.
\newblock \emph{Proc. Amer. Math. Soc.}, 136\penalty0 (3):\penalty0 1113--1123,
  2008.

\bibitem[Farb and Margalit(2012)]{farb_margalit}
Benson Farb and Dan Margalit.
\newblock \emph{A primer on mapping class groups}, volume~49 of \emph{Princeton
  Mathematical Series}.
\newblock Princeton University Press, Princeton, NJ, 2012.

\bibitem[Freedman and Scharlemann(2018)]{freedman}
Michael Freedman and Martin Scharlemann.
\newblock Powell moves and the goeritz group, 2018.
\newblock URL \url{https://arxiv.org/abs/1804.05909}.

\bibitem[Powell(1980)]{powell}
Jerome Powell.
\newblock Homeomorphisms of {$S\sp{3}$} leaving a {H}eegaard surface invariant.
\newblock \emph{Trans. Amer. Math. Soc.}, 257\penalty0 (1):\penalty0 193--216,
  1980.

\bibitem[Scharlemann(2003)]{scharlemann_notes}
Martin Scharlemann.
\newblock Heegaard splittings of 3-manifolds.
\newblock In \emph{Low dimensional topology}, volume~3 of \emph{New Stud. Adv.
  Math.}, pages 25--39. Int. Press, Somerville, MA, 2003.

\bibitem[Scharlemann(2004)]{scharlemann}
Martin Scharlemann.
\newblock Automorphisms of the 3-sphere that preserve a genus two {H}eegaard
  splitting.
\newblock \emph{Bol. Soc. Mat. Mexicana (3)}, 10\penalty0 (Special
  Issue):\penalty0 503--514, 2004.

\bibitem[Volodin et~al.(1974)Volodin, Kuznecov, and Fomenko]{volodin}
I.~A. Volodin, V.~E. Kuznecov, and A.~T. Fomenko.
\newblock The problem of the algorithmic discrimination of the standard
  three-dimensional sphere.
\newblock \emph{Uspehi Mat. Nauk}, 29\penalty0 (5(179)):\penalty0 71--168,
  1974.
\newblock Appendix by S. P. Novikov.

\bibitem[Waldhausen(1968)]{waldhausen}
Friedhelm Waldhausen.
\newblock Heegaard-{Z}erlegungen der {$3$}-{S}ph\"{a}re.
\newblock \emph{Topology}, 7:\penalty0 195--203, 1968.

\bibitem[Zupan(2020)]{zupan}
Alexander Zupan.
\newblock The {P}owell conjecture and reducing sphere complexes.
\newblock \emph{J. Lond. Math. Soc. (2)}, 101\penalty0 (1):\penalty0 328--348,
  2020.

\end{thebibliography}
\end{document}